%% file: main.tex
\documentclass{article}
\usepackage{arxiv}

\usepackage{generic}
\usepackage{cite}
\usepackage{amsmath,amssymb,amsfonts}
\usepackage{lipsum}
\usepackage{amsfonts}
\usepackage{graphicx}
\usepackage{epstopdf}
\usepackage{algorithmic}
\usepackage{textcomp}
\usepackage{comment} 
\usepackage{lipsum} 
\usepackage{graphicx}
\usepackage{fullpage} 
\usepackage{amssymb, amsfonts, amsmath, }
\usepackage{amsmath,amsthm}
\usepackage{bbm}
\usepackage{epsfig, latexsym, graphicx}
\usepackage{soul}
\usepackage{tikz} 
\usepackage{todonotes}
\usepackage{algorithm}
\usepackage{comment}
\usepackage{tabularx}
\usepackage[american]{circuitikz}

\usepackage{graphicx} 
\usepackage{subcaption}
\usepackage{algorithmic} 
\usepackage{enumitem}
\usepackage{booktabs}
\usepackage{url}
\usetikzlibrary{shapes,arrows}

\DeclareMathOperator*{\argmin}{arg\,min}
\newcommand{\R}{\mathbb{R}}
\newcommand{\x}{\mathbf{x}}
\newcommand{\z}{\mathbf{z}}
\newcommand{\ahat}{\hat{\mathbf{a}}}
\newcommand{\fp}{\hat\x} 
\newcommand{\y}{\mathbf{y}}

\newcommand{\fh}{\hat{f}}

\newcommand{\dfdx}{\nabla f(\x)}
\newcommand{\dfdy}{\nabla f(\y)}
\newcommand{\dfdxt}{\nabla f(\xt)}
\newcommand{\dxdt}{\dot{\x}}
\newcommand{\xt}{\x(t)}
\newcommand{\Zeta}{Z}
\newcommand{\Ab}{\mathbf{A}}

\usepackage{accents}

\usepackage[export]{adjustbox}

\theoremstyle{plain}
\newtheorem{theorem}{Theorem}[section]

\newtheorem{lemma}[theorem]{Lemma}

\theoremstyle{definition}
\newtheorem{definition}[theorem]{Definition}

\theoremstyle{remark}

\begin{document}

\title{An Equivalent Circuit Workflow for Unconstrained Optimization}

\author{Aayushya Agarwal$^{1*}$, Carmel Fiscko$^{1*}$, Soummya Kar$^{1}$, Larry Pileggi$^{1}$, and Bruno Sinopoli$^{2}$
\thanks{\hspace{-2em}$*$These authors contributed equally.}%
\thanks{\hspace{-2em}$^{1}$Aayushya Agarwal, Carmel Fiscko, Soummya Kar, and Larry Pileggi are with the Dept. of Electrical and Computer Engineering at Carnegie Mellon University at 5000 Forbes Ave, Pittsburgh, PA 15213. {\tt\small \{aayushya,cfiscko,soummyak,pileggi\}@andrew.cmu.edu }}%
\thanks{\hspace{-2em}$^{2}$Bruno Sinopoli is with the Dept. of Electrical and Systems Engineering at Washington University in St. Louis, MO at 1 Brookings Dr, St. Louis, MO 63130. {\tt\small bsinopoli@wustl.edu }}%
}
\maketitle

\begin{abstract}
 We introduce a new workflow for unconstrained optimization whereby objective functions are mapped onto a physical domain to more easily design algorithms that are robust to hyperparameters and achieve fast convergence rates. Specifically, we represent optimization problems as an equivalent circuit that are then solved solely as nonlinear circuits using robust solution methods. The equivalent circuit models the trajectory of component-wise scaled gradient flow problem as the transient response of the circuit for which the steady-state coincides with a critical point of the objective function.
 The equivalent circuit model leverages circuit domain knowledge to methodically design new optimization algorithms that would likely not be developed without a physical model.
 We incorporate circuit knowledge into optimization methods by 1) enhancing the underlying circuit model for fast numerical analysis, 2) controlling the optimization trajectory by designing the nonlinear circuit components, and 3) solving for step sizes using well-known methods from the circuit simulation. We first establish the necessary conditions that the controls must fulfill for convergence. We show that existing descent algorithms can be re-derived as special cases of this approach and derive new optimization algorithms that are developed with insights from a circuit-based model. The new algorithms can be designed to be robust to hyperparameters, achieve convergence rates comparable or faster than state of the art methods, and are applicable to optimizing a variety of both convex and nonconvex problems. 
\end{abstract}



\input{Intro}

\input{Previous_works}
\input{Problem_Formulation}

\input{Meta-Alg.tex}

\input{equivalent_circuit_model.tex}
\input{modifying_ec.tex}

\input{Control}

\input{Discretization.tex}

\input{examples.tex}

\input{Algo.tex}

\input{results.tex}

\input{Conclusion}
\newpage
\appendix

\input{Appendix.tex}

\bibliographystyle{siamplain}
\bibliography{main}
\end{document}

%% file: Intro.tex
\section{Introduction}

Optimization is a key problem across all areas of science and engineering. Originally, optimization methods were inspired by a physical analogy of a ball rolling down a hill. This physical analogy provided additional domain knowledge that inspired heuristics and extensions, like acceleration.

\begin{figure}
\centering
\subfloat[Workflow for standard optimization methods.]{\label{fig:optimization_workflow}\includegraphics[width=0.4\linewidth]{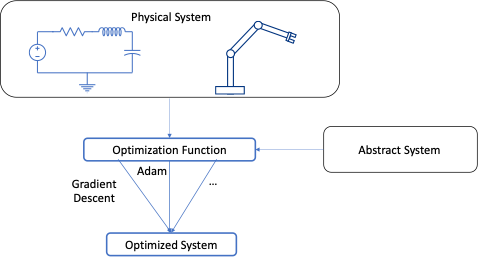}}\label{mult_ec}
\subfloat[ECCO workflow uses an equivalent circuit model.]{\label{fig:ecco_workflow}\includegraphics[width=0.4\linewidth]{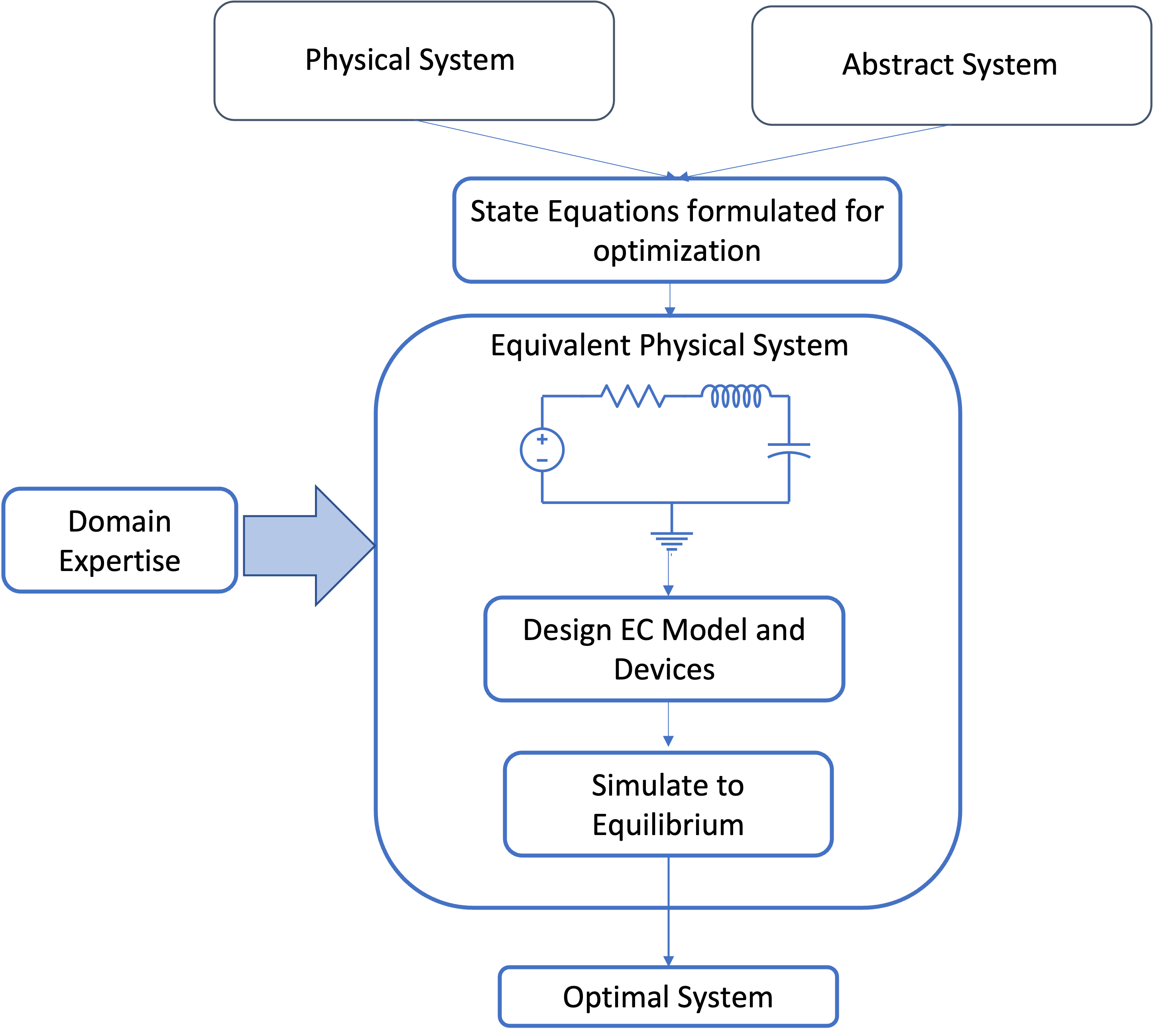}}\label{mult_trans}
\label{fig:testfig}
\end{figure}

Optimization in practice, however, prefers to abstract out the relevant equations, thus severing the optimization technique from intuition that the application may provide. For example, consider the chart shown in Figure \ref{fig:optimization_workflow}.  Physical systems such as a circuits, and abstract systems, such as neural networks, are first mathematically represented in the form of a general optimization problem. The user then selects an optimization technique such as gradient descent to optimize the mathematical representation to reach an optimized system. Methods such as Adam \cite{kingma2014adam} and extra-gradient descent \cite{korpelevich1976extragradient} are generally not tied to any physical counterpart and instead rely on the experience of researchers to select hyperparameters. In abstracting physical systems to a set of equations, this workflow loses insights from the original physical domain that can easily help design step sizes and shape trajectories to better optimize the system.


We present a new method that instead solves optimization problems by modeling them as physical processes. This enables derivation of new optimization algorithms by solving the optimization problems within the domain of the chosen physics application. In this workflow, shown in Figure \ref{fig:ecco_workflow}, we first convert the optimization problem to a continuous-time gradient flow method, where the trajectory of optimization variables is a set of state-space equations whose steady state coincides with a critical point. As many physical processes are naturally described using a state space formulation, gradient flow provides an opportunity to map optimization methods onto physical models, thereby unifying the ideas in physics with general optimization. In our work, the state space formulation is mapped onto an equivalent circuit (EC) whose transient response models the trajectory of the optimization variables and the steady-state coincides with a critical point. We select electrical circuits due to the large breadth of knowledge on circuit analysis and simulation. The optimization problem is then solved solely in the domain of circuits by designing new circuit-inspired optimization methods that capitalize on the additional domain knowledge. 

To solve the optimization problem in the circuit domain, we design the EC model and the nonlinear devices and then simulate the circuit to its steady state (i.e., the critical point of the optimization problem), as shown in Figure \ref{fig:ecco_workflow}. In the process, we integrate concepts of energy and charge, along with an existing wealth of knowledge on circuits to design the circuit to reach steady-state faster and devise simulation methods that complement the EC.

Our approach introduces a new structured framework, which we introduce as Equivalent Circuit Controlled Optimization (ECCO), that uses circuit knowledge to design new algorithms that control the trajectory and step-sizes of the optimization method. ECCO provides three opportunities for practitioners to incorporate circuit domain knowledge into optimization methods: 1) devise the underlying EC model to reach steady state fast, 2) design nonlinear capacitors that control the trajectory of the optimization, and 3) devise new discrete time-stepping simulation algorithms to select step sizes based on circuit principles. 

For the optimization practitioner or researcher, this circuit based framework has several key advantages. First, demonstrating that gradient flows are mathematically equivalent to circuits enables the vast circuits knowledge to be used to solve optimization problems, both analytically and in practice. 
Decades of research in circuit simulation have produced sophisticated simulation techniques that combine circuit physics, numerical tricks, and domain knowledge into industry-standard tools to efficiently simulate systems of millions of variables. By representing optimization problems as circuits, the optimization community may leverage additional domain-insights to easily derive new optimization methods. Second, this work is a \emph{unifying workflow} that provides a structure to construct new optimization trajectories in continuous-time that are then solved using numerical integration methods. In this work, we establish parallels between optimization iterations and numerical integration methods, which provides a new powerful approach that utilizes vast research in numerical methods to solve optimization problems. Unlike previous works that solve gradient-flow \cite{behrman1998efficient,franca2018dynamical,antipin1994minimization}, we provide an end-to-end design of the continuous-time trajectory as well as numerical methods to efficiently solve to a steady-state. 
We demonstrate how standard existing algorithms can be derived via different choices of controls or numerical integration methods, and thus how the circuits can be used to propose new algorithms for different purposes. Third, this framework is ideal for both generalization and refinement; for example, future work could draw specific circuits could for structured objective functions, e.g. loss functions with regularizers, to further propose new algorithms.

\textbf{Our main contributions in this paper} are (1) developing the EC model for component-wise scaled gradient flow and proposing circuit-intuitive modifications, (2) deriving a second order and a first order control scheme for the gradient flow based on the EC, and (3) proposing an adaptive numerical integration algorithm for optimization that choose step sizes by bounding the estimated local truncation error (LTE) accrued and guaranteeing numerical stability. We find that the ECCO-derived algorithms can be designed to achieve fast convergence rates to critical points and be highly robust to hyperparameters. To the best of our knowledge, this is the first work formally using circuit theory for gradient flow, and using circuit simulation techniques to solve general optimization problems. No circuit background knowledge is needed to implement the proposed algorithms.

In Section \ref{sec:problem_formulation}, we specify the problem formulation of the scaled gradient flow and introduce necessary convergence conditions. In Section \ref{meta-alg} we introduce the ECCO workflow at a high level. The EC is derived in Section \ref{sec:ec_model}. Controlling the EC is discussed in Section \ref{sec:control}, as well as derivation of two new control schemes. Discretization of the system, i.e. numerical solution methods, is discussed in Section \ref{sec:discretization}. Examples applying these ideas to rederive existing optimization methods are shown in Section \ref{sec:examples}. Our new proposed algorithms are summarized in Section \ref{sec:algo}, and they are demonstrated in simulation in Section \ref{sec:results}.

%% file: Previous_works.tex
\section{Related Work}
\label{sec:previous_works}

Discrete, iterative optimization algorithms such as gradient descent and its variants are well-loved techniques for their ease of implementation \cite{boyd2004convex}. While convergence properties of these algorithms can be shown, deeper investigations consider continuous-time optimization trajectories where ODE theory have been used to analyze the optimization variable and design new new optimization methods.

Continuous-time, dynamical system representations of first-order optimization methods, often called gradient flow methods, have long been studied \cite{behrman1998efficient}, \cite{antipin1994minimization}, \cite{attouch1996dynamical}, \cite{brown1989some}, including for neural networks \cite{baldi1995gradient}. These works apply control mechanisms in the aim of analyzing and improving on gradient descent methods. A recent surge of interest in ODE realizations of descent methods have given strong convergence analyses of methods such as distributed techniques \cite{swenson2021distributed}, momentum methods \cite{muehlebach2021optimization} \cite{franca2018dynamical}, stochastic gradient descent \cite{latz2021analysis}, as well as analyzing ADMM \cite{franca2018admm}. These works provide insight into choosing appropriate parameters such as momentum values in Nestrov's acceleration \cite{su2014differential}. 
 Strong connections between gradient descent and linear control systems have been established, showing that these methods can accomplish the same goals \cite{hardt2016gradient}. First-order optimization methods have been shown to be interpretable as elements such as PID controllers and lag compensators \cite{hu2017control}. These ideas have naturally extended to the distributed optimization paradigm \cite{wang2010control}, \cite{sundararajan2019canonical}. 

One important method in control systems is Lyapunov theory, widely used to show stability of dynamical systems, which has emerged as a tool to show convergence of continuous-time optimization techniques \cite{wilson2018lyapunov}, \cite{wilson2021lyapunov} and acceleration methods\cite{polyak2017lyapunov}, \cite{hustig2019robust}. Fewer works have used Lyapunov theory to propose improved optimization algorithms \cite{taylor2018lyapunov}. A key issue in applying Lyapunov theory to optimization is that traditional techniques in control systems use domain knowledge of the system - for example, the energy of the system - yet standard optimization tools work with just an objective function, perhaps with some known properties. For optimization to take full advantage of Lyapunov theory, some physical interpretation is desired. One work \cite{hu2017dissipativity} created an analogy with energy dissipation in physics, yielding insights into Lyapunov function construction and convergence analysis. We are only aware of one other reference for using circuits to model an optimization problem \cite{boyd2021distributed}, but that work focuses on a distributed optimization setup and defining a proximal operator. Nevertheless, there exists a gap in using circuit intuition for continuous-time optimization. 

While useful theoretical results can be established with the continuous-time system, solving the ODE with a computer generally necessitates some numerical integration scheme.  Recent advances have used more sophisticated explicit and implicit integration techniques to simulate the continuous system \cite{pmlr-v97-muehlebach19a}, \cite{wadiaoptimization}, \cite{lin2016distributed}. 
Many engineering simulation tools, such as SPICE \cite{nagel1971computer}, often exploit domain-specific behavior of physical systems to solve the underlying ODEs much faster. SPICE is a circuit simulation engine, developed in 1971, that analyzes the transient response of circuits by solving an underlying stiff, nonlinear ODE. Unlike generic ODE solvers, SPICE uses physics and circuits knowledge to develop intuitive heuristics that solve the underlying ODE and scale upwards to billions of variables. Key insights of passivity and quasi-steady-state behavior of energy-storage devices have lead to insights for numerical integration algorithms.  Additionally, intuitive modifications of the circuit have been proposed, in the form of homotopy methods such as pseudo-transient analysis \cite{Nilsson2000-jy}, Gmin stepping \cite{Nilsson2000-jy} and source stepping \cite{Nilsson2000-jy} to easily simulate the circuit behavior.

Our equivalent circuit framework closely follows circuit simulation practices that develops heuristics based on insights from circuits. This enables us to not only utilize various control methodologies as previous works have done, but also utilize circuit simulation methods which have demonstrated robust and scalable convergence.

Specifically, we utilize knowledge of circuit models to design nonlinear devices in the EC that ensures proper energy dissipation to reach a steady-state. In addition, we adopt dynamic time-stepping algorithms from SPICE to numerically solve for the steady-state. Finally, we realize an intuitive homotopy method that is applied to the equivalent circuit model which experimentally demonstrates improved convergence. 


%% file: Problem_Formulation.tex
\section{Problem Formulation}
\label{sec:problem_formulation}
In this work, we consider the following unconstrained optimization problem:
\begin{gather}
    \min_{\x} f(\x),\label{prob}\\
    \x^* \in \argmin_{\x} f(\x).
\end{gather}
where $\x\in \R^n$ and $f: \R^n\to\R$. It is known \cite{brown1989some} that  \eqref{prob} may be solved via the gradient flow ODE initial value problem (IVP),
\begin{equation}
    \dxdt(t)=-\dfdxt,\quad \x(0)=\x_0, \label{gradient flow}
\end{equation}
where $\dxdt(t)$ refers to the time derivative of $\xt$. In this work we consider scaled gradient flow IVP, where the RHS of \eqref{gradient flow} is multiplied by some invertible matrix $Z^{-1}$.
\begin{equation}
    \dxdt(t)=-Z(\xt)^{-1} \nabla f(\xt),\quad \x(0) = \x_0. \label{scaled gradient flow}
\end{equation}
Note that in the future sections with circuit notation we will move $Z$ to the LHS for convenience and the ODE will be of the form  $Z(\xt)\dxdt(t) = -\dfdxt$.

The solution to an ODE IVP is computed with the integral,
\begin{equation}
    \x(t) = \x(0) + \int_0^t \dxdt(\tau) d\tau.\label{soln}
\end{equation}

One of the objectives in the ECCO workflow is to design $Z$ such that $\x(t)\to \x^*$ as $t\to\infty$. If $Z\equiv \mathbf{I}_{n}$, i.e. an identity matrix of size $n\times n$, then we call the ODE \emph{uncontrolled gradient flow}. 

\subsection{Identities}
The following identities will be used within this paper.

\begin{lemma}\label{identity1}
The following holds:
\begin{equation}
    -\frac{d}{dt}\nabla f(\xt) = \nabla^2 f(\xt)Z(\xt)^{-1}\nabla f(\xt).
\end{equation}
\end{lemma}
\begin{proof}
    \begin{align*}
        -\frac{d}{dt}\nabla f(\xt)&=-\frac{d\dfdxt}{d\xt} \frac{d\xt}{dt}\\
        &=-\nabla^2 f(\xt)(-Z(\xt)^{-1}\dfdxt)\\
        &=\nabla^2 f(\xt) Z(\xt)^{-1}\dfdx
    \end{align*}
\end{proof}

\begin{lemma}\label{identity2}
    The following holds:
    \begin{equation}
        -\frac{1}{2}\frac{d}{dt}\|\dfdxt\|^2 = \dfdxt^{\top}\nabla^2 f(\xt) Z(\xt)^{-1}\dfdxt. \label{ddt identity}
    \end{equation}
\end{lemma}
\begin{proof}
As a consequence of Lemma \ref{identity1} we find:
    \begin{align*}
        -\frac{1}{2}\frac{d}{dt}\|\dfdxt\|^2= -\dfdxt^{\top}\frac{d}{dt}\Big[\dfdxt\Big]= \dfdxt^{\top}\nabla^2 f(\xt) Z(\xt)^{-1}\dfdxt.
    \end{align*}

\end{proof}

\subsection{Newton Method Gradient Flow}
One strategy to design $Z(\xt)^{-1}$ is to choose the control that maximizes the decrease in gradient over time. In Section \ref{sec:control}, we will show that this choice is well-motivated via the circuit formulation. Thus,
\begin{align}
    &\max_{Z} -\frac{1}{2}\frac{d}{dt}\|\dfdxt\|^2,\label{dense z setup 1}\\
    &=\max_{Z} \dfdxt^{\top}\nabla^2 f(\xt) Z(\xt)^{-1}\dfdxt,\label{dense z setup 2}\\
    &=\max_{Z}\dfdxt^{\top} Q(\xt)\Lambda(\xt)W(\xt)Q(\xt)^{\top}\dfdxt, \label{dense z setup 3}\\
    &=\max_{\{w_{ii}\}}\sum_{i=1}^n \lambda_i(\xt) w_{ii}(\xt) \nabla f(\xt)^{\top}q_i(\xt)q_i(\xt)^{\top}\dfdxt,\label{dense z setup 4}
\end{align}
where \eqref{dense z setup 2} uses the identity in Lemma \ref{identity2} and \eqref{dense z setup 3} uses the eigenvalue decomposition $\nabla^2 f(\xt) \equiv Q(\xt)\Lambda(\xt)Q(\xt)^{\top}$ and similarity transformation $Z(\xt)^{-1}\equiv Q(\xt)W(\xt) Q(\xt)^{\top}$ where $W$ is a diagonal matrix with diagonal entries $w_{ii}$. To make \eqref{dense z setup 4} tractable, we will assert that each $w_{ii}$ must be nonzero and the set $\{w_{ii}\}$ has a finite sum.

The objective now is to design the set of $w_{ii}$ to solve \eqref{dense z setup 4}. Note that for uncontrolled gradient flow it holds that,
\begin{align}
    \lambda_{min}\|\nabla f(\xt)\|^2\leq \dfdxt^{\top}\nabla^2 f(\xt) \dfdxt \leq \lambda_{max}\|\nabla f(\xt)\|^2,
\end{align}
where $\lambda_{min}$ and $\lambda_{max}$ are the minimum and maximum eigenvalues of the Hessian. The controlled gradient flow likewise satisfies,
\begin{align}
    \min_i(\lambda_{i}w_{ii})\|\nabla f(\xt)\|^2\leq \dfdxt^{\top}\nabla^2 f(\xt) Z(\xt)^{-1}\dfdxt \leq \max_i(\lambda_{i}w_{ii})\|\nabla f(\xt)\|^2.
\end{align}
One strategy may be to choose $W$ to increase this upper bound to potentially enable better performance; however, this can lead to a logical fallacy. Consider a control where $w_{ii}$ for $i$ corresponding to $\lambda_{max}$ is heavily weighted. The maximum eigenvalue of the controlled system will be increased, thus successfully increasing the upper bound. There is no guarantee, however, that this bound can be attained. Note that in \eqref{dense z setup 4} the term $q_i(\xt)q_i(\xt)^{\top}$ is a projection matrix. If the current $\dfdxt$ is orthogonal to or badly aligned with the eigenvector associated with the maximum eigenvector, \eqref{dense z setup 4} will be very far from the theoretical upper bound.

A better design is to ensure that the gradient flow can progress in the worst-case scenario, i.e. when $\dfdxt$ is parallel to the eigenvector with $\lambda_{min}$. Selecting $w_{ii} \triangleq \lambda_{i}^{-1}$ changes the bound to,
\begin{align}
    \|\nabla f(\xt)\|^2\leq \dfdxt^{\top}\nabla^2 f(\xt) Z(\xt)^{-1}\dfdxt \leq \|\nabla f(\xt)\|^2.
\end{align}
This design of $W$ ensures a negative time derivative on the gradient norm for every $\dfdxt$, thus ensuring that the gradient shrinks to zero. 
\begin{equation}
    \frac{1}{2}\frac{d}{dt}\|\dfdxt\|^2 = -\dfdxt^{\top}\nabla^2 f(\xt) Z(\xt)^{-1}\dfdxt = -\|\nabla f(\xt)\|^2.
\end{equation}
This design is a gradient flow version of the Newton method, as $Z(\xt)^{-1}\triangleq \nabla^2 f(\xt)^{-1}$ when $W$ is transformed to $Z$. This result aligns with standard knowledge in optimization literature to the Newton method's utility when the Hessian and its inverse may be computed.

In practice, inverting the Hessian requires a prohibitively slow $\mathcal{O}(n^3)$ operations per query (for a dense matrix), which is the motivation for continued study into alternate methods. In this work, we will restrict $Z$ to the class of diagonal matrices that scale with respect to the canonical basis. There will be no need to compute the eigendecomposition or invert the Hessian.  

\subsection{Assumptions}
In this work, we will assume that the following set of assumptions are satisfied.

\begin{enumerate}[wide=\parindent,label=\textbf{(A\arabic*)}]
    \item $f\in C^2$ and $\inf_{\x\in\R^n}f(\x)>-R$ for some $R>0$. \label{a1}
    \item $f$ is coercive, i.e., $\lim_{\|\x\|\to\infty} f(\x) = +\infty$. \label{a2}
    \item $\nabla^2 f(\x)$ is non-degenerate. \label{a3}
    \item (Lipschitz and bounded gradients): for all $\x,\y\in\R^n$,  $\Vert \dfdx-\dfdy\Vert\leq L\Vert\x-\y\Vert$, and $\Vert\dfdx\Vert\leq B$ for some $B>0$. \label{a4}
    \item $Z(\x)^{-1}$ is diagonal for all $\x$ and $d_2>Z_{ii}(\x)^{-1}>d_1$ for all $i,\ \x$ and for some $d_1,\ d_2>0$.\label{z def}
\end{enumerate}


Note that $f$ is not assumed to be convex. Coercivity guarantees that there exists a finite global minimum for $f(\x)$ that is differentiable (\cite{peressini1988mathematics}). 

\begin{definition}
We say $\x$ is a \emph{critical point} of $f$ if it satisfies $\nabla_{\x}f(\x) = \vec{0}$. Let $S$ be the set of \emph{critical points}, i.e. $S=\{\x\ |\ \dfdx = \vec{0}\}$.
\end{definition}

The coercivity and differentiability of $f$ guarantee that any minima are within the set $S$.

\subsection{Component-Wise Scaled Gradient Flow}

To begin this optimization idea, we first establish convergence of gradient flow as presented in \eqref{scaled gradient flow} but with $Z$ restricted to the class of diagonal matrices. It can be shown that for any objective function $f$ and any scaling function $Z$ satisfying the stated assumptions, the gradient flow IVP will converge to some $\x$ within the set of critical points. 

\begin{theorem} \label{convergence theorem}
Under \ref{a1} - \ref{z def}, consider the component-wise scaled gradient flow IVP in \eqref{scaled gradient flow}. Then,
\begin{equation}
    \lim_{t\to\infty}\Vert \dfdxt\Vert = 0.
\end{equation}
\end{theorem}

The proof is provided in the supplementary material. Given Theorem \ref{convergence theorem}, the next objective in this paper is to design $Z$ that satisfies \ref{z def} and yields faster convergence than the base case of $Z=\mathbf{I}$. To tackle this objective, we next demonstrate that \eqref{scaled gradient flow} is mathematically equal to the transient response of an electrical circuit; therefore the circuit can be used to design the control matrix $Z$.

%% file: Meta-Alg.tex
\section{Method}
\label{meta-alg}

\begin{figure}
    \centering
    \includegraphics[width=0.9\columnwidth]{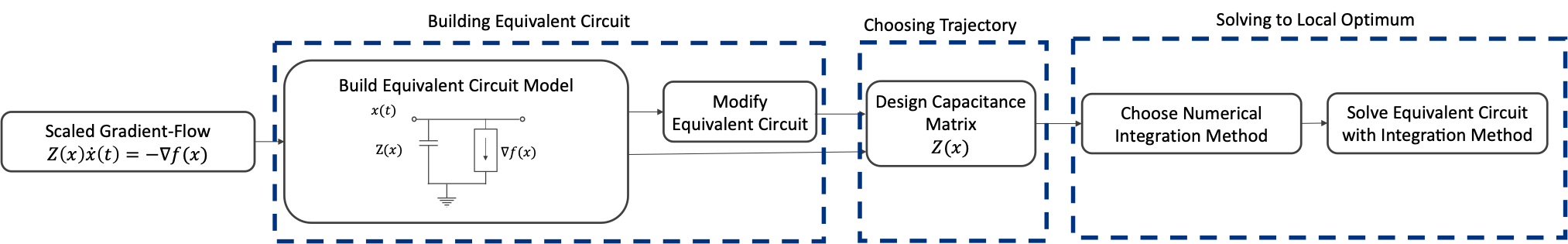}
    \caption{Solving the optimization problem using the ECCO workflow consists of three main steps: Designing the circuit model, choosing a trajectory and solving the circuit using a numerical integration method}
    \label{fig:ecco_meta_alg}
\end{figure}

The goal of the ECCO workflow (Figure \ref{fig:ecco_workflow}) is to design new optimization algorithms that utilize a physics-based circuit model. In this workflow, the scaled gradient flow equations are mapped to an equivalent circuit model, which is then designed and solved using domain knowledge to reach an optimal solution. This workflow promotes deriving insights from the equivalent circuit  that would not necessarily be considered without the physical model.
We outline the modular steps within the ECCO workflow (Figure \ref{fig:ecco_meta_alg}) that methodically incorporate circuit knowledge into new optimization methods.
This framework is ideal for customizing circuit-inspired methods for specific steps to achieve different effects such as wall-clock speed, ease of implementation, and iteration accuracy.

Solving an optimization problem in the ECCO workflow takes three major steps: (1) modeling the objective function with electrical circuitry, (2) controlling the circuit for fast convergence, and (3) numerically solving the circuit for its steady state as shown in the overall workflow (Figure \ref{fig:ecco_meta_alg}). 
We discuss the properties that each of these steps must generally fulfill in the following sections, as well as the relevant tradeoffs.

First, we will show that for any optimization problem \eqref{scaled gradient flow} satisfying the given assumptions, we can draw an \emph{equivalent circuit} (EC) whose transient response is identical to the controlled gradient flow equation. We will derive an all-purpose EC and show how optional modifications of the circuit may be added for better convergence properties. For any version of the circuit, we will show that EC steady states coincide with critical points of the optimization problem. 

Second, we derive a functional form for $Z$ using the EC model.  In the circuit domain, this is identical to designing a capacitor model. Using knowledge of charge and energy conservation, we propose physics-based techniques to choose the capacitances and control the trajectory of the optimization problem. We derive two approaches: one using the full Hessian, and an alternate first-order method. Both of these control schemes compel the circuit to a steady state, and thus likewise cause the ODE to trend to a fixed point.

Finally, the value of the steady-state must be found by solving the ODE via numerical methods. The transient response of the circuit is simulated to steady-state using numerical methods. We develop an adaptive time-stepping method for simulating the EC model that ensures accuracy and stability. We discuss heuristics such as initial step sizes derived from circuit simulation and adapt them for use in optimization algorithms. In addition, we present two explicit integration methods, Forward-Euler and fourth-order Runge-Kutta (RK4) to be used in ECCO.

%% file: equivalent_circuit_model.tex
\section{Designing the Equivalent Circuit}
\label{sec:ec_model}

\subsection{All-Purpose EC}
First, we must design an EC whose transient response is equal to \eqref{scaled gradient flow}. For $\x\in\R^n$, a multidimensional gradient flow may be represented by a set of $n$ coupled sub-ECs, each representing the response of one $\x_i(t)$. 

\begin{figure}
    \centering
    \includegraphics[width=0.5\columnwidth]{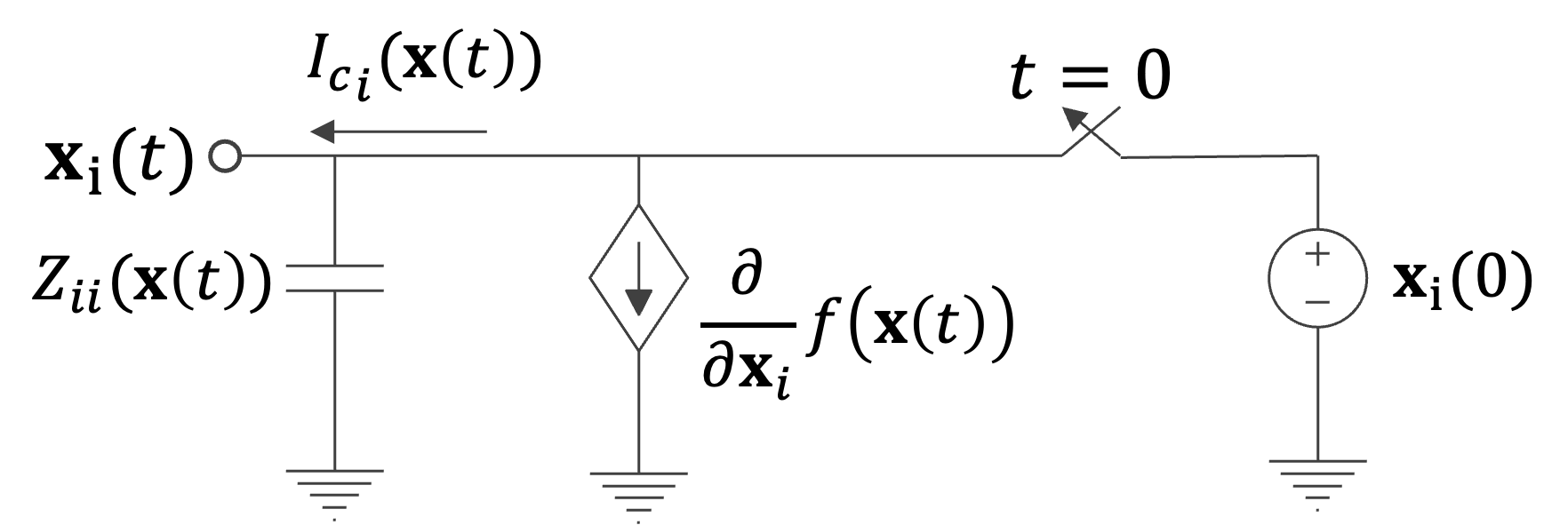}
    \caption{Equivalent Sub-Circuit Model}
    \label{fig:ECorig}
\end{figure}

The sub-circuit in Figure \ref{fig:ECorig} is a physical representation of the transient waveform of a single dimension in the scaled gradient flow. Each sub-EC model is composed of four circuit elements: a nonlinear capacitor, voltage-controlled current source (VCCS), an ideal switch and a DC voltage source. We next derive the transient voltage-current equations of the circuit shown in Figure \ref{fig:ECorig} and show that they are equivalent to the scaled gradient flow \eqref{scaled gradient flow}.

Initially for $t<0$, the ideal switch in each sub-circuit is closed, thus directly connecting the DC voltage source to the capacitor and VCCS. During this time, the voltage is set by the voltage source at $\x_i(0)$. Note that this value is the initial condition given in \eqref{scaled gradient flow}. At $t=0$, the ideal switch is switched opened, thereby disconnecting the DC voltage source and letting the capacitor and VCCS drive the node voltage.

The circuit then may be analyzed using nodal analysis \cite{pillage1995electronic}, which is a formal method that characterizes Kirchhoff’s Current Law (KCL) at each node. 
Applying KCL to each sub-EC for $t>0$, we see the currents leaving the node of the sub-circuit are the current from the non-linear capacitor, $I_{c_i}(\x(t))$, and the voltage-controlled current source, $\partial f(\xt)/\partial \x_i$. The response of the entire circuit is modeled by:
\begin{align}
    I_c(\x(t)) + \nabla f(\x(t)) &=0\ \forall\ t\geq0,     \label{eq:KCL_orig} \\
    \x(t) &= \x(0)\ \forall\ t\leq0. \label{kcl_orig_lzero}
\end{align}

Each capacitor in the sub-EC satisfies the voltage-current equation, 
\begin{equation}
    I_{c_i}(\x(t)
    ) = Z_{ii}(\x(t))\dot{\x}_i(t)
    \label{eq:cap_eq}
\end{equation}
The voltage response across the capacitor, $\x_i(t)$, represents the trajectory of a single state variable. The VCCS is controlled by the voltage vector, $\xt$, and produces a current equal to the gradient component, $\partial f(\xt)/\partial \x_i$.

Adding the capacitor current definition \eqref{eq:cap_eq} back to the KCL equation produces the following, which is identical to the scaled gradient flow \eqref{scaled gradient flow}:
\begin{equation}
    Z(\x(t))\dot{\x}(t) + \nabla f(\x(t)) = 0.\label{charge}
\end{equation}

By Theorem \ref{convergence theorem}, the scaled gradient-flow will converge to a critical point of the objective function; as the EC models the same equations, it will also reach a steady state that is the critical point of the objective function. This logic may be confirmed by observing the behavior of the circuit. 

The steady-state for the EC is achieved when the the voltage across the capacitors stops changing. At this state, the currents through the capacitors become zero as $\dot{\x}(t)=0$, and the system reaches an equilibrium. Applying KCL to the EC, we find that the current from the VCCS must also equal zero and therefore $\nabla f(\x)=0$. This condition is equivalent to the definition of a critical point for the objective function; therefore, we can confirm that a steady state of the EC is a critical point of the optimization problem.


The EC model establishes a physical model from which to derive all subsequent methods. In particular, the capacitor models the dynamics of the optimization variable, and the speed of convergence in continuous time can be controlled via the capacitance value, $Z(\x)$. A large capacitor, $Z_{ii}$ causes the node-voltage to change slower. The relative value between the sub-EC capacitances, $Z_{ii}$, dictates changes in the trajectory of the uncontrolled gradient-flow. To choose these values, we further develop our understanding of the optimization trajectory's behavior via an adjoint circuit model.

\subsection{Adjoint Circuit Model}

\begin{figure}
    \centering
\includegraphics[width=0.5\columnwidth]{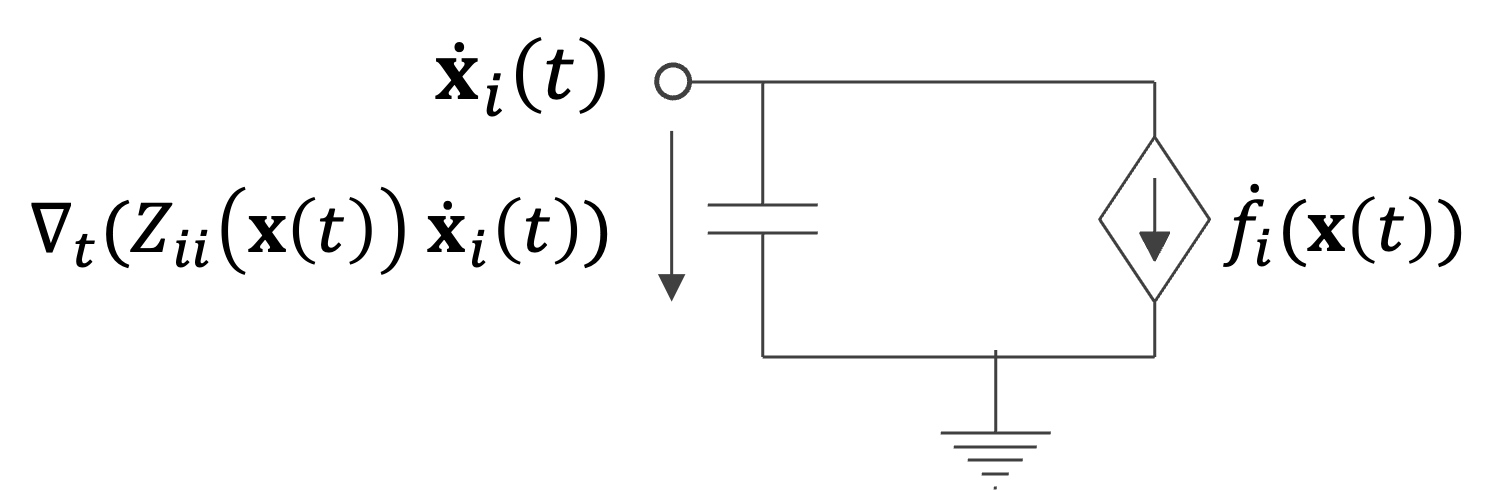}
    \caption{Adjoint Equivalent Circuit}
    \label{fig:ECadj}
\end{figure}

To gain further insight into the energy stored in the capacitors, we next construct an \emph{adjoint} EC that models the time derivative of the original gradient flow \eqref{scaled gradient flow}:

\begin{equation}
    \frac{d}{dt}(\Zeta(\x(t))\dot{\x}(t)) = - \frac{d}{dt}\nabla f(\x(t)) \;\;\; \dot{\x}(t<0) = \dot{\x}(0) \label{adjoint equation}
\end{equation}
A sub-EC for \eqref{adjoint equation} is shown in Figure \ref{fig:ECadj}, where now the node voltage is $\dot{\x}(t)$. This time derivative circuit model will be used to analyze the convergence properties of the original gradient flow. Importantly, $\xt$ in both the original EC and the adjoint EC are the same variable and have the same trajectory under the given assumptions.

Similar to the EC model, the adjoint circuit is composed of an adjoint capacitor model and an adjoint VCCS. Adjoint capacitor $\bar{c}_i$ produces the current
\begin{equation}
    \bar{I}_{\bar{c}_i}(t) = \frac{d}{dt}(\Zeta_{ii}(\x(t))\dot{\x}_i(t)),
\end{equation}
which is the negative of the current produced by the adjoint VCCS.

The adjoint EC will be key to the development of our optimization theory as it allows us to study the electrical energy related to the optimization variable within its trajectory. We will now introduce the concepts of charge, energy, and passivity that will be used in future sections to design capacitor models. 


The electrical charge stored in each adjoint capacitor is defined as , 
\begin{equation}
    \bar{Q}_{c_i}(t) = \Zeta_{ii}(\x(t))\dot{\x_i}(t),
\end{equation}
where charge is derived as the time integral of current. The energy stored by the capacitor in the adjoint circuit can be derived from the charge by,
\begin{equation}
    \bar{E}_{c_i}(t) = \frac{\bar{Q}_{c_i}^2(t)}{2\Zeta_{ii}(\x(t))}.
\end{equation}

The energy stored in each adjoint capacitor is proportional to the square of charge stored $\bar{Q}_{c_i}^2(t)$.

Again, we can discuss the steady state condition which occurs when $\dxdt(t) = 0$, and whose existence can be confirmed from Theorem \ref{convergence theorem}. 
Initially, the adjoint capacitors store an electrical charge equal to $\bar{Q}_c(t) = Z(\x(t))\dot{\x}(0)$. During the transient response, the initial charge dissipates from the adjoint capacitors to reach a steady-state condition where $\dxdt(t)=0$ for $t>t'$. At this state, the charge in the adjoint capacitors will be $\bar{Q}_c(t)=0 \forall t>t'$, indicating that all the initial energy stored in the adjoint capacitors dissipated through a set of adjoint VCCS elements whose power consumption is greater than zero during the transient response. 


\begin{definition}
    The power consumption of an adjoint VCCS element, $\bar{P}_{i}^{VCCS}$, with a voltage potential $\dot{\x}_i$ and a current $\frac{d}{dt}\frac{\partial \nabla f(\xt)}{\partial \x_i} $ is :
    \begin{equation}
        \bar{P}_{i}^{VCCS}=\frac{d}{dt}\frac{\partial \nabla f(\xt)}{\partial \x_i} \dot{\x}_i(t)
    \end{equation}
\end{definition}

%% file: modifying_ec.tex
\label{sec:modifying_ec}
\subsection{Modifying the EC Model via Homotopy}

The equivalent circuit model of the scaled gradient-flow establishes a physical model from which to derive new insights and optimization methods. During the first stage of the workflow, we can modify the all-purpose EC in Figure \ref{fig:ECorig} to induce desirable effects into the optimization trajectory. For example, circuit simulation regularly modifies circuits to find a DC operating point \cite{pillage1995electronic} through the use of homotopy methods. 

Homotopy is a general class of solution methods whereby a difficult numerical problem is approached by solving a series of incremental, simpler problems. Consider the set of equations $F(w)=0$ and an easier set of equations $E(w)=0$, representing a trivial problem. The homotopy approach defines a series of sub-problems, \cite{allgower1980homotopy}
\begin{equation}
    H(w,\gamma) = F(w,\gamma) + E(w,\gamma)\equiv0,
\end{equation}
where $\gamma\in [0,1]$ is scalar quantity embedded into the sets of equations that control the mixed contribution of $F$ and $E$. The impact of $\gamma$ is designed such that $H(w,\gamma = 0)=E(w)$ and $H(w,\gamma = 1)=F(w)$. By incrementing $\gamma$ from 0 to 1 in discrete steps of $\Delta \gamma$, we can define a series of sub-problems that start from the easier problem $E$ and finish at the original difficult problem $F$. Each sub-problem $H(w,\gamma)$ is solved using the solution to the previous sub-problem, $H(w^*,\gamma-\Delta \gamma)=0$, as an initial condition. This traces a continuous path in the solution space. 
Well-designed homotopy methods often have a smoother trajectory to an optimal than solving $F(w)$ directly.

Given a goal problem, $F(w)=0$, there are numerous ways to embed $\gamma$ and define $E(w)$. However, an understanding of the physics of the underlying model can provide insights into developing a homotopy method. 

\subsubsection{Source Stepping}
Using knowledge of circuit behavior, we can propose a homotopy method to modify the EC model based on the technique of \emph{source stepping}, which has a natural analogy in circuit simulation. 

Consider the voltage source on the right side of the original EC in Figure \ref{fig:ECorig}. During the period where $t<0$, this source fully charges the capacitors and sets the initial node voltages. At $t=0$, we abruptly disconnect the voltage source, letting the EC be driven by the scaled gradient flow equations. Instead of an abrupt switch, we can use a homotopy method that incrementally reduces the contribution of the DC source.

\begin{figure}
    \centering
    \includegraphics[width=0.5\columnwidth]{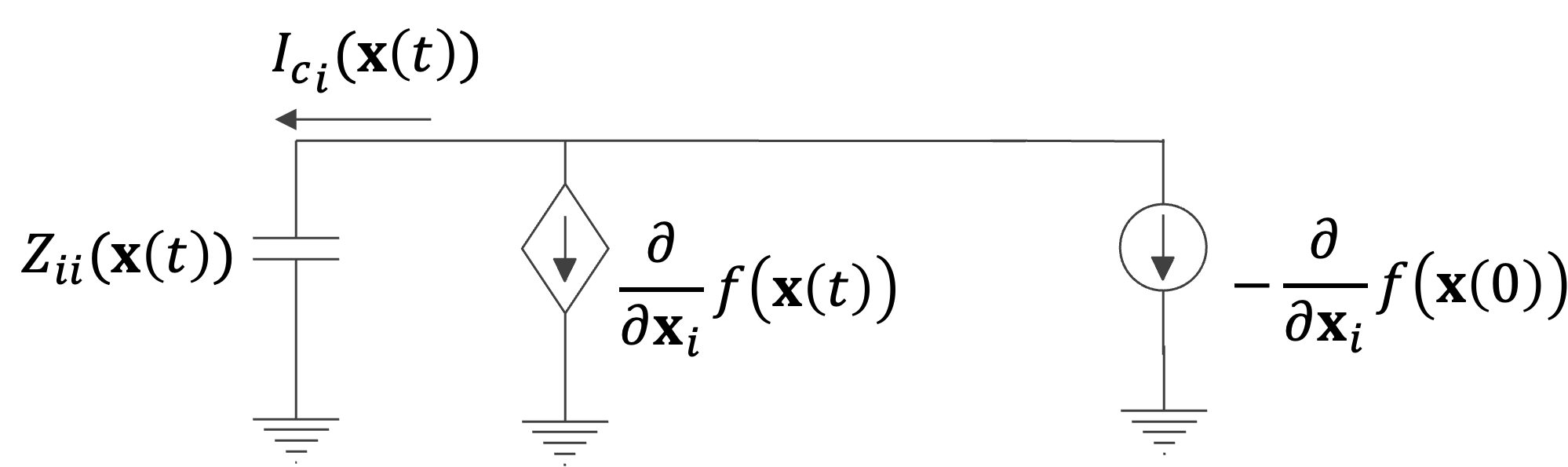}
    \caption{EC Model with DC Current Source for $t<0$}
    \label{fig:ec_current_source}
\end{figure}

 To iteratively reduce the contribution of the DC source, we first transform the DC voltage source to its current source equivalent. This is accomplished via the substitution principle \cite{pillage1995electronic}, and yields a current source that supplies $-\nabla f(x(0))$ for $t<0$. This new EC model in Figure \ref{fig:ec_current_source}, which is electrically identical to the model in Figure \ref{fig:ECorig} for $t<0$, represents our easier problem $E(\xt) = Z(\xt)\dxdt(t)+ \nabla f(\x(t))- \nabla f(\x(0)) = 0$, which has a trivial solution of $\x^*=\x(0)$.

 In this homotopy method, we gradually reduce the contribution of the current source and solve for the steady state of the new system using the steady state of the previous iteration as an initial condition. With a small enough decrement in the current source, we can ensure that at each iteration the next steady state is close to its initial conditions. This traces a smooth path in the solution space, where the solution to each sub-problem is close to the previous sub-problem. 



To define the source stepping approach, we linearly mix the $E$ and $F$ terms with a scalar homotopy factor, $\gamma_{ss}$ and define a series of problems, $H(x,\gamma_{ss})$ as
\begin{equation}
\begin{split}
    H(\x, \gamma_{ss}) &= \gamma_{ss}F(\x) + (1-\gamma_{ss})E(\x),\\
    &= \gamma_{ss}(Z(\x(t))\dot{\x}(t) + \nabla f(\x(t))) + (1-\gamma_{ss})(Z(\xt)\dxdt(t)+ \nabla f(\x(t)) - \nabla f(\x(0))), \\
    &= Z(\x(t))\dot{\x}(t) + \nabla f(\x(t)) - (1-\gamma_{ss}) f(\x(0)) \equiv 0. 
\end{split}
\end{equation}

\begin{figure}
    \centering
    \includegraphics[width=0.5\columnwidth]{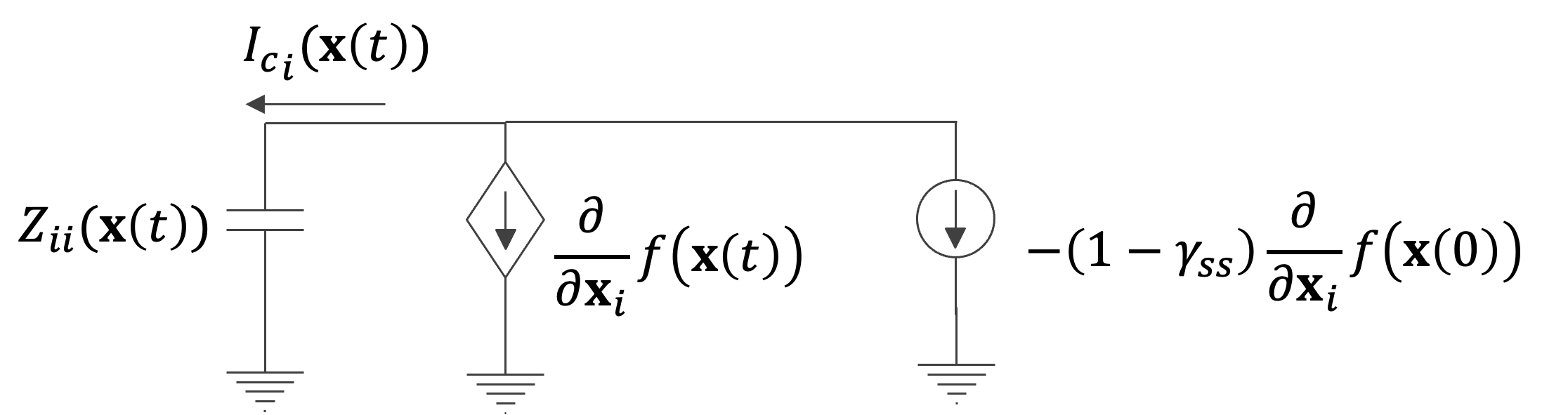}
    \caption{EC Model for Source-Stepping}
    \label{fig:ec_source_stepping}
\end{figure}

The EC model is modified according to Figure \ref{fig:ec_source_stepping}. The gradient flow is likewise modified to the following:
\begin{equation}
    \Zeta(\x)\dot{\x}(t) = -\nabla f(\x(t)) + (1-\gamma_{ss})\nabla f(\x(0)).
\end{equation}

The full algorithm for the homotopy approach is outlined in Algorithm \ref{source-stepping_overview}

\begin{algorithm}
\caption{Source Stepping ECCO}
\label{source-stepping_overview}
\textbf{Input: } $\nabla f(\cdot)$,$\x(0)$, $\Delta \gamma_{ss}\in (0,1)$
\begin{algorithmic}[1]
\STATE{$\gamma_{ss} \gets 0$}

\STATE{\textbf{do while} $\gamma_{ss} > 0$}

\STATE{\hspace*{\algorithmicindent} Solve for $\x_{ss}$: $\Zeta(\xt)\dot{\x}(t) = -\nabla f(\x(t)) + (1-\gamma_{ss})-\nabla f(\x_0)$}

\STATE{\hspace*{\algorithmicindent}$\gamma_{ss} \gets \gamma_{ss}+\Delta \gamma_{ss}$}
\RETURN $\x$
\end{algorithmic}
\end{algorithm}



The source stepping method uses the initial conditions of the optimization to define a modified EC with a known feasible steady-state and gradually deforms the underlying circuit to the original EC. This process provides a smooth trajectory for optimization, which we can capitalize during the numerical integration process by taking larger time-steps to reach the steady-state of each homotopy iteration.

%% file: Control.tex
\section{Controlling the EC}
\label{sec:control}

The next step in the ECCO workflow is to design a scaling matrix $Z$ that enables quick convergence to the fixed point of the gradient flow. Instead of viewing this as a numerical problem, we leverage the EC and time-derivative EC models and treat choosing $Z$ as designing the nonlinear circuit capacitances.  In the circuit sense, fast convergence to a critical point of the optimization problem is equivalent to fast convergence to a stable steady-state of the dynamical system in continuous time. Thus, the capacitances can be chosen as a function of $\xt$ to control for fast convergence to  steady state. Through the circuit properties it is known that smaller capacitances correspond to faster convergence to a steady-state. By modifying the individual capacitances, i.e. $Z_{ii}$, the overall optimization trajectory can also be changed. 

We pause here to emphasize that any capacitance that satisfies Assumption \ref{z def} will yield asymptotic convergence. In Section \ref{sec:examples}, examples will show how the existing optimization methods of gradient descent and Nesterov's accelerated method can be derived from the gradient flow equations based on a simple construction of $Z$. In the rest of this section, however, we will focus on construction of new $Z$ with the explicit purpose of fast convergence of the continuous time ODE. 

To choose $Z$ values, we analyze the energy and charge stored in the adjoint capacitors. During the transient response of the time-derivative EC, the initial charge in the adjoint capacitors ($\bar{Q}_c(\x(0))$) is  discharged through adjoint VCCS elements to reach a steady-state conditions where $\bar{Q}_c=0$. In this work, we design the capacitors models, $Z(\x)$, to quickly discharge the adjoint capacitors and achieve fast convergence to a zero-energy steady-state, i.e., the critical point of the objective function.


The method for choosing $Z$ is another step in the workflow where a user may customize the implementation. In the following sections we present two models for $Z$ based on the idea of maximizing the rate of charge dissipation of the adjoint capacitor that may be used by a practitioner directly. One scheme uses the full Hessian for applications in which the Hessian is easily computed, and the other option is a first order approximation of the Hessian-based result. If neither of these schemes suit a desired application, a practitioner may derive other valid controls from the EC; if the proposed $Z$ satisfies the boundedness conditions in Assumption \ref{z def}, then convergence will hold.


\subsection{Second Order Control}
The squared charge stored in the adjoint circuit capacitor at some time $t$ is equal to $\Vert \bar{Q}_c(\xt)\Vert^2$, which by \eqref{charge} is equal to  $\Vert Z(\xt)\dot{\x}(t))\Vert^2$ and by \eqref{scaled gradient flow} is equal to $\Vert \nabla f(\xt)\Vert^2$. Thus to quickly \emph{dissipate} charge, we choose $Z$ by solving the following optimization problem, which maximizes the negative time gradient of the charge:

\vspace{-0.5cm}
\small
\begin{align}
    &\max_{Z}- \frac{1}{2} \frac{d}{dt}\Vert \bar{Q}_c(\xt) \Vert^2,\\
    &=\max_{Z}- \frac{1}{2} \frac{d}{dt}\Vert \dfdxt \Vert^2,\\
    &=\max_{Z} - \dfdxt^{\top} \frac{d}{dt} \dfdxt, \label{passivity_max_Z}\\
    &=\max_Z \dfdxt^{\top}\nabla^2 f(\xt)Z^{-1}\dfdxt,\label{max dvdt}
\end{align}
\normalsize
where \eqref{max dvdt} is found via Lemma \ref{identity2}. 


To gain further intuition on the role of $Z$, we can expand the optimization problem \eqref{passivity_max_Z} as:
\begin{align}
    &\max_{Z}  Z(\x)\dot{\x}(t) \frac{d}{dt} \dfdxt. \label{eq:passivity_max_Z_VCCS}, \\
    &\max_{Z}  Z(\x)\bar{P}^{VCCS} . \label{eq:power_max_Z_VCCS}
\end{align}
From the circuit perspective, maximizing the rate of charge in \eqref{max dvdt} is identical to maximizing the power dissipation through the adjoint VCCS elements (the power dissipation, $\bar{P}^{VCCS}$ is equal to $\dot{\x}(t) \frac{d}{dt} \dfdxt$) in \eqref{eq:power_max_Z_VCCS} to quickly reach the zero-energy steady-state. Given the assumption that $Z_{ii}>0$, the goal of the optimization should be to align $Z(\x)$ with adjoint VCCS components, where $\bar{P}^{VCCS}_i>0$. In this work, we choose a direction that increases the total power dissipation through dissipative adjoint VCCS elements and ensures that the term $Z(\x)\dot{\x}(t) \frac{d}{dt} \dfdxt$ is positive. 

The optimization problem \eqref{max dvdt} must yield a diagonal $Z$. To this end, we expand $\dfdxt$ to a diagonal matrix and shrink $Z^{-1}$ to a vector. Define $G(\xt)$ be a diagonal matrix where the diagonal elements are the gradient $G_{ii}(\xt) = \frac{\partial f(\xt)}{\partial \x_i(t)}$. Let $\mathbf{z}$ be a vector where $\mathbf{z}_i=Z_{ii}^{-1}$. Then \eqref{max dvdt} is equal to,
\begin{align}
    &=\max_{\z} \dfdxt^{\top}\nabla^2 f(\xt)G(\xt)\z-\frac{\delta}{2} \Vert \z\Vert^2,
\end{align}
where a regularization term with $\delta>0$ has been added for tractability. Taking the derivative $\frac{\partial}{\partial \z}$:
\begin{gather}
    G(\xt)\nabla^2f(\xt) \dfdxt -\delta \z \equiv \vec{0}.\\
    \z = \frac{1}{\delta} G(\xt) \nabla^2 f(\xt)\dfdxt.\label{final z}
\end{gather}
Any negative value of $\z_i$ occurs when $ G_{ii}(\xt)\nabla^2f(\xt) \dfdxt = Z_{ii}(\x)\dot{\x}_i(t) \frac{d}{dt} \dfdxt$ is negative. This indicates that the VCCS element is generating power rather than dissipating it (i.e., $\bar{P}^{VCCS}_i<0$). To suppress the trajectory from going towards the non-passive adjoint VCCS elements, we truncate $\z_i$ to a small value. In the circuit model, a small $z_i$ places large capacitance at the corresponding node that causes the node-voltage to change slower. To ensure that $Z^{-1}$ performs as least as well as the base case where $Z^{-1}=I$ (and to retain positivity and invertibility of $Z$), any $\z_i<1$ is truncated to $\z_i=1$. This also guarantees each capacitor in the EC model has a positive capacitance ($Z_{ii}>0$), which is necessary for the capacitors to behave passively.

For some $\delta>0$, the final construction of the control matrix $Z(\xt)^{-1}$ is:
\begin{equation}
    Z_{ii}(\xt)^{-1} = \max\{\delta^{-1} [G(\xt) \nabla^2 f(\xt)\dfdxt]_i,1\}.\label{z true}
\end{equation}

We find that truncating any $\z_i<1$ to $\z_i=1$ performs well in practice as this ensures that $Z^{-1}$ performs at least as well as the base case where $Z^{-1}=I$, i.e. uncontrolled gradient flow. This truncation choice also ensures positivity and invertibility of $Z$. For some $\delta>0$, the final construction of the control matrix $Z(\xt)^{-1}$ is thus:
\begin{equation}
    Z_{ii}(\xt)^{-1} = \max\{\delta^{-1} [G(\xt) \nabla^2 f(\xt)\dfdxt]_i,1\}.\label{z true}
\end{equation}

Equation \eqref{z true} is a second order method as it uses Hessian information.

\begin{lemma} \label{bounded z true}
    Assume \ref{a1}-\ref{a4} hold. The definition of $Z(\xt)^{-1}$ in \eqref{z true} satisfies \ref{z def}.
\end{lemma}
\begin{proof}
    See Appendix \ref{zbounded}.
\end{proof}
\begin{lemma} \label{z true complexity}
    Given the gradient and Hessian, the computation complexity to evaluate $Z^{-1}$ at a specific $\xt$ is $\mathcal{O}(n^2)$.
\end{lemma}
\begin{proof}
    See Appendix \ref{complexity full}.
\end{proof}

In comparison to Newton methods, this $Z$ does not require a Hessian inversion step. 

\subsection{First Order Control}
\label{sec:acontrol}
The proposed control scheme may be used directly as in \eqref{z true}; however, it requires computation of the full Hessian. For applications in which the Hessian in unavailable or expensive to compute, such as machine learning, we present an approximation that may be computed from only gradient information.

Define the finite difference approximation, i.e. the first order Taylor expansion, of the optimization trajectory $\frac{d}{dt}\dfdxt$,
\begin{align}
    \ahat(\xt)\triangleq \frac{\dfdxt - \nabla f(\x (t-\Delta t))}{\Delta t}. \label{ahat}
\end{align}
The limit as $\Delta t\to 0$ is exactly equal to desired quantity $\frac{d}{dt}\nabla f(\xt)$, making this an apt approximation for small $\Delta t$. Based on \eqref{ahat}, approximation for $Z^{-1}$ can be defined:
\begin{equation}
     \widehat{Z}_{ii}(\xt)^{-1}=\sqrt{\max\{[-\delta^{-1}G(\xt)\ahat(\xt)]_i,1\}}. \label{z approx}
\end{equation}

The derivation is shown in Appendix \ref{z approx deriv}.

\begin{lemma}\label{bounded z approx}
    Assume \ref{a1}-\ref{a4} hold. The definition of $Z(\xt)^{-1}$ in \eqref{z approx} satisfies \ref{z def}.
\end{lemma}
\begin{proof}
    See Appendix \ref{zhat bounded}.
\end{proof}
\begin{lemma}\label{z approx complexity}
    Given the gradients, the computation complexity to evaluate $\widehat{Z}^{-1}$ at a specific $\xt$ is $\mathcal{O}(n)$.
\end{lemma}
\begin{proof}
    See Appendix \ref{complexity approx}.
\end{proof}

We form $\widehat{Z}$ by approximating the trajectory of the optimization variable at a specific $t$ and solving for $Z$ as a function of $\ahat$; we do not estimate the Hessian matrix, and therefore $\widehat{Z}$ is not a quasi-Newton method. The proposed iterative update using \eqref{z approx} is a scaled gradient or step-size normalized descent method, similar in spirit to normalized gradient flow (in continuous time) and step-size or momentum scaling methods methods such as Adam, RMSprop, and Adagrad.

%% file: Discretization.tex
\section{Solving the EC with Numerical Integration}
\label{sec:discretization}

Given the EC model and the choice of $Z$, the next step in the ECCO workflow is to solve for the steady state. To accomplish this goal, we simulate the transient response of the EC to steady-state using  numerical integration methods, which are widely used in circuit simulation. In this work, we will identify techniques used in circuit simulation that are currently overlooked in optimization, and adapt them to solve the scaled gradient flow problem.  

In simulating the trajectory of $\xt$, we march through time in discrete time-steps, $\Delta t$. The state at each time-point is given by \eqref{eq:discretization_state}.
\begin{align}
    &\x(t+\Delta t) =\x(t) + \int_{t}^{t+\Delta t} \dot{\x}(t) dt  =  \x(t) + \int_{t}^{t+\Delta t} -\Zeta(\x(t))^{-1} \nabla f(\x(t)) dt. 
    \label{eq:discretization_state}
\end{align}

In general, the integral on the RHS will not have an analytical solution and must instead be approximated over a time step $\Delta t$ using numerical integration techniques. The time-step mirrors the role of a step size, i.e. learning rate, in optimization methods, which dictates the size of the approximation to the objective function. 



\begin{figure}
\centering
\subfloat[EC for \eqref{eq:multisim_example}.]{\label{mult_ec}\includegraphics[width=0.35\linewidth]{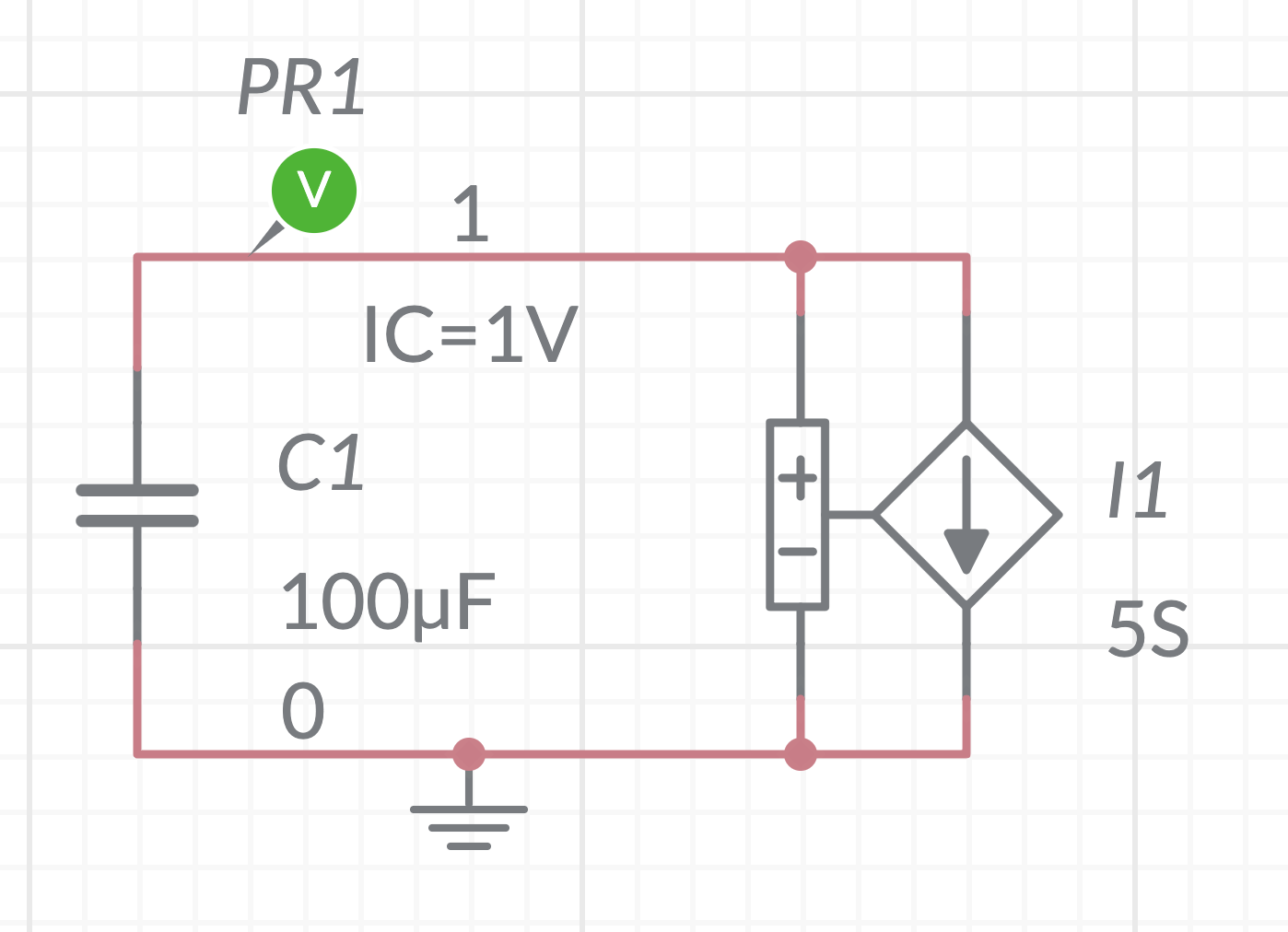}}
\subfloat[Transient response for EC in Multisim.]{\label{mult_trans}\includegraphics[width=0.35\linewidth]{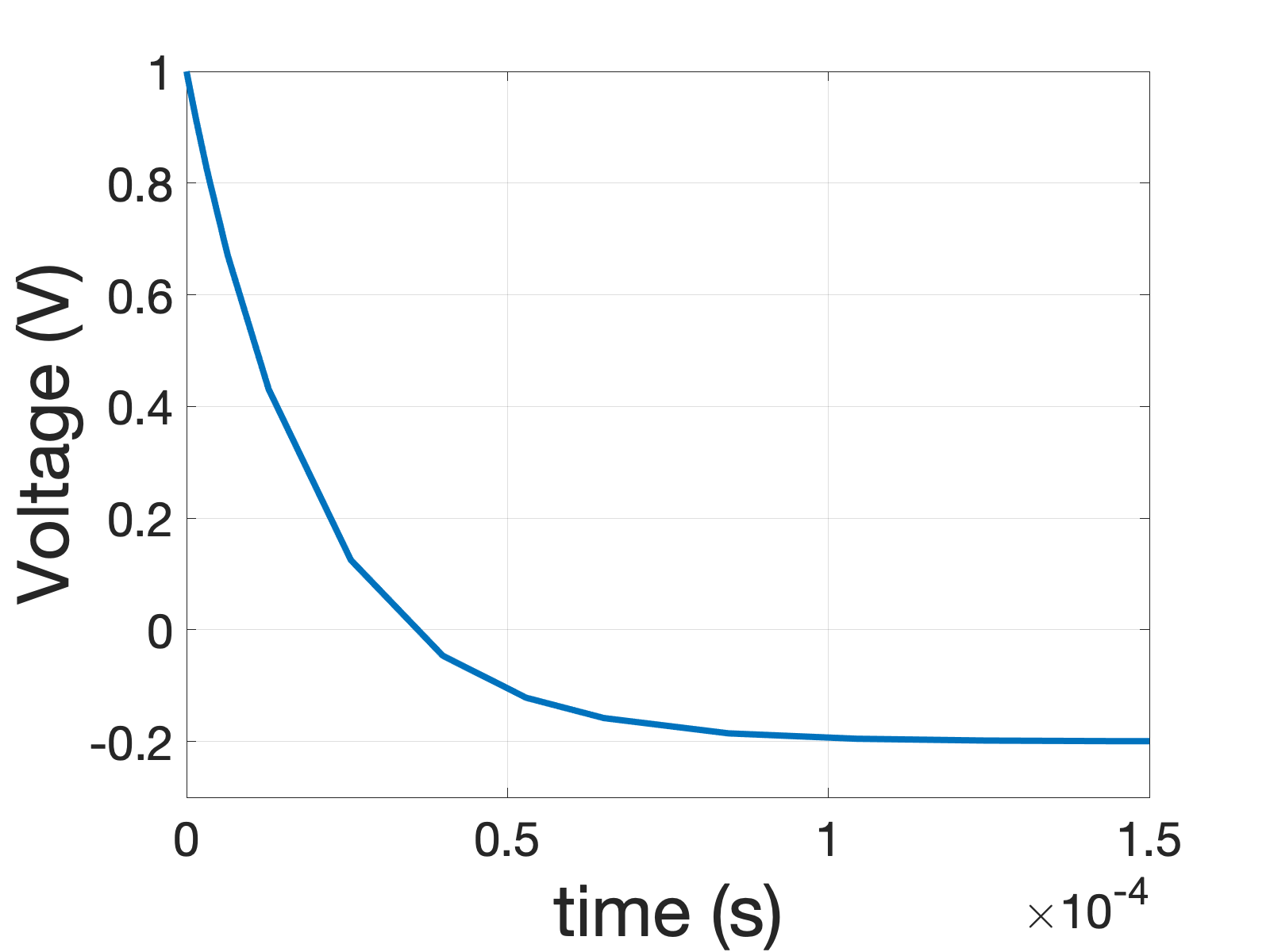}}
\caption{Solving \eqref{eq:multisim_example} via EC using Multisim circuit simulator.} 
\label{fig:testfig}
\end{figure}

One approach for solving for the critical point of the scaled gradient flow is to directly use commercial circuit simulation softwares or standard numerical solvers (such as the MATLAB ODE toolbox) to simulate the transient response of the EC. For example, we can solve the following scalar optimization problem by building an EC model and solving using a commercial circuit simulator (Multisim \cite{guide1989multisim}):
\begin{equation}
    min_{\x} \frac{5}{2}\x^2 +\x. 
    \label{eq:multisim_example}
\end{equation}
The EC model for the optimization problem \eqref{eq:multisim_example}, shown in Figure \ref{mult_ec}, is initialized at a voltage of $\x(0)=1V$, i.e. an initial condition of $\x=1$. The transient response, shown in Figure \ref{mult_trans}, settles to a steady-state voltage of $\x=-0.2V$, indicating that the solution to the optimization problem in \eqref{eq:multisim_example} is $\x^*=-0.2$. 

Internally, the Multisim engine uses a semi-implicit trapezoidal numerical integration method \cite{guide1989multisim} to solve for the transient response. While general purpose circuit simulators, such as Multisim, are designed to provide accurate transient responses of large-scale stiff circuits, the goal of the EC approach for optimization is to reach the desired steady state in the fewest time steps. 
In this section, we provide a general algorithm for an adaptive time step selection that incorporates principles of numerical integration, circuit simulation, as well as optimization. This allows us to generalize step size selection as a choice of numerical integration to achieve varying effects on the optimization method. We present examples of two common numerical integration methods as well as an initial time step selection based on circuit principles. 

\subsection{Adaptive Time Step Selection}
Any integration technique used to approximate \eqref{eq:discretization_state} will have numerical properties that are influenced by the size of the time step $\Delta t$. Step sizes that are too large mean that the approximation to the integral in \eqref{eq:discretization_state} will be used over too large of a period, yielding inaccurate results or divergence. If the step sizes are too small, however, computing the steady state will be very slow due to the number of time-points evaluated. Combining ideas from circuit simulation and optimization, we propose that the time step must satisfy conditions that preserve numerical \emph{accuracy} and \emph{stability}.



\begin{figure}
    \centering
    \includegraphics[width=0.7\linewidth]{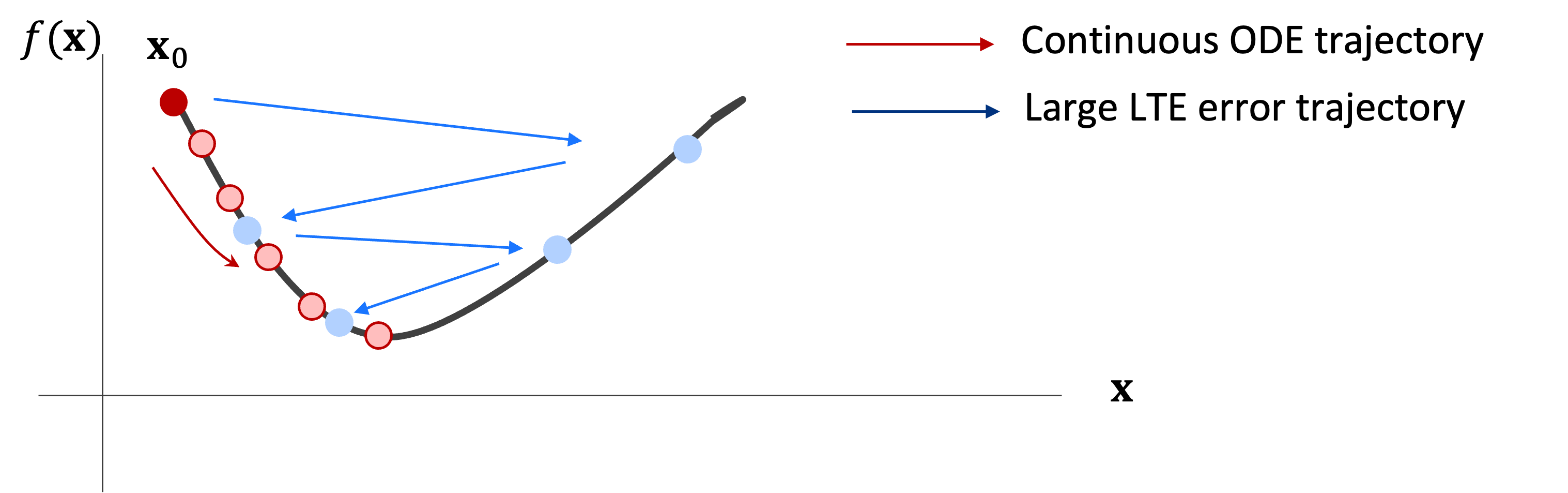}
    \caption{Effect of different time step sizes in discretizing  the Continuous ODE trajectory}
    \label{fig:lte}
\end{figure}

\paragraph{Accuracy Condition} 
To ensure that the discrete approximation is not diverging from the continuous-time trajectory, we track the Local Truncation Error (LTE)\cite{pillage1995electronic}, which measures the accuracy of a single numerical integration step. The calculation of LTE will depend on the specific numerical integration technique used. 
In the context of optimization, a step size with a large LTE risks diverging from the continuous-time waveform, resulting in numerical oscillations or yielding meaningless values as illustrated in Figure \ref{fig:lte}.


While we do not require a highly accurate calculation of the entire gradient flow trajectory, the estimated values must still be accurate enough to be usable. A tolerance $\eta$ must be selected; in general we find that $\eta \in (0.001, 1)$ performs well. The final accuracy condition becomes, 
\begin{equation}
    max(|LTE|) \leq \eta. \label{LTE cond}
\end{equation}

\paragraph{Stability Condition}


The approximations made by numerical integration methods at each time step are capable of accumulating over a time-window, thereby diverging from the intended steady-state. To resolve this, standard circuit methods ensure that the signal approaches steady state at each time-step. Unlike in circuits, the steady state value is unknown in \eqref{scaled gradient flow}; however, we know that the steady-state yields a minimum of the objective function. Thus, we can introduce a stability check that ensures monotonic decrements to the objective function at each time-step. For $\x(t+\Delta t)\neq \x(t)$,
\begin{equation}
    f(\x(t+\Delta t)) < f(\x(t)). \label{conv cond 1}
\end{equation}

This monotonic condition is analogous to the Armijo condition \cite{wolfe1969convergence}, commonly used in line search methods:
\begin{equation}
    f(\x(t+\Delta t))< f(\x(t))+c\Delta t\nabla f(\xt),\label{conv cond}
\end{equation}
for some $0\leq c<1$. As \eqref{conv cond} is standard in optimization and a generalization of \eqref{conv cond 1}, we select it as the stability criterion. This check ensures that the discretized trajectory goes toward a local minimum, thus converging to steady state and avoiding numerical instability issues.

\begin{figure}
    \centering
    \includegraphics[width=0.3\linewidth]{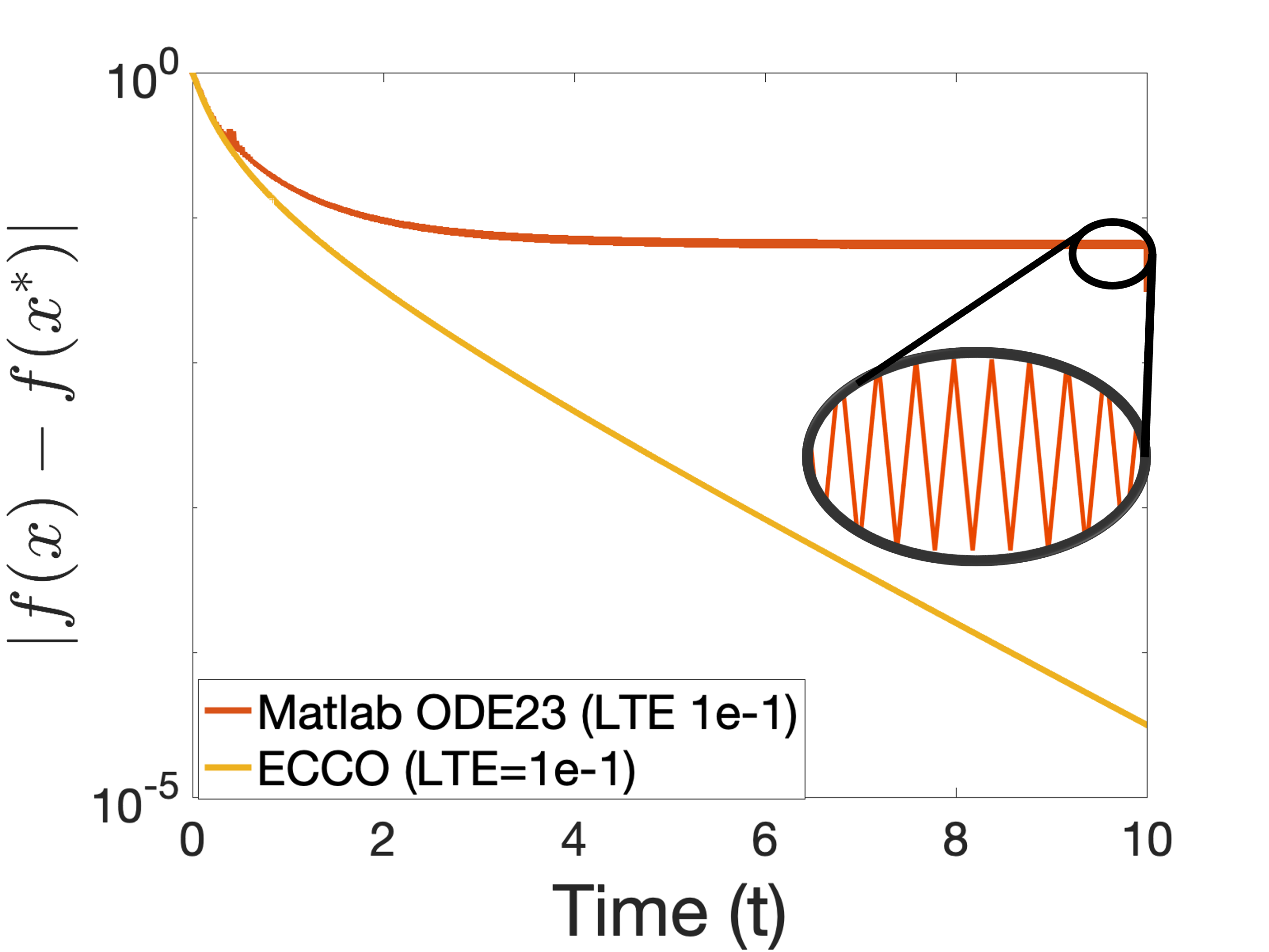}
    \caption{Simulating the transient response of the EC model of the Rosenbrock function using Matlab ODE23 (an explicit integration method) and FE discretization with EATSS}
    \label{fig:bad matlab}
\end{figure}

Withouth this stability criteria, off-the-shelf explicit ODE solvers are ill-suited to solving gradient flows. This is exemplified in Figure \ref{fig:bad matlab}, which logs the error while optimizing the Rosenbrock test function \cite{andrei2008unconstrained} via gradient flow. The error trajectory produced by MATLAB's ODE23 implementation (an explicit method) fails to converge to the steady state and instead oscillates due to numerical instability of the explicit integration method. In comparison, solving the same problem with the added stability check enables computation of the correct answer.

\begin{algorithm}
\caption{Error Aware Time Step Search (EATSS)}
\label{time step algo}
\textbf{Input: }{$f(\cdot)$,  $\nabla f(\cdot)$, $\x(t)$, $Z(\xt)$, $(\Delta t)_0$, $\alpha\in(0,1)$, $\beta>1$, $\eta>0$, $0\leq c < 1$}
\begin{algorithmic}[1]
\STATE{$\Delta t \gets (\Delta t)_t$}
\WHILE{$\max(LTE) < \eta$ and  $f(\x(t+\Delta t))< f(\x(t))+c\Delta t\nabla f(\xt)$}
\STATE{$\Delta t = \beta \Delta t$}
\ENDWHILE
\WHILE{$\max(LTE) \geq \eta$ or  $f(\x(t+\Delta t))\geq f(\x(t))+c\Delta t\nabla f(\xt)$}
\STATE{$\Delta t = \alpha \Delta t$}
\ENDWHILE
\RETURN $\Delta t$ 
\end{algorithmic}
\end{algorithm}

The largest step size that satisfies both the accuracy and stability checks can be found with a backtracking line search as shown in Algorithm \ref{time step algo}. The time step is initialized to $(\Delta t)_0$. The first while statement scales up $\Delta t$ until it violates the LTE condition \eqref{LTE cond} or the stability condition \eqref{conv cond}. The second while statement diminishes the time step until both conditions are satisfied. This algorithm is generic to any integration method, meaning that the user can customize the integration method to the application.

Before execution of the EATSS routine, there are a few design choices that must be specified. First, the time step search can either be initialized according to \eqref{tfe-naive} at each iteration, or to the $\Delta t$ used at the last iteration. If the objective function is known to have a small Lipschitz parameter $L$, then the change in gradient will be small, and the latter option may be more efficient in the subsequent search routine. 

Next, the LTE bound $\eta$ must be specified. In circuit simulation, the LTE is typically kept on the order of $10^{-5}$ \cite{Nilsson2000-jy}, but in optimization we can trade off high accuracy for speed; we find in simulation that $\eta \in (0.001, 1)$ works well. In optimization literature the Armijo parameter is generally set to $10^{-4}$, which we find also works well in simulation.

It is important to note that there will exist a $\Delta t$ that will satisfy both the accuracy and stability conditions. Regardless of the choice of integration method, as $\Delta t \rightarrow 0$, the approximation error of the numerical integration approaches 0. At a small enough $\Delta t$, the truncation error will thus satisfying the accuracy criteria. Additionally, each time-step takes a step in the direction of the negative slope of the objective function. A small enough $\Delta t$ will yield a $f(\x(t+\Delta t))$ that is smaller than $f(\x(t))$ and thus satisfy the stability criteria. 



\subsection{Applying Numerical Integration Methods to Optimization}
To approximate a point in the discretized trajectory, $\x(t+\Delta t)$, the integral on the RHS of \eqref{eq:discretization_state} must be approximated using a numerical integration technique. The specific choice of implementation can provide varying effects in the optimization trajectory and is often a compromise between accuracy and complexity.



A practitioner must choose between an explicit or implicit numerical integration method. Implicit techniques approximate the RHS of \eqref{eq:discretization_state} using the current state, $\x(t+\Delta t)$. While implicit integration methods are often more accurate, they often require queries to  nonlinear solvers which are computationally expensive.  In comparison, explicit methods approximate quantities using previous iterations and are less accurate for a given $\Delta t$ but have faster computation time.

The choice of integration technique is often application dependent. For example, circuit simulators often use trapezoidal integration (an implicit method) to accurately simulate the transient response. For optimization applications, however, the priority is fast computation of the final steady state, rather than accurately tracing the whole transient trajectory. In the following subsections, we propose two explicit integration techniques adapted for use in optimization algorithms and to solve \eqref{eq:discretization_state} and \eqref{scaled gradient flow}.

\subsubsection{Forward Euler Integration}
\label{sec:FE_integration}
Forward Euler is a simple explicit integration method that approximates the integral in \eqref{eq:discretization_state} using the state at the previous time-step, $\Delta t$ as:
\begin{equation}
    \int_{t}^{t+\Delta t} \Zeta(\x(t))^{-1}  \nabla f(\x(t)) dt  \approx \Delta t \Zeta(\x(t))^{-1}  \nabla f(\x(t)).
\end{equation}

The Forward-Euler method generally requires a small time step to ensure a  stable integration approximation \cite{pillage1995electronic}. 

The LTE for Forward-Euler can be estimated as: \cite{pillage1995electronic}
\begin{align}
    LTE = 0.5 \Delta t | Z^{-1}(\x(t))\nabla f(\xt) - Z^{-1}(\x(t+\Delta t))\nabla f(\x(t+\Delta t))\vert.
\end{align}
In the context of gradient descent, the LTE for FE is a measure of how fast the gradient changes over a single time step. By maintaining $LTE\leq\eta$, the upper bound for the step size adaptively changes to prevent the gradient from changing too quickly. This LTE control combined with checking for monotonic improvement in the objective function allow FE to be implemented accurately and stably.

In addition to LTE control, another circuit simulation idea that can be applied to discretization is the initial choice of time step $(\Delta t)_0$. The time step search method in Algorithm \ref{time step algo} scales an initial time step until the accuracy and stability conditions are satisfied; clearly a better initialization will reduce the time spent searching. One choice for a stable FE time-step is defined as: \cite{rohrer1981passivity}


\begin{align}
    \Delta t_{FE} &\leq 2 \frac{\xt^{\top}\nabla f(\x(t))}{ \sum_{i=1}^n Z_{ii}(\xt)^{-1} I_{c_i}^2(\xt)}. \label{tfe-naive}
\end{align}

This time-step ensures the linearized circuit at each time-point is passive, which as described in \cite{rohrer1981passivity}, equates to a numerically stable FE integration that avoids divergence. This bound can often be conservative, leading to smaller $\Delta t$ than desired.
For the application of optimization, accurate computation of intermediate points on the trajectory is less important than reaching the final steady-state. 
In  Algorithm \ref{time step algo}, \eqref{tfe-naive} can be used as an initial estimate for $\Delta t$ which we then maximize subject to the accuracy and stability conditions. Depending on choice of implementation, \eqref{tfe-naive} can be used as every initialization, or it can be used only in the first iteration with subsequent iterations initialized to the last used $\Delta t$. 




\subsubsection{Runge-Kutta Integration}
Another common numerical integration used in ODE methods is the explicit Runge-Kutta integration \cite{butcher1993rungeestimating}, which uses a linear combination of previous time-steps to estimate the integral in \eqref{eq:discretization_state}. These methods work by estimating the slope between $\xt$ and $\x(t+\Delta t)$ using a series of midpoints. In certain applications, this provides better accuracy of the numerical integration step than FE while still being cheap to compute. A 4th-order Runge-Kutta (RK4) integration for the scaled gradient flow is shown below:
\begin{equation}
    \int_{t}^{t+\Delta t} -Z(\x(t))^{-1}f(\x(t)) dt  \approx \frac{\Delta t}{6}(k_1 + k_2 + k_3 + k_4),
\end{equation}
where the values of $k_1, k_2, k_3$ and $k_4$ are evaluated at each time-step as:

\begin{equation*}
    \begin{aligned}[c]
    k_1 &= -Z(\xt)^{-1}\nabla f(\x(t))\\
    k_3 &= -Z\left(\x(t)+\frac{\Delta t}{2}k_2\right)^{-1}\nabla f\left(\x(t)+\frac{\Delta t}{2}k_2\right)
    \end{aligned}
    \hspace{0.5in}
    \begin{aligned}[c]
        k_2 &= -Z\left(\x(t)+\frac{\Delta t}{2}k_1\right)^{-1}\nabla f\left(\x(t)+\frac{\Delta t}{2}k_1\right)\\
        k_4 &= -Z\left(\x(t)+\Delta t k_3\right)^{-1}\nabla f\left(\x(t)+\Delta t k_3\right)
    \end{aligned}
\end{equation*}

The LTE of RK4 method for the scaled-gradient flow is estimated as:
\begin{align*}
    LTE &= \frac{1}{30} (10\x(t-3\Delta t) + 9\x(t-2\Delta t) -18\x(t-\Delta t)-\x(t)) \\
    &\:\:\:+\frac{1}{10}Z(x(t-3\Delta t))\Delta t(\nabla f(\x(t-3\Delta t)) \\
    &\:\:\:+ 6Z(x(t-2\Delta t))\nabla f(\x(t-2\Delta t)) +3Z(x(t-\Delta t))\nabla f(\x(t-\Delta t)))
\end{align*}

The derivation of the estimate above is provided in \cite{butcher1993rungeestimating}.
In the context of optimization, the Runge-Kutta explicit integration method uses past state information to estimate the current state. This is similar to momentum-based techniques, in which previous state and gradient information contribute to computing the next state.

In many applications, RK4 offers a considerable speed-up in comparison to FE by allowing larger time-steps that the accuracy and stability conditions. However, RK4 may not be suited if the objective function has a large Lipschitz parameter $L$, which correlates to VCCS elements that may change drastically over a short time-span. In these applications, single-step integration methods such as FE are often preferred over multi-step integration methods such as RK4, since the series of midpoints in RK4 may create an inaccurate integration approximation \cite{multistep_integration}.


%% file: examples.tex
\section{Examples}
\label{sec:examples}

The EC provides a structured methodology to incorporate circuit domain knowledge into optimization. The three mains steps are: building an equivalent circuit model, designing a capacitance/scaling matrix, and solving for the steady state. In this section we show how the scaled gradient flow and EC model can be used to rederive standard optimization techniques, thus linking our work to existing results. In Sections \ref{sec:control} and \ref{sec:discretization} we use the EC to improve the control and discretization steps and propose new optimization algorithms in Section \ref{sec:algo}.

\subsection{Gradient Descent with Constant Step Size}
Gradient descent is a commonly used optimization technique derived from the Taylor approximation of the objective function. From an initial state $\x_0$, the states are updated iteratively according to, 
\begin{equation}
    \x_{k+1}=\x_{k}-\alpha \nabla f(\x_{k}), \label{eq:gradient descent discrete}
\end{equation}
where $\alpha>0$ is a fixed step size.

\begin{figure}
    \centering
    \includegraphics[width=0.3\columnwidth]{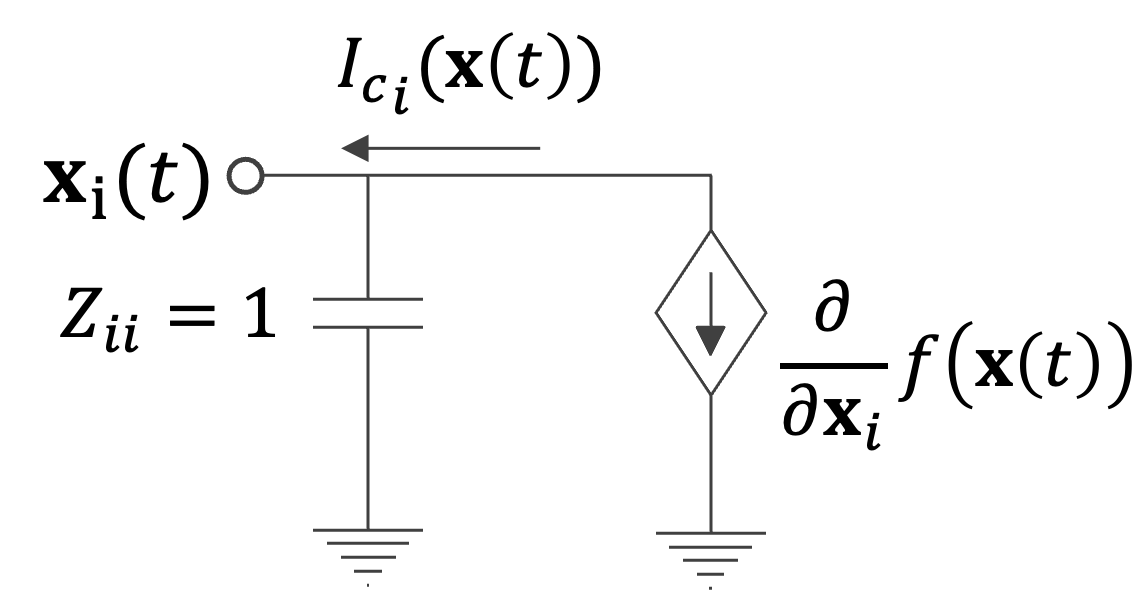}
    \caption{EC with unit capacitance}
    \label{fig:ex-unitcap}
\end{figure}

Note that \eqref{eq:gradient descent discrete} is equivalent to the FE discretization of,
\begin{equation}
    \dot{\x}(t) = -\nabla f(\x(t)), \label{eq:gradient descent cont}
\end{equation}
 with $\Delta t \triangleq \alpha$. Define $Z(\xt)^{-1} \triangleq \mathbf{I}_{n}$ to see that \eqref{eq:gradient descent cont} is in the same form as the standard form scaled gradient flow \eqref{scaled gradient flow}. The EC for this method thus has unit capacitances and is shown in Figure \ref{fig:ex-unitcap}.

\subsection{Gradient Descent with Varying Step Size}
Consider now the following modification to \eqref{eq:gradient descent discrete}:
\begin{equation}
    \x_{k+1}=\x_k -\beta_k \nabla f(\x_k), \label{eq:gd-vary}
\end{equation}
where $\beta_k$ is a variable step size at each iteration. For example, $\beta_k$ may be inversely proportional to $k$. 

To replicate the varying step size in ECCO, we can set $\Delta t$ to s fixed value and define the capacitor values to be $Z_{ii}(\x_k)^{-1}\triangleq \beta_k/\Delta t$. The FE discretization of this scaled gradient flow yields \eqref{eq:gd-vary}. 

\subsection{Heavy Ball Method}

The heavy ball method is a momentum method which uses two previous iterates to form the next guess: 
\begin{equation}
    \x_{k+1} = \x_k-\alpha_k\nabla f(\x_k)+\beta_k(\x_k-\x_{k-1}). \label{heavy ball}
\end{equation}

Using a different choice of discretization, we can replicate the heavy ball method in the ECCO framework. The EC can be defined with a unit capacitance as in Figure \ref{fig:ex-unitcap} and the uncontrolled gradient flow $\dxdt(t) = -\nabla f(\x(t))$. The heavy ball update rule \eqref{heavy ball} can be found by discretizing the flow with a second order Adam-Bashford integration method \cite{vanintroduction}. In general, this technique estimates the next point in a trajectory by interpolating two previous states:
\begin{equation}
    \x(t+\Delta t)=\x(t)-\Delta t(k_1 \nabla f(\x(t)) + k_2 \nabla f(\x(t-\Delta t))), \label{adambash}
\end{equation}
where $k_1$ and $k_2$ are scalar quantities. To match \eqref{heavy ball} and \eqref{adambash}, first note,
\begin{equation}
    -\nabla f(\x(t-\Delta t))=\dxdt(t-\Delta t)\approx (\Delta t)^{-1}(\x(t) - \x(t-\Delta t)),
\end{equation}
and then define $\Delta tk_1 = \alpha_k$ and $k_2 = \beta_k$ and ensure the initial conditions match.

\subsection{Nestrov's Accelerated Method}

Nestrov's accelerated method also uses momentum to mimic acceleration in physical systems. The updates are formed as,
\begin{align}
    \mathbf{y}_{k+1} &= \x_k+\beta_k(\x_k-\x_{k-1}),\\
    \x_{k+1} &=\mathbf{y}_{k+1}-\alpha_{k+1}\nabla f(\mathbf{y}_{k+1}),\\
    &\Rightarrow \x_{k+1} = \x_k+\beta_k(\x_k-\x_{k-1}) -\alpha_{k+1}\nabla f(\x_k+\beta_k(\x_k-\x_{k-1})).
    \label{eq:nestrovs}
\end{align}

In the ECCO framework, we can again define gradient flow $\dxdt_i(t) = -\nabla f(\x(t))$ where $Z=\mathbf{I}_n$, i.e. unit capacitance. Nestrov's acceleration can be shown to be a part of the family of 2nd-order Runge-Kutta discretization of the ODE \cite{dormand1980family}, which is implemented as:
\begin{align}
    &k_1 =  -\Delta t \nabla f(\x(t)) \\
    &k_2 =  - \frac{\Delta t}{2} \nabla f(\x(t) - c_1 k_1) \\
    &\x(t+\Delta t) = \x(t) + c_1 k_1 +  c_2 k_2
\end{align}

Using a finite differences approximation for $k_1$:
\begin{align}
        k_1 = -\Delta t \nabla f(\x(t)) = \Delta t \dot{\x}(t)\approx \Delta t \frac{\x_{k} - \x_{k-1}}{\Delta t}= \x_{k} -\x_{k-1}.
\end{align}
 
The 2nd order Runge-Kutta can be written as:
\begin{equation}
    \x(t+\Delta t) = \x(t) + c_1 (\x_{k} -\x_{k-1}) - c_2 \frac{\Delta t}{2} \nabla f(\x(t) - c_1 (\x_k - \x_{k-1})).
\end{equation}

With the choice of $c_1$, $c_2$ and $\Delta t$ that satisfy:
\begin{align}
    c_1&=\beta_k \\ 
    \Delta t c_2 &= 2 \alpha_{k+1},
\end{align}
we have achieved an iterative update that is identical to Nestrov's acceleration.

%% file: Algo.tex
\section{Algorithm Design}
\label{sec:algo}
In the previous sections, we have outlined three modular steps and their implementations in the ECCO workflow to incorporate circuit knowledge into building optimization methods. In this section, we demonstrate end-to-end implementations of new optimization algorithms derived from the ECCO workflow using the specific methods from the previous sections.

\subsection{Implementing All-Purpose EC Model}

After constructing the unmodified EC model of the scaled-gradient flow (Figure \ref{fig:ECorig}), a user can customize the capacitance matrix and numerical method to achieve varying effects in the optimization process. A user must decide to implement either the second order \eqref{z true} or the first order control scheme \eqref{z approx} to model the capacitors that is evaluated at each time-point. If the Hessian can be easily computed then the second order method is preferable because no error will be incurred from the trajectory approximation. For applications in which the Hessian is unavailable or inconvenient, e.g. training neural networks, then the speed of the first order control will likely out-weigh the effects of the approximation error. 
 We find that normalizing the diagonal elements of $Z_{ii}(\xt)^{-1}$ works well in practice as the step size will dictate the magnitude of the change in $\xt$ per iteration. 

Each of the suggested algorithms will use the EATSS routine in Algorithm \ref{time step algo} to choose time steps for the specific numerical integration process; see Section \ref{sec:discretization} for details on implementation and hyperparameters.

Algorithm \ref{ecco fe} describes an implementation using FE integration where step sizes are selected by the EATSS routine. As discussed in Section \ref{sec:FE_integration}, an initial time step can be pre-computed for FE via \eqref{tfe-naive} as shown in line \ref{ecco fe initial t}. In the main loop, the FE update itself is shown on line \ref{ecco fe disc}.

Algorithm \ref{ecco rk4} shows a similar routine but with the RK4 routine. In this algorithm, the initial time-step ($\Delta t_0$) is selected as the time-step of the last iteration. After calculating the control and time step size, the intermediate information for RK4 is calculated on lines \ref{ecco rk4 begin}-\ref{ecco rk4 middle} and the update is computed on \ref{ecco rk4 end}.

\begin{algorithm}
\caption{ECCO+FE+EATSS}
\label{ecco fe}
\textbf{Input: } $\x(0)$, $\alpha\in (0,1)$, $\beta>1$, $\eta>0$, $\delta>0$, $c\in [0,1)$
\begin{algorithmic}[1]
\STATE{$t\gets 0$}
\STATE{$(\Delta t)_0$ according to \eqref{tfe-naive}}\label{ecco fe initial t}
\STATE{\textbf{do while }$\|f(\x(t-\Delta t))-f(\x(t))\|> \epsilon$}
\STATE{\hspace*{\algorithmicindent}Choose $Z_{ii}(\x(t))^{-1}$ according to \eqref{z true} or \eqref{z approx}}
\STATE{\hspace*{\algorithmicindent} $\Delta t \gets$ EATSS($f(\cdot)$,  $\nabla f(\cdot)$, $\x(t)$, $Z(\xt)$, $(\Delta t)_t$, $\alpha$, $\beta$, $\eta$, $c$)}
\STATE{\hspace*{\algorithmicindent}$\x(t+\Delta t)=\x(t)-\Delta t Z(\xt)^{-1}\dfdxt$}\label{ecco fe disc}
\STATE{\hspace*{\algorithmicindent}Take step $t\gets t+\Delta t$}
\STATE{\hspace*{\algorithmicindent}$(\Delta t)_t\gets \Delta t$}
\RETURN $\x$
\end{algorithmic}
\end{algorithm}

\begin{algorithm}
\caption{ECCO+RK4+EATSS}
\label{ecco rk4}
\textbf{Input: } $\x(0)$, $\alpha\in (0,1)$, $\beta>1$, $\eta>0$, $\delta>0$, $c\in [0,1)$
\begin{algorithmic}[1]
\STATE{$t\gets 0$}
\STATE{\textbf{do while }$\|f(\x(t-\Delta t))-f(\x(t))\|> \epsilon$}
\STATE{\hspace*{\algorithmicindent}Choose $Z_{ii}^{-1}(\x(t))$ according to \eqref{z true} or \eqref{z approx}}
\STATE{\hspace*{\algorithmicindent} $\Delta t \gets$ EATSS($f(\cdot)$,  $\nabla f(\cdot)$, $\x(t)$, $Z(\xt)$, $(\Delta t)_t$, $\alpha\in(0,1)$, $\beta>1$, $\eta>0$)}

\STATE{\hspace*{\algorithmicindent}$k_1 = -Z(\x(t))^{-1} \nabla f(\xt)$}\label{ecco rk4 begin}
\STATE{\hspace*{\algorithmicindent}$k_2 = -Z(\x(t) + \frac{\Delta t}{2} k_1)^{-1} \nabla f\left(\x(t)+\frac{1}{2}\Delta tk_1\right)$}
\STATE{\hspace*{\algorithmicindent}$k_3 = -Z(\x(t)+\frac{\Delta t}{2} k_2)^{-1} \nabla f\left(\x(t)+\frac{1}{2}\Delta tk_2\right)$}
\STATE{\hspace*{\algorithmicindent}$k_4 = -Z(\x(t) + \Delta t k_3)^{-1} \nabla f\left(\x(t)+\Delta t k_3\right)$}\label{ecco rk4 middle}
\STATE{\hspace*{\algorithmicindent}$\x(t+\Delta t)=\x(t)+\frac{1}{6}\Delta t (k_1 + k_2 + k_3 + k_4)$}\label{ecco rk4 end}
\STATE{\hspace*{\algorithmicindent}Take step $t\gets t+\Delta t$}
\STATE{\hspace*{\algorithmicindent}$(\Delta t)_t\gets \Delta t$}
\RETURN $\x$
\end{algorithmic}
\end{algorithm}

A key benefit of these proposed algorithms is robustness to the hyperparameters. The convergence of the continuous-time trajectory is independent of $\delta$, so while its value can slow the convergence speed (large $\delta$ recovers general gradient flow, small $\delta$ will be ``ignored" with normalization) it cannot cause divergence. $(\Delta t)_0$, $\alpha$, and $\beta$ only affect the speed of the time search step.  $\eta$ and $c$ are the discretization hyperparameters that will affect the number of steps taken, as they set conditions on the shape of the variable's trajectory. The stability check ensures convergence. In simulation we find $\eta$ can be changed by orders of magnitude with little effect on the optimization trajectory. We find that $\eta=0.1$, $c=10^{-4}$, $\delta=1$, and normalizing $Z^{-1}$ at each iteration works well.

Additionally, the hyperparameters $\alpha$ and $\beta$ do not affect convergence as they only scale up or scale down the time-steps within the EATSS subroutine. In practice, a good selection of the initial time-step (either by \ref{ecco fe initial t} or using the last time-step) does not require multiple iterations of EATSS, making the entire algorithm robust to the choice of $\alpha$ and $\beta$.


\subsection{Homotopy-Driven Circuit Model}

We can propose another optimization method using the source-stepping approach described in Section \ref{sec:ec_model}. Controlled by the scalar homotopy factor $\gamma_{ss}$, the source-stepping homotopy method iteratively deforms the EC model as shown by Figure \ref{fig:ec_source_stepping}. In this example, we use the source stepping method with a unit capacitance matrix (i.e. $\Zeta=I$), and a FE numerical integration method. The entire algorithm is shown in Algorithm \ref{source_stepping_algo}. In this method, the outer loop (in line 4), increments the homotopy factor, while the inner loop solves the resulting deformed EC to its steady-state. Importantly, we assign the initial conditions of each homotopy problem as the solution of the previous homotopy problem (line 10). 


\begin{algorithm}
\caption{Source Stepping ECCO+FE+EATSS}
\label{source_stepping_algo}
\textbf{Input: } $\x(0)$, $\alpha\in (0,1)$, $\beta>1$, $\eta>0$, $\delta>0$, $c\in [0,1)$, $\Delta \gamma_{ss}\in (0,1)$
\begin{algorithmic}[1]
\STATE{$\gamma_{ss} \gets 0$}
\STATE{$I_{c_0} \gets I_c(\x(0))$}
\STATE{\textbf{do while} $\gamma_{ss} > 0$}
\STATE{\hspace*{\algorithmicindent}$t\gets 0$}
\STATE{\hspace*{\algorithmicindent}\textbf{do while }$\|f(\x(t-\Delta t))-f(\x(t))\|> \epsilon$}

\STATE{\hspace*{\algorithmicindent}\hspace*{\algorithmicindent} $\Delta t \gets$ EATSS($f(\cdot)$,  $\nabla f(\cdot)$, $\x(t)$, $Z(\xt)$, $(\Delta t)_t$, $\alpha\in(0,1)$, $\beta>1$, $\eta>0$)}
\STATE{\hspace*{\algorithmicindent}\hspace*{\algorithmicindent}$\x(t+\Delta t)=\x(t)-\Delta t (Z(\xt)^{-1}\dfdxt - \gamma_{ss}I_{c_0}$}\label{ecco fe disc}

\STATE{\hspace*{\algorithmicindent}\hspace*{\algorithmicindent}Take step $t\gets t+\Delta t$}
\STATE{\hspace*{\algorithmicindent}\hspace*{\algorithmicindent}$(\Delta t)_t\gets \Delta t$}
\STATE{\hspace*{\algorithmicindent}$\x(0)\gets \x$}
\STATE{\hspace*{\algorithmicindent}$\gamma_{ss} \gets \gamma_{ss}+\Delta \gamma_{ss}$}
\RETURN $\x$
\end{algorithmic}
\end{algorithm}

%% file: results.tex
\section{Experiments}
\label{sec:results}

In this section, we demonstrate using the algorithms presented Section \ref{sec:algo} to solve optimization problems. We first demonstrate the applicability of the ECCO algorithms on known test functions to evaluate performance to known critical points. The second example is a power systems application which demonstrates a difficult optimization function with a computable Hessian. This example also illustrates an application for which source stepping is especially useful due to the initial conditions. Finally, we demonstrate ECCO on a neural network example to classify MNIST data and perform an experiment on robustness to hyperparameters.

All ECCO experiments used the hyperparameters $\eta = 0.1$, $\delta = 1$, $\alpha = 0.9$, $\beta = 1.1$, $c=10^{-4}$. The EC-based optimization methods are compared against ADAM, gradient-descent with Armijo-line search, BFGS, and SR1. In each of the experiments, the comparison methods are optimally tuned for the fastest convergence. Convergence is defined as $\x_k$ such that $|f(\x_k)-f(\x_{k+1})|<1e-4$. 


\subsection{Optimization Test Functions}
We first demonstrate the applicability of ECCO optimization algorithms on optimization test functions with varying properties; for example, on multi-modal (non-convex) functions versus functions with flat valleys \cite{jamil2013literature}. 
 The tested scenarios were (1) the convex Booth function \cite{testfunc}, (2) the non-convex multi-modal Himmelblau function \cite{0070289212}, (3) the non-convex Rosenbrock function \cite{testfunc}, and (4) the nonconvex three-hump function \cite{testfunc}.
 
The functions were optimized using 
the following four EC-based optimization methods highlighted in Table \label{table:simulation_optimizers}.

\begin{center}
\begin{tabular}{||c c c c||} 
 \hline
 Optimization Name & EC Model & $\Zeta(\x)$ & Numerical Integration\\ [0.5ex] 
 \hline\hline
 1 & Original EC & Full Hessian & Forward-Euler \\ 
 \hline
 2 & Original EC & Approximate Hessian & Forward-Euler \\
 \hline
 3 & Original EC & Full Hessian & Runge-Kutta \\
 \hline
 4 &  Original EC & Approximate Hessian & Runge-Kutta \\
 \hline
\end{tabular}
\label{table:simulation_optimizers}
\end{center}

\begin{figure}
    \centering
    \includegraphics[width=.8\columnwidth]{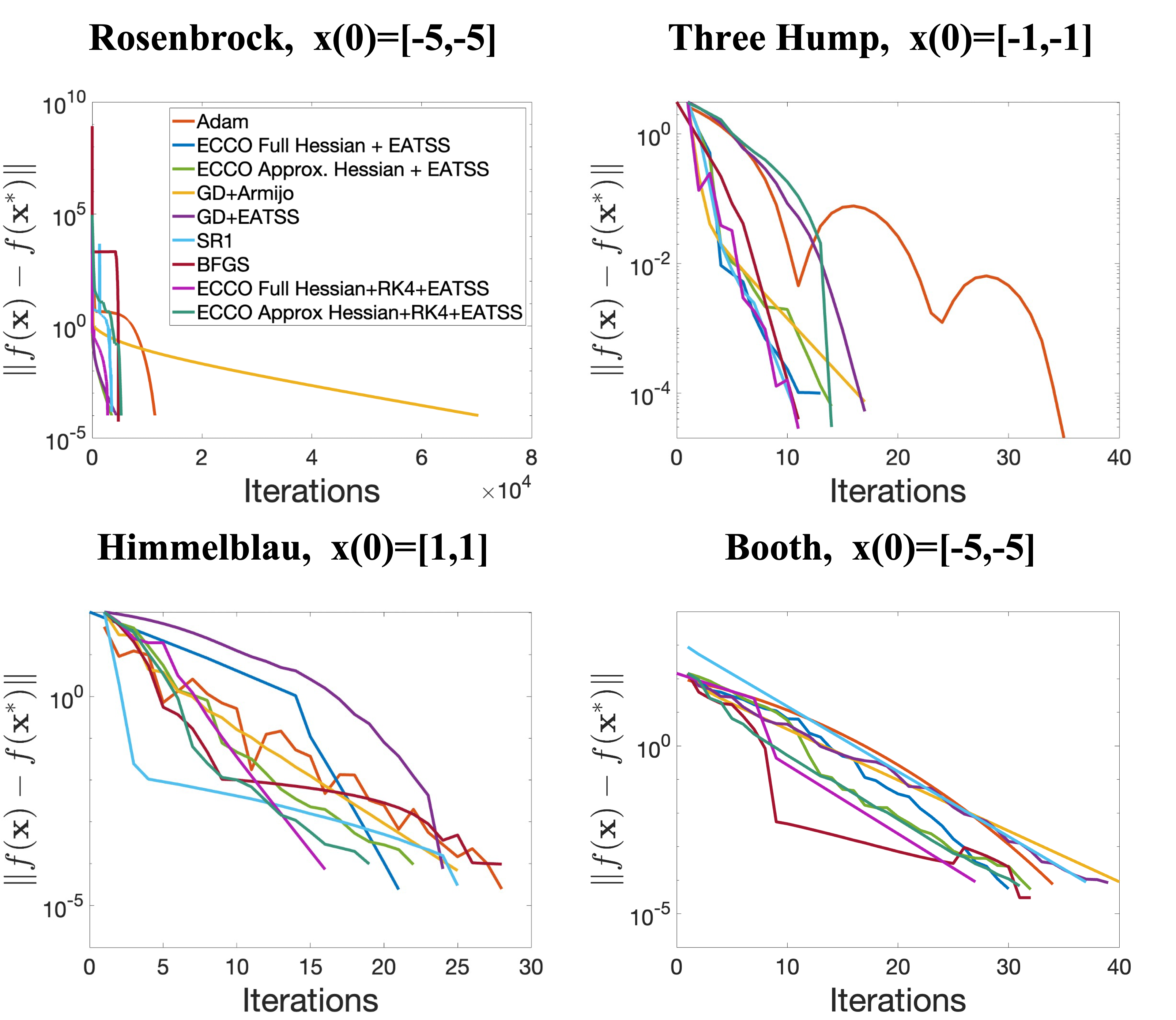}
    \caption{Comparison of Optimization Methods on Test Functions Using ECCO Optimization Methods from Table \ref{table:simulation_optimizers}}
    \label{fig:test_functions}
\end{figure}

The results are shown in Figure \ref{fig:test_functions} All four ECCO optimization methods achieve faster or on-par convergence compared to second-order methods such as BFGS and SR1, while outperforming first-order methods (Adam and gradient descent with Armijo line search). Unlike the comparison methods which are optimally tuned for fast convergence, the ECCO optimization methods did not require any hyperparameter tuning to achieve these convergence results. Note that the approximate $\hat{Z}$ \eqref{z approx} approximates the full Hessian $Z$ \eqref{z true}, as methods 1 and 2, as well as 3 and 4 in Table \ref{table:simulation_optimizers} achieve similar performances. Also, the Runge-Kutta also outperforms the Forward-Euler since it uses four previous time-steps for a momentum-like effect. 

\subsection{Optimizing Power Systems}
Next we apply ECCO to a power systems optimization problem to demonstrate its performance on studying the feasibility of a power grid by finding the smallest total constraint violations for the system. In planning applications for power systems, this helps identify areas where additional assets can be built, and in operations, this study finds if the dispatch for the power grid will be operational.

The optimization problem is $min_\x \frac{1}{2} \| f(\x) \| ^2$ where $f(\x)$ defines the power-flow constraints at each bus (defined in \cite{FOSTER2022108486}). The equation $f$ is a multi-modal, non-convex function and optimizes a 118-bus power grid network with 291 variables. 

Using the ECCO workflow, we can design optimization methods that are best suited to the application with the goal of reaching the lowest total violations $\|f(X)\|^2$. In this application, evaluating the Hessian is generally easy to compute due to the sparsity of the power grid. As a result, we optimize this objective using ECCO with a full Hessian \eqref{z true} and ECCO using source-stepping with a full Hessian. We also use a single-step FE integration since we expect a high Lipschitz parameter, $L$, for the optimization problem.  We optimize the 118-bus network for both planning and operation applications and demonstrate the benefit of both methods against other second-order methods including  BFGS \cite{fletcher2013practical}, SR1 \cite{conn1991convergence} and Newton-Raphson \cite{fletcher2013practical}.

\subsubsection{Planning Application}
In the planning scenario, the optimization problem is solved to identify additional assets that need to be built to satisfy a future grid scenario. Since this represents a potential grid topology in the future, the initial conditions are generally not close to the optimum (initial conditions are usually assumed to be flat start, where the voltages are 1 $p.u.$). In this scenario, we optimize a future scenario of the IEEE-118 bus network that has a higher loading factor of 1.5 with additional installed generation.

 The two methods from the ECCO workflow (ECCO with Full Hessian + FE + EATSS and ECCO with Source-Stepping), are both used to optimize this scenario as shown in Figure \ref{fig:power_systems_planning}. In comparison to other second-order methods, both methods achieve a lower total feasibility violations. However, source-stepping is not suited for this scenario because the initial conditions are not close to the optimal solution.


\subsubsection{Operation Application}
In the operation scenario, we solve the optimization problem during each dispatch to identify possible infeasibilities in the grid. The operation case is solved every 5-10 minutes and uses the previous dispatch as an initial condition. In general, there is limited change in the network between successive dispatches, and the initial conditions are close to the optimal solution. In this experiment, we optimize the 118-bus system for a operation scenario where a single line is tripped. We use the solution from the pre-fault conditions as an initial condition, which is generally close to the optimal solution.

\begin{figure}
\centering
\begin{minipage}{.5\columnwidth}
  \centering
  \includegraphics[width=.8\columnwidth]{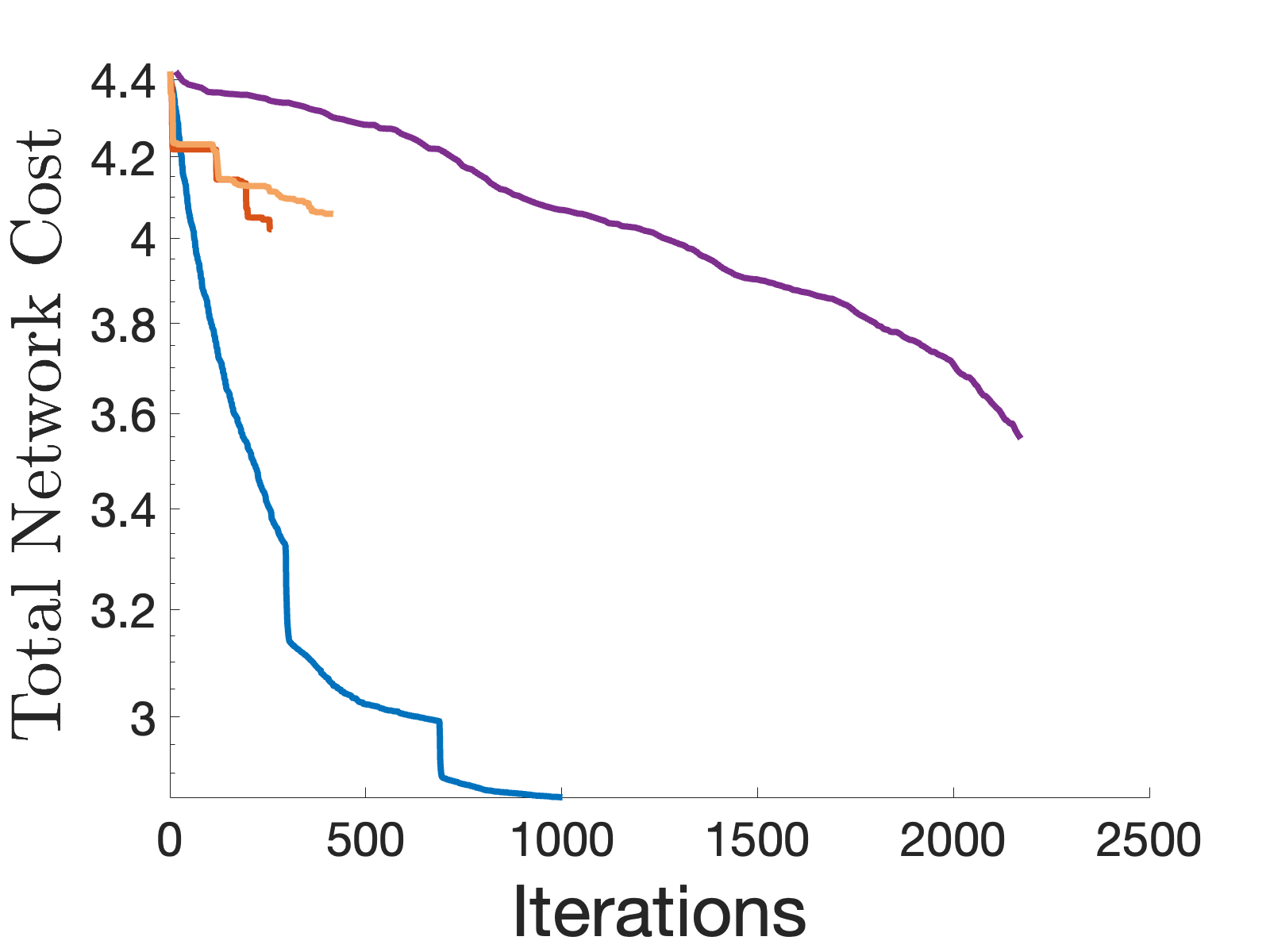}
  \caption{Optimizing the Planning Scenario for Stressed IEEE-118 Bus System with Flat Start Initial Conditions}
  \label{fig:power_systems_planning}
\end{minipage}%
\begin{minipage}{.5\columnwidth}
  \centering
  \includegraphics[width=.8\columnwidth]{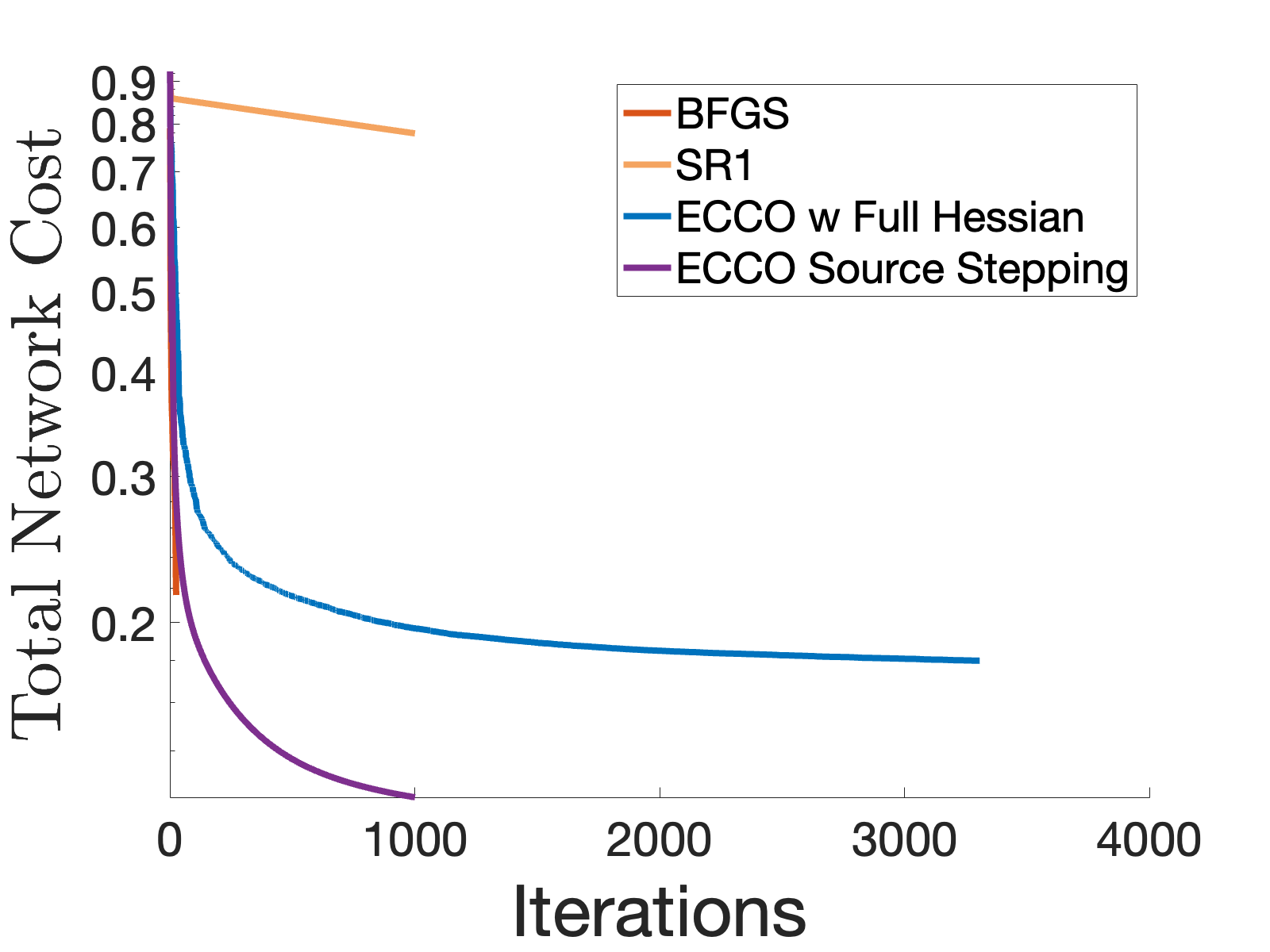}
  \caption{Optimizing the Operation Scenario for Stressed IEEE-118 Bus System with Initial Conditions from Prior Dispatch}
  \label{fig:power_systems_operations}
\end{minipage}
\end{figure}

As shown in Figure \ref{fig:power_systems_operations}, the source-stepping method provides faster convergence for the operation since the initial condition is close to the final solution. In comparison to other second-order methods, the source-stepping algorithm converges faster to a state with a smaller total violation.

This particular application demonstrates that we can design optimization methods using the ECCO workflow that leverage knowledge of initial conditions and the application to select a better method.

\subsection{Training a Neural Network}

We next design an optimization method using ECCO to train a neural network to classify MNIST data \cite{deng2012mnist}. See the Appendix \ref{more nn} for details on this experiment. The optimization is designed to provide fast convergence and robustness to hyperparameters. However, evaluating the Hessian in this application is generally costly, thereby suggesting we use a first-order method. Additionally, a high Lipshitz parameter indicates we use a FE numerical integration. With this goal, we choose an ECCO optimization method with the approximate Hessian control \eqref{z approx} and FE numerical integration with EATSS. The approximate Hessian control provides fast convergence to steady-state in continuous time, while EATSS provides robustness to the following hyperparameters  $\theta=[\alpha, \beta, \eta, \delta,c]$.

\begin{figure}
\centering
\begin{minipage}{.5\columnwidth}
  \centering
  \includegraphics[width=.8\columnwidth]{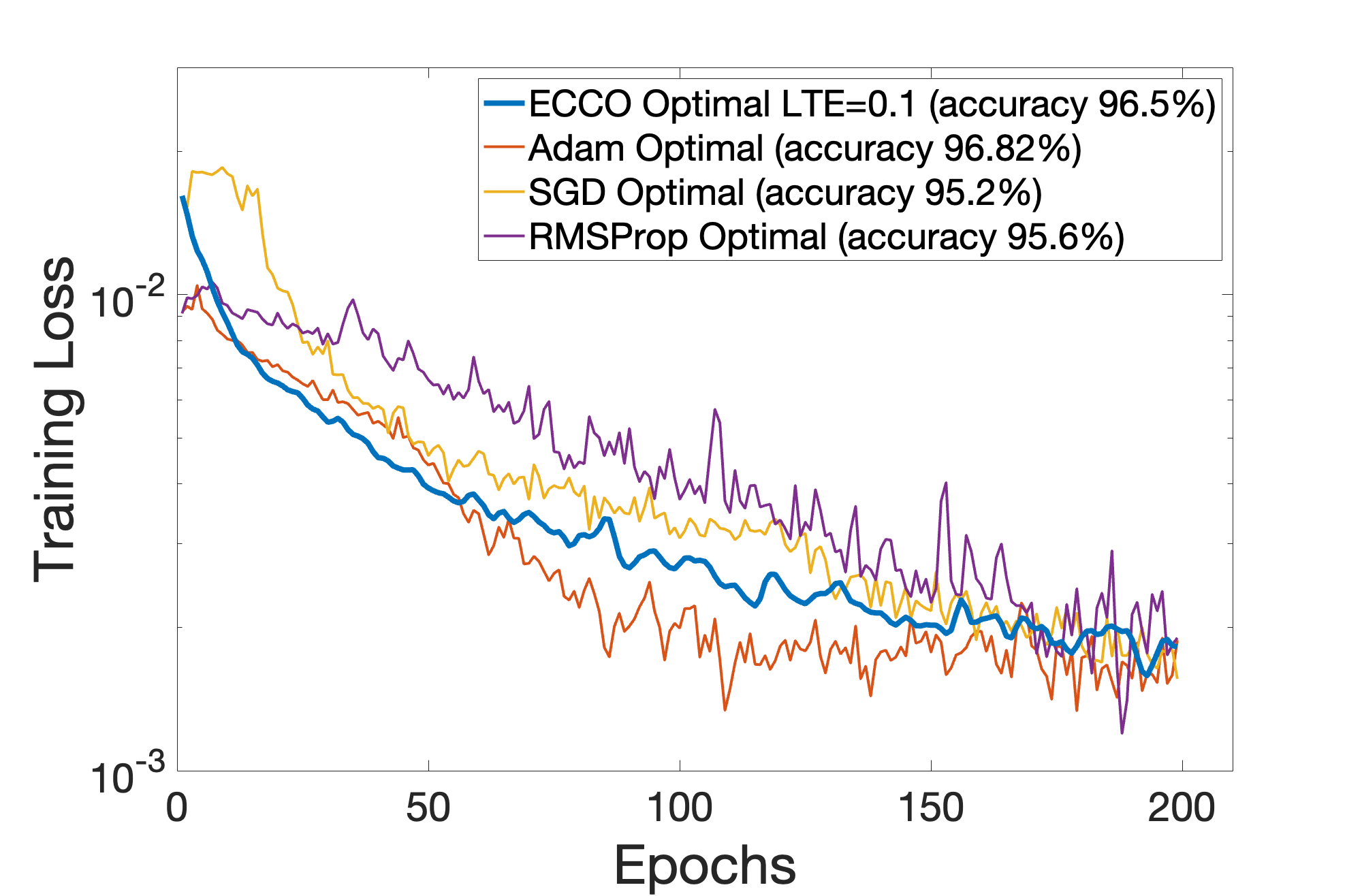}
    \caption{\small{Training of neural network for MNIST data using optimally tuned first order ECCO \eqref{z approx}, Adam, GD and RMSProp.}}
  \label{fig:supp_training_loss_mnist}
\end{minipage}%
\begin{minipage}{.5\columnwidth}
  \centering
  \includegraphics[width=.8\columnwidth]{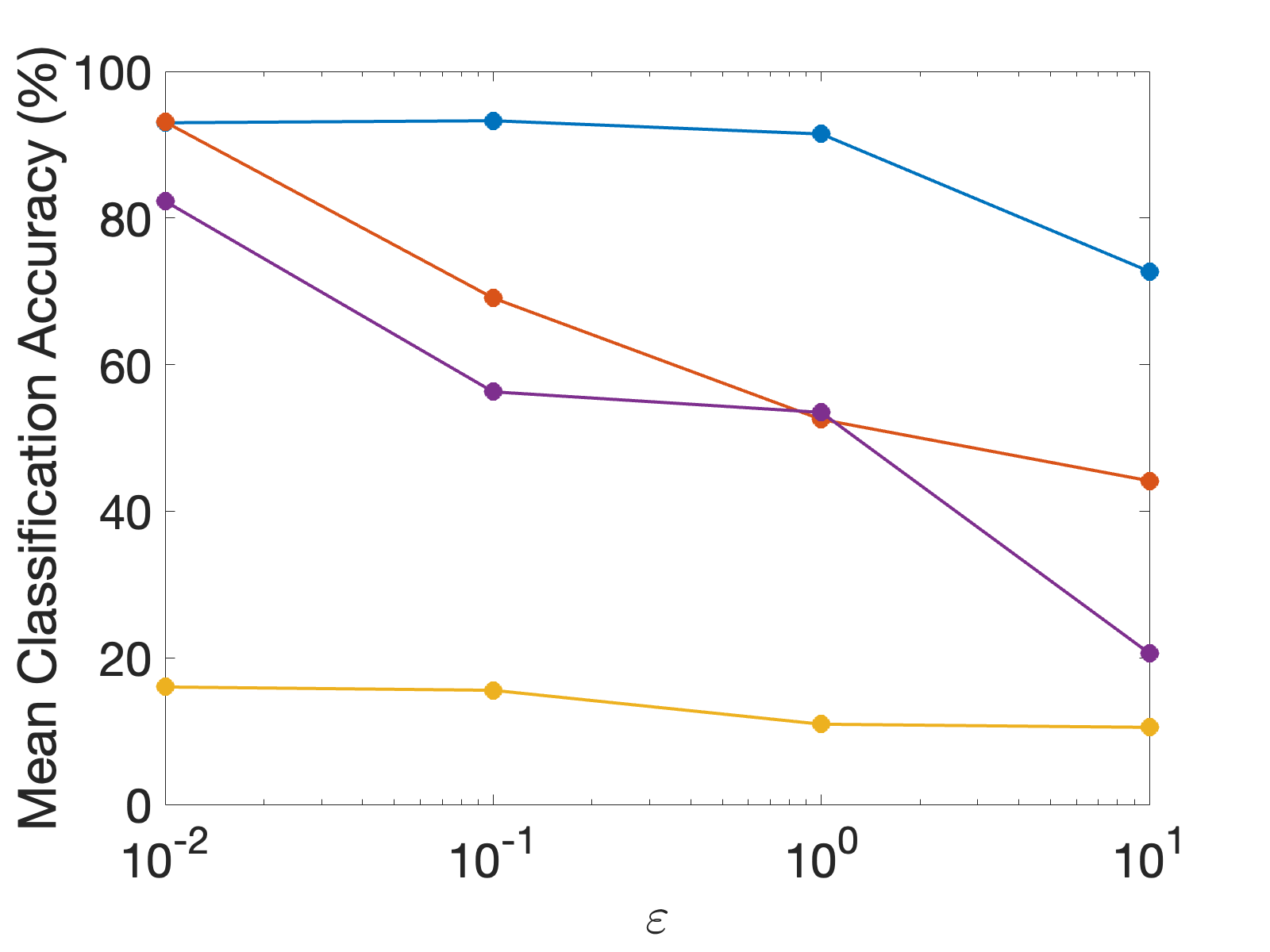}
  \caption{Mean Classification Accuracy ($\%$) of MNIST Data Using ADAM, SGD, RMSProp and ECCO (Approximate Hessian + FE) with Hyperparameters Randomly Selected from Ball $\varepsilon$.}
  \label{fig:ml_varepsilon_test}
\end{minipage}
\end{figure}

We compared first order ECCO \eqref{z approx} with FE+EATSS to other first-order methods used in training neural networks including stochastic gradient descent (SGD) with a fixed step size, Adam, and RMSProp. We perform a grid search to find optimal hyperparameters for the  the comparison methods. See Appendix \ref{hyper search} for details. Using a minibatch size of 1000, all four methods train the neural network for 200 epochs and achieve a similar test classification accuracy of $96\%$. Figure \ref{fig:supp_training_loss_mnist} shows the result of the training loss.


To demonstrate robustness of the ECCO optimization method, we perturb the optimally tuned hyperparameter vectors within a normalized ball of radius $\varepsilon$. The hyperparameter values were sampled from a uniform distribution as $\tilde\theta \sim U(\max\{{\theta}^*-\frac{\varepsilon}{{\theta}^*}, \underline{\theta}\},\min\{{\theta}^*+\frac{\varepsilon}{{\theta}^*},\bar{\theta} \}) $, where ${\theta}^*$ was the optimal hyperparameter values and $[\underline{\theta} ,\bar{\theta}]$ are the operating bounds of the hyperparameters. Note, the hyperparameters were selected from a radius $\varepsilon/\theta^*$, indicating that the size is relative to value of the optimal hyperparameters.

The network was then trained with the randomly selected hyperparameters, $\tilde{\theta}$ for 200 epochs. The test classification accuracy was recorded as we vary the normalized ball of radius, $\varepsilon$ by four magnitudes. The mean classification accuracy for different $\varepsilon$ was recorded in Figure \ref{fig:ml_varepsilon_test}. We find that while Adam, SGD, and RMSProp are disrupted by the hyperparameter perturbation, ECCO was less sensitive to changes in the hyperparameter. In fact, we found the hyperparameters of ECCO could be changed by three orders of magnitude and still yield classification accuracies over 70\%. 

These experiments demonstrate that ECCO can achieve convergence behavior on par with Adam without comparable computation and data requirements to select hyperparameters, i.e. cross-validation. This means ECCO is better than Adam at problem generalization and  distribution shift.

%% file: Conclusion.tex
\section{Conclusion}
\label{sec:conclusion}


In this paper we propose ECCO, a novel workflow for unconstrained optimization by modeling gradient flow as a physical process to design new optimization methods for unconstrained convex and nonconvex problems using knowledge from the physical dynamical systems. In this work we map gradient flow onto the domain of circuits by constructing an equivalent circuit that is mathematically identical to the gradient flow problem. This circuit provides natural, physics-based intuition to control the trajectory of the solution to a steady-state, which coincides with the critical point of the objective function. We use the equivalent circuit to propose new techniques to design the trajectory and step-size of the optimization trajectory based on circuit principles. We define a methodology to systematically incorporate circuit knowledge by 1) designing the underlying circuit model, 2) controlling the trajectory by designing nonlinear capacitance models that minimizing the total stored charge in the circuit and 3) developing fast numerical integration methods which ensure stability and accuracy. Using the ECCO workflow, we can design new optimization methods that demonstrate robust convergence of non-convex test functions without hyper-parameter tuning and is capable of outperforming state-of-the art optimization methods. By bridging the fields of optimization and circuit simulation, the framework enables a new class of methods of physics-empowered optimization and future innovation.


%% file: Appendix.tex
\section{Proof of Theorem 1}
\label{main proof}

\begin{align}
    \frac{d}{dt}f(\xt)&=\langle \dfdxt,\dxdt(t)\rangle \\
    &= -\dfdxt^{\top} Z(\xt)^{-1}\dfdxt\\
    &\leq -d_1\Vert \dfdxt\Vert^2 \label{bounded dfdt}
\end{align}
Where \eqref{bounded dfdt} holds by Assumption \ref{z def}. 

Eq. \eqref{bounded dfdt} implies that the objective function is non-increasing along $\xt$.

For $t>0$ consider $\int_0^t \Vert\dfdxt\Vert^2dt$. By \eqref{bounded dfdt},
\begin{align}
    \int_0^t\Vert\dfdxt\Vert^2dt&\leq\frac{1}{d_1}(f(\x(0))-f(\x(t)))\\
    &\leq \frac{2R}{d_1} \label{bounded int norm}
\end{align}
Where \eqref{bounded int norm} holds by Assumption \ref{a1}.

Then for all $t>0$,
\begin{equation}
     \int_0^t\Vert\dfdxt\Vert^2dt\leq \frac{2R}{d_1}<\infty
\end{equation}

For all $\x$, $\y\in \R^n$,
\begin{align}
    \big \vert\Vert \dfdx \Vert^2 - \Vert \dfdy \Vert^2 \big \vert &\leq \big \vert \Vert \dfdx \Vert + \Vert \dfdy \Vert \big \vert \ \big \vert \Vert \dfdx \Vert -\Vert \dfdy \Vert \big \vert\\
    &\leq 2B \big \vert \Vert \dfdx \Vert - \Vert \dfdy \Vert \big \vert\\
    &\leq 2B \Vert \dfdx - \dfdy\Vert\\
    &\leq 2BL \Vert \x-\y \Vert
\end{align}

Therefore we can conclude that $\Vert \dfdx\Vert^2:\R^n\to\R$ is Lipschitz and hence uniformly continuous.

By Assumptions \ref{a4} and \ref{z def}, $\Vert \dxdt(t)\Vert$ is bounded and so $\xt$ is uniformly continuous in $t$.

Therefore, the composition $\Vert \dfdxt \Vert^2: \R_+\to\R$ is uniformly continuous in $t$. 

Since $\int_0^t\Vert\dfdxt\Vert^2dt<\infty$ and $\Vert\dfdxt\Vert^2$ is a uniformly continuous function of $t$, we conclude that $\lim_{t\to\infty}\Vert\dfdxt\Vert=0$.

\section{Proof of Lemma \ref{bounded z true}}\label{zbounded}

\begin{proof}
Note that the boundedness condition in Assumption \ref{z def} will hold if $0< \sum_{i=1}^n Z_{ii}(\xt)^{-1}<\infty$.

Lower bound:
\begin{align}
    Z_{ii}^{-1}(\xt) = \max\{\delta^{-1} [G(\xt) \nabla^2 f(\xt)\dfdxt]_i,1\}\geq 1>0.
\end{align}
Upper bound:
\begin{align}
     \sum_{i=1}^n Z_{ii}(\xt)^{-1} &\leq \mathbf{1}^{\top}\Big(\delta^{-1}G(\xt)\nabla^2f(\xt)\dfdxt\Big)+\sum_{i=1}^n 1\\
     &=\delta^{-1}\mathbf{1}^{\top}G(\xt)\nabla^2f(\xt)\dfdx +n\\
     &=\delta^{-1}\dfdxt^{\top}\nabla^2f(\xt)\dfdx +n\\
     &\leq \delta^{-1}\lambda_{\max}\big(\nabla^2 f(\xt)\big)\Vert \dfdxt\Vert^2+n\label{bound w lambda}\\
     &\leq \delta^{-1}\lambda_{\max}\big(\nabla^2 f(\xt)\big)B^2+n<\infty
\end{align}
Where $\lambda_{\max}\big(\nabla^2 f(\xt)\big)$ exists and is finite due to the existence of symmetric $\nabla^2 f$ \ref{a1}.
\end{proof}

\section{Derivation of \eqref{z approx}}\label{z approx deriv}
Recall the approximation of $\frac{d}{dt}\dfdxt$ as in \eqref{ahat} and then use Lemma \ref{identity1} to find,
\begin{align}
    \ahat(\xt)&= \frac{\dfdxt - \nabla f(\x (t-\Delta t))}{\Delta t}\nonumber\\
    &\approx \frac{d}{dt}\dfdxt\\
    &=-\nabla^2 f(\xt) Z(\xt)^{-1} \dfdxt.
\end{align} 
Pre-multiply by $\delta^{-1}G(\x(t))$:
\begin{align}
    \delta^{-1}G(\xt)\ahat(\xt)\approx -\delta^{-1}G(\xt) \nabla^2 f(\xt) Z(\xt)^{-1} \dfdxt. \label{deriv 1}
\end{align}
As $Z$ is a diagonal matrix, the expression on the RHS of \eqref{deriv 1} can be rewritten as,
\begin{align}
    &-\delta^{-1}G(\xt) \nabla^2 f(\xt) Z(\xt)^{-1} \dfdxt\nonumber\\
    &=- \delta^{-1}G(\xt)\nabla^2 f(\xt)\begin{bmatrix}
    Z_{11}(\xt)^{-1}\frac{\partial f(\xt)}{\partial \x_1}\\ \vdots\\ Z_{nn}(\xt)^{-1}\frac{\partial f(\xt)}{\partial \x_n}
    \end{bmatrix} \label{deriv 2}\\
     &=-\delta^{-1}G(\xt)\nabla^2 f(\xt)\begin{bmatrix}
        \frac{\partial f(\xt)}{\partial x_1} & 0 &\dots & 0\\
        0  & \frac{\partial f(\xt)}{\partial x_2} & \ddots & 0\\
        \vdots & \vdots & \ddots & \vdots \\
        0 & 0 & \dots & \frac{\partial f(\xt)}{\partial x_n}
    \end{bmatrix}\begin{bmatrix}
        Z_{11}(\xt)^{-1}\\ \vdots \\ Z_{nn}(\xt)^{-1}
    \end{bmatrix}\\
    &=-\delta^{-1}G(\xt)\nabla^2 f(\xt)G(\xt)\z(\xt). \label{deriv 4}
\end{align}
In comparison, the objective is to compute:
\begin{align}
    \z(\xt)&= \delta^{-1}G(\xt) \nabla^2 f(\xt) \dfdxt \label{deriv 2}\\
    &= \delta^{-1}G(\xt) \nabla^2 f(\xt)G(\xt) \mathbf{1}\label{deriv 3}.
\end{align}
Where \eqref{deriv 2} follows from \eqref{final z} and \eqref{deriv 3} holds because by definition $G(\xt)_{ii} = [\dfdxt]_i$.

Define the shorthand notation for the expressions $\Ab\triangleq G(\xt)\nabla^2 f(\xt)G(\xt)$ and $\z\triangleq \z(\xt)$. Then \eqref{deriv 4} becomes $-\delta^{-1}\Ab\z$ and similarly \eqref{deriv 3} becomes $\delta^{-1}\Ab\mathbf{1}$. From \eqref{deriv 1} and \eqref{deriv 4},
\begin{equation}
    -\delta^{-1}\Ab \z \approx \delta^{-1}G(\xt)\ahat(\xt) \label{deriv 6}
\end{equation}
Use the definition of $\z$ in \eqref{deriv 3} and substitute into \eqref{deriv 6}:
\begin{align}
    \z &= \delta^{-1}\Ab \mathbf{1} \Rightarrow -\delta^{-1}\Ab \z = -\delta^{-2}\Ab\Ab\mathbf{1} \approx   \delta^{-1}G(\xt)\ahat(\xt)
\end{align}
Pre-multiply by $-\mathbf{1}^{\top}$ and note that $\Ab = \Ab^{\top}$ :
\begin{align}
    (\delta^{-1}\Ab\mathbf{1})^{\top}(\delta^{-1}\Ab\mathbf{1}) &\approx -\delta^{-1}\mathbf{1}^{\top}G(\xt)\ahat(\xt)\\
    \z^{\top}\z&\approx-\delta^{-1}\mathbf{1}^{\top}G(\xt)\ahat(\xt) \label{deriv 5}
\end{align}

One possible solution to \eqref{deriv 5} is thus,
\begin{align}
    &\z_i^2(\xt)\triangleq -[\delta^{-1}G(\xt)\ahat(\xt)]_i\\
    &\Rightarrow \z_i(\xt) = \sqrt{\max\{-[\delta^{-1}G(\xt)\ahat(\xt)]_i, 0\}}.
\end{align}

Similarly to the derivation of \eqref{z true}, we truncate $\z_i(\xt)\geq 1$ to finish the derivation.

\section{Proof of Lemma \ref{bounded z approx}} \label{zhat bounded}

\begin{proof}
Note that the boundedness condition in Assumption \ref{z def} will hold if $0< \sum_{i=1}^n \hat{Z}_{ii}(\xt)^{-1}<\infty$.

Lower bound:
\begin{align}
    \hat{Z}_{ii}(\xt)^{-1}=\max\{\sqrt{-\delta^{-1}[G(\xt)\ahat(\xt)]_i},1\}\geq 1 >0.
\end{align}
Upper bound:
\begin{align}
    \hat{Z}_{ii}(\xt)^{-1} &\leq \max\{\sqrt{-\delta^{-1}[G(\xt)\ahat(\xt)]_i},1\}^2\\
    \sum_{i=1}^n \hat{Z}_{ii}(\xt)^{-1} &\leq  \sum_{i=1}^n -\delta^{-1}[G(\xt)\ahat(\xt)]_i+\sum_{i=1}^n 1\\
    &= \delta^{-1}(\Delta t)^{-1}\Big(\dfdxt^{\top} \nabla f(\x(t-\Delta t)) -\dfdxt^{\top}\dfdxt \Big)+n\\
    &=\delta^{-1}(\Delta t)^{-1}\dfdxt^{\top}\Big( \nabla f(\x(t-\Delta t)) -\dfdxt \Big)+n\\
    &\leq\delta^{-1}(\Delta t)^{-1}\left\vert\dfdxt^{\top}\Big(  \nabla f(\x(t-\Delta t)) -\dfdxt \Big)\right\vert+n\\
    &\leq \delta^{-1}(\Delta t)^{-1}\Vert \dfdxt\Vert \Vert \nabla f(\x(t-\Delta t)) -\dfdxt\Vert+n\\
    &\leq \delta^{-1}(\Delta t)^{-1} BL\Vert \x(t-\Delta t)-\x(t)\Vert+n\\
    &\leq \delta^{-1}(\Delta t)^{-1} BL\max_{\x^*\in S}\Vert \x(0)-\x^*\Vert+n<  \infty
\end{align}
Where the last line holds because $\Vert \x(0)-\x^*\Vert$ is bounded as $\x(0)$ is given to be finite and $\x^*$ is finite due to \ref{a2}.
\end{proof}

\section{Computation Complexity}
\subsection{Proof of Lemma \ref{z true complexity}} \label{complexity full}
Given the gradient and Hessian, the per-iteration computation complexity to evaluate $Z$ according to \eqref{z true} is $\mathcal{O}(n^2)$. This can be verified as each element $i$ requires $\mathcal{O}(n)$ computations:
\begin{equation}
    Z_{ii}^{-1}(\xt) = \max\left\{\delta^{-1} \frac{\partial f(\xt)}{\partial \x_i(t)}\sum_{j=1}^n \frac{\partial^2 f(\xt)}{\partial \x_i \partial \x_j}\frac{\partial f(\xt)}{\partial \x_j}, 1\right\}
\end{equation}

\subsection{Proof of Lemma \ref{z approx complexity}} \label{complexity approx}
Given the current and previous gradient, the per-iteration computation complexity to evaluate $\widehat{Z}$ according to \eqref{z approx} is $\mathcal{O}(n)$. This can be verified as each element $i$ requires $\mathcal{O}(1)$ computations:
\begin{equation}
    Z_{ii}^{-1}(\xt) = \max\left\{-\delta^{-1}\frac{1}{\Delta t} \frac{\partial f(\xt)}{\partial \x_i(t)} \left(\frac{\partial f(\xt)}{\partial \x_i}-\frac{\partial f(\x(t-\Delta t))}{\partial \x_i}\right), 1\right\}
\end{equation}

\section{More Details on Test Functions Experiment}\label{all toy graphs}

\subsection{All Graphs}
\begin{figure}[h]
\centering

\begin{subfigure}{0.3\textwidth}
    \includegraphics[width=\textwidth]{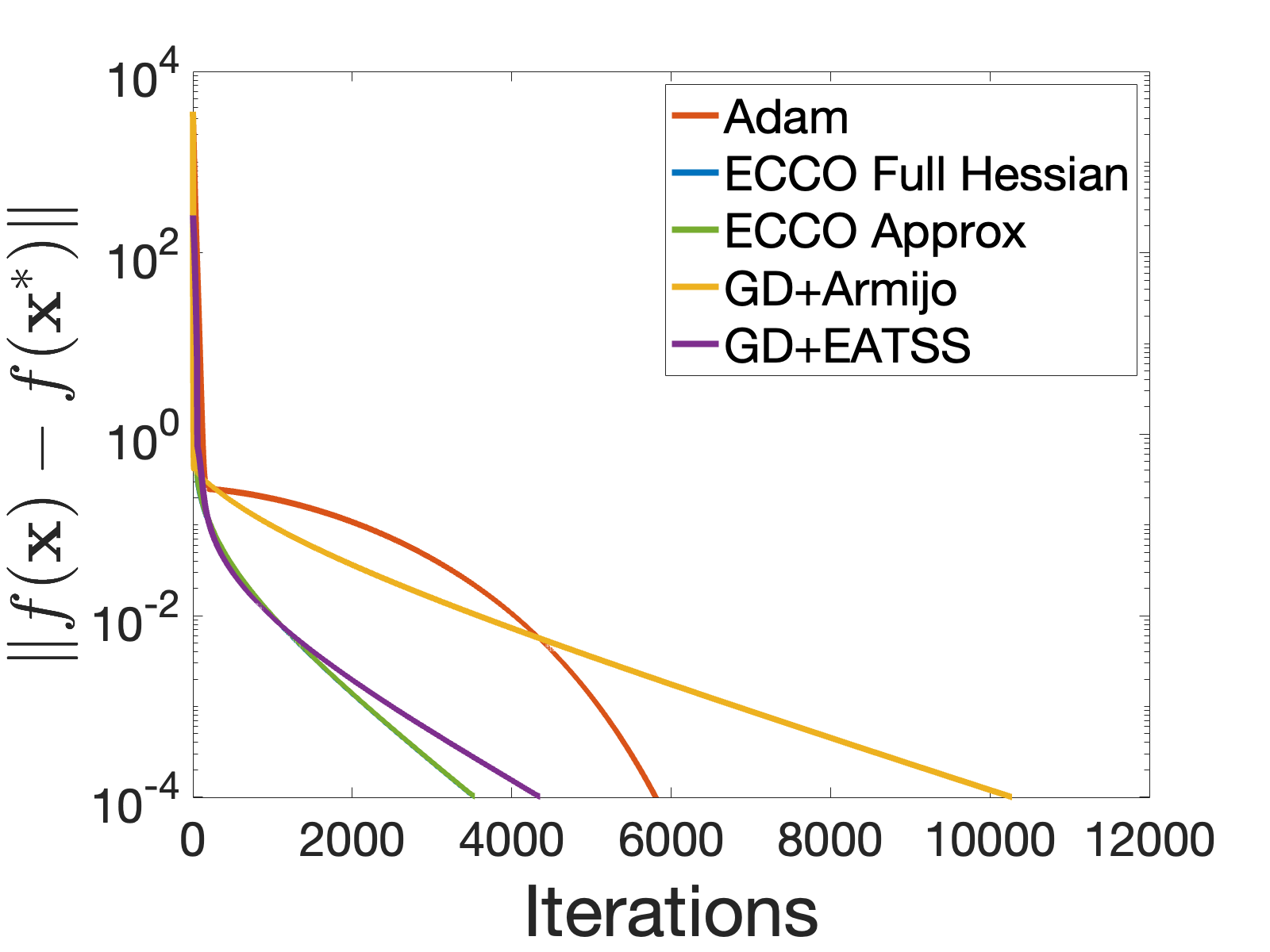}
    \caption{Rosenbrock $\x(0)=(-2,-2)$}
    \label{fig:a}
\end{subfigure}
\hfill
\begin{subfigure}{0.3\textwidth}
    \includegraphics[width=\textwidth]{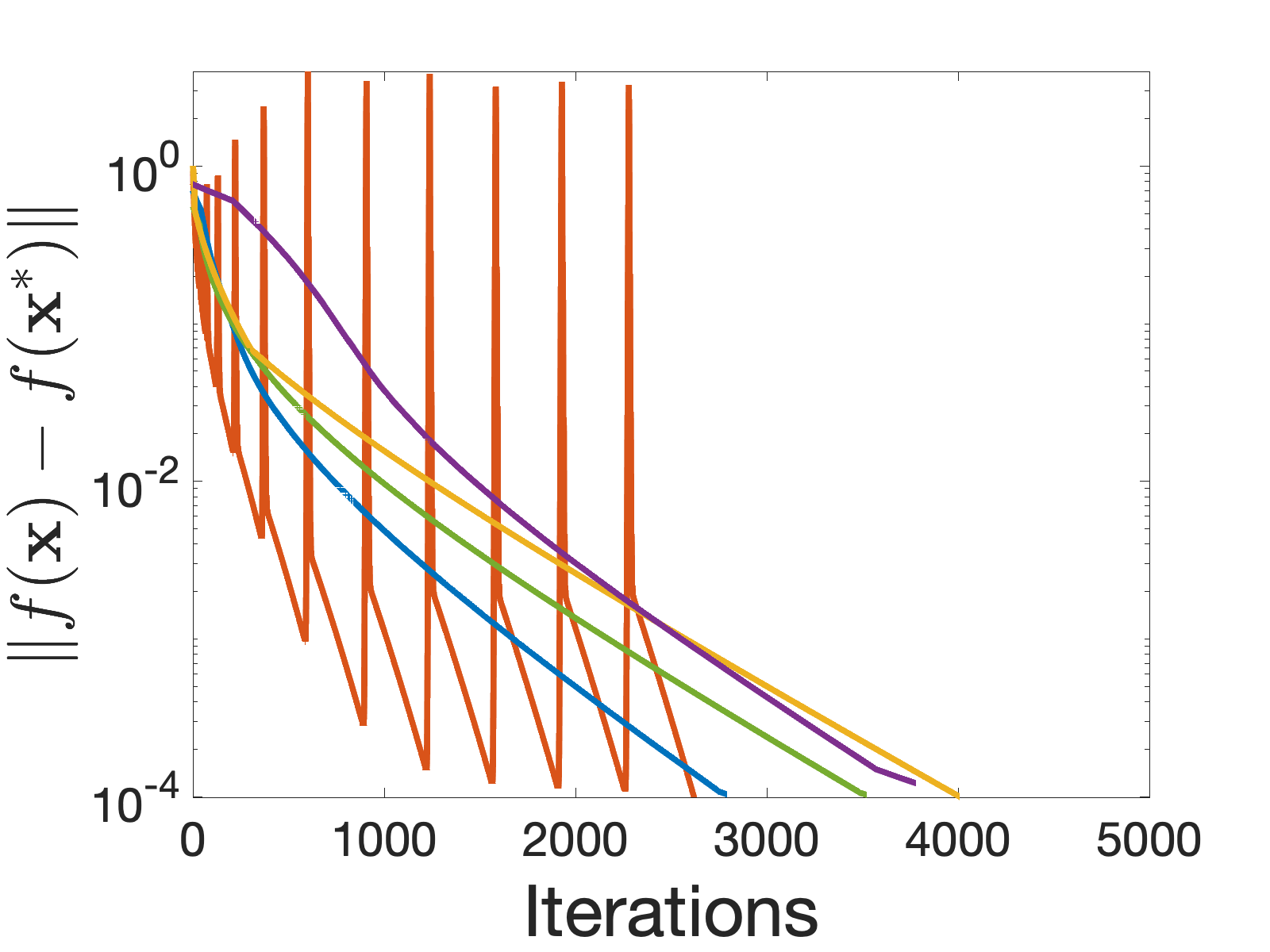}
    \caption{Rosenbrock $\x(0)=(0,0)$}
    \label{fig:b}
\end{subfigure}
\hfill
\begin{subfigure}{0.3\textwidth}
    \includegraphics[width=\textwidth]{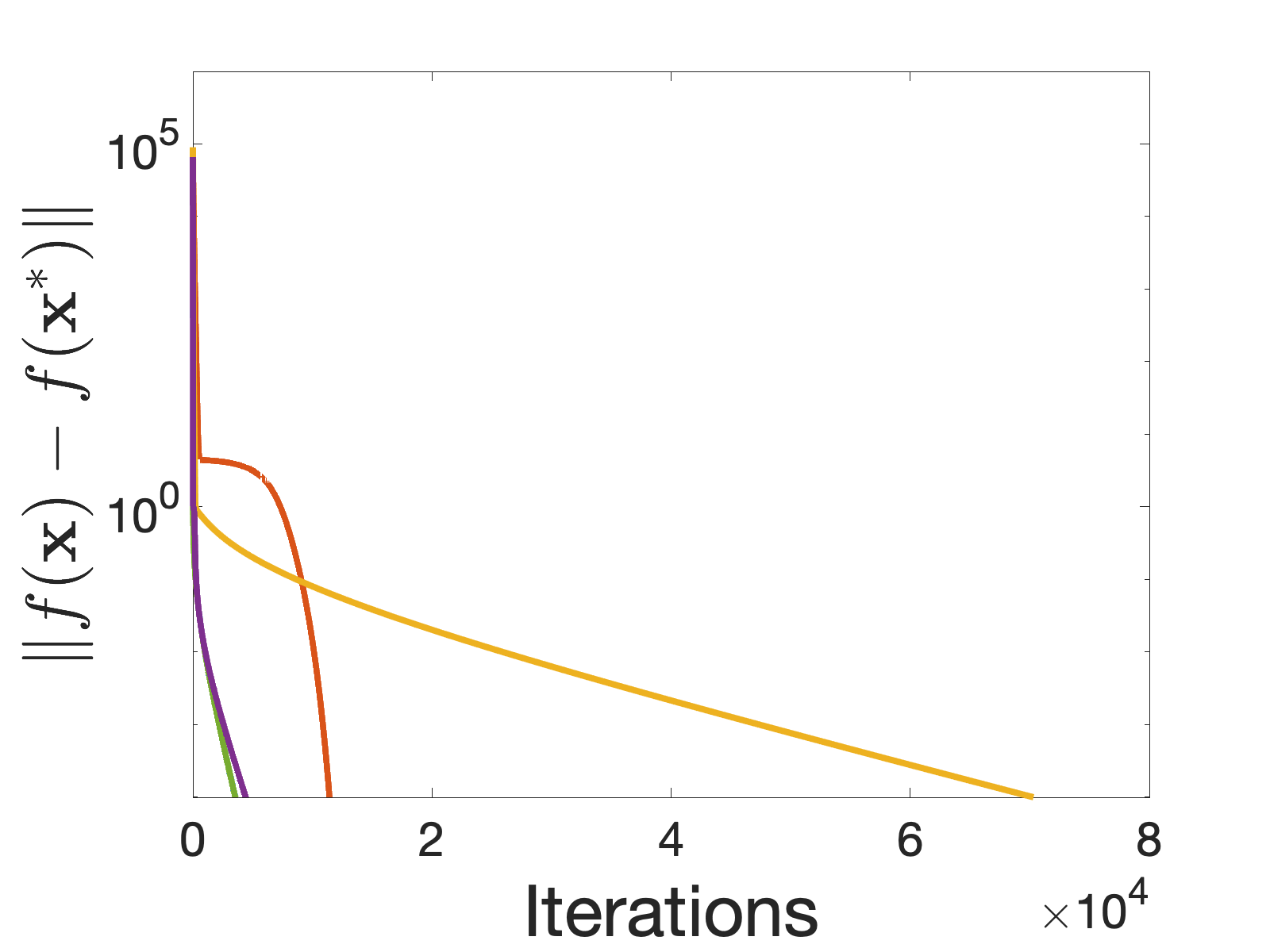}
    \caption{Rosenbrock $\x(0)=(-5,-5)$}
    \label{fig:c}
\end{subfigure}
\hfill
\begin{subfigure}{0.3\textwidth}
    \includegraphics[width=\textwidth]{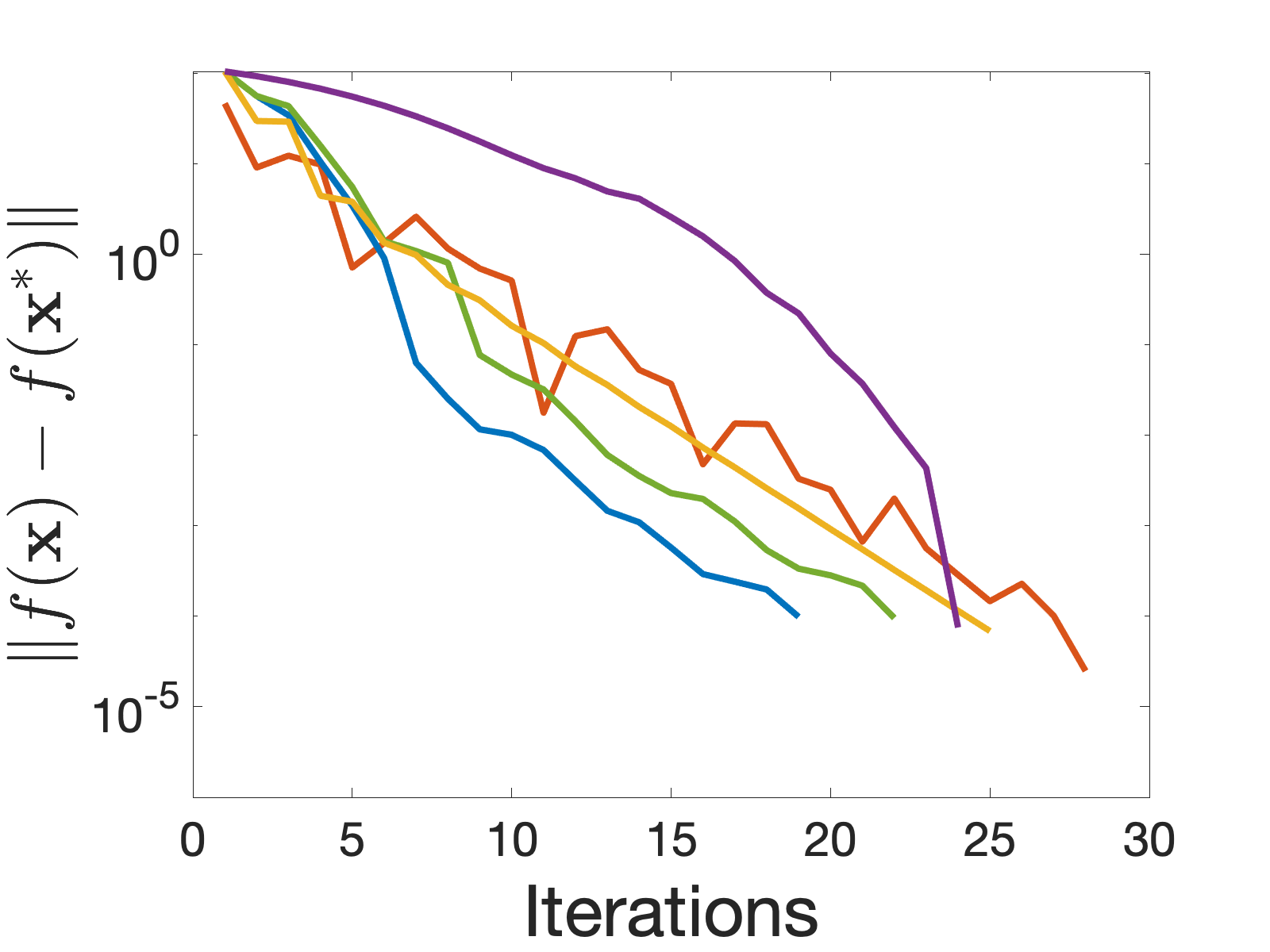}
    \caption{Himmelblau $\x(0)=(1,1)$}
    \label{fig:d}
\end{subfigure}
\hfill
\begin{subfigure}{0.3\textwidth}
    \includegraphics[width=\textwidth]{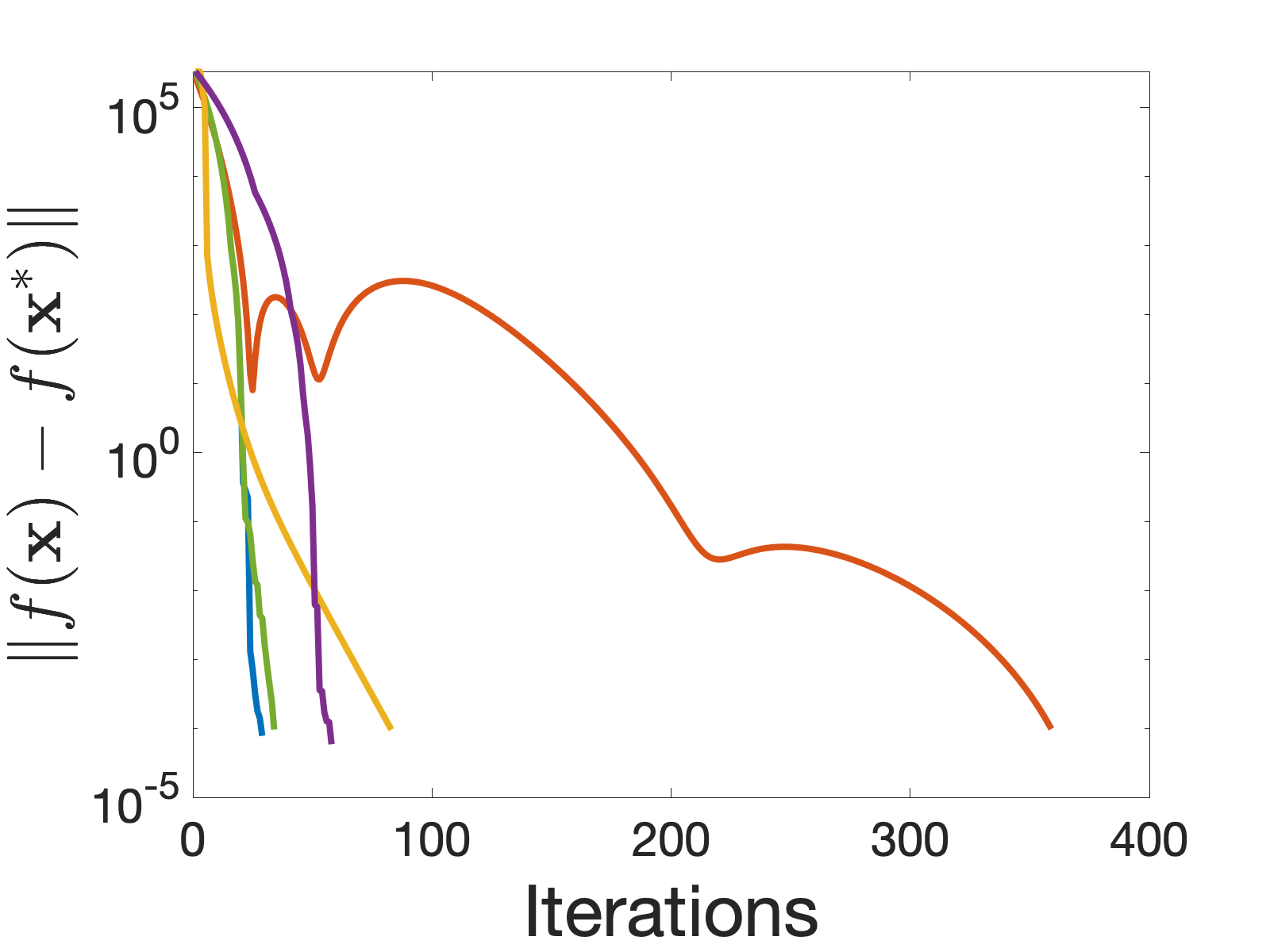}
    \caption{Himmelblau $\x(0)=(20,20)$}
    \label{fig:e}
\end{subfigure}
\hfill
\begin{subfigure}{0.3\textwidth}
    \includegraphics[width=\textwidth]{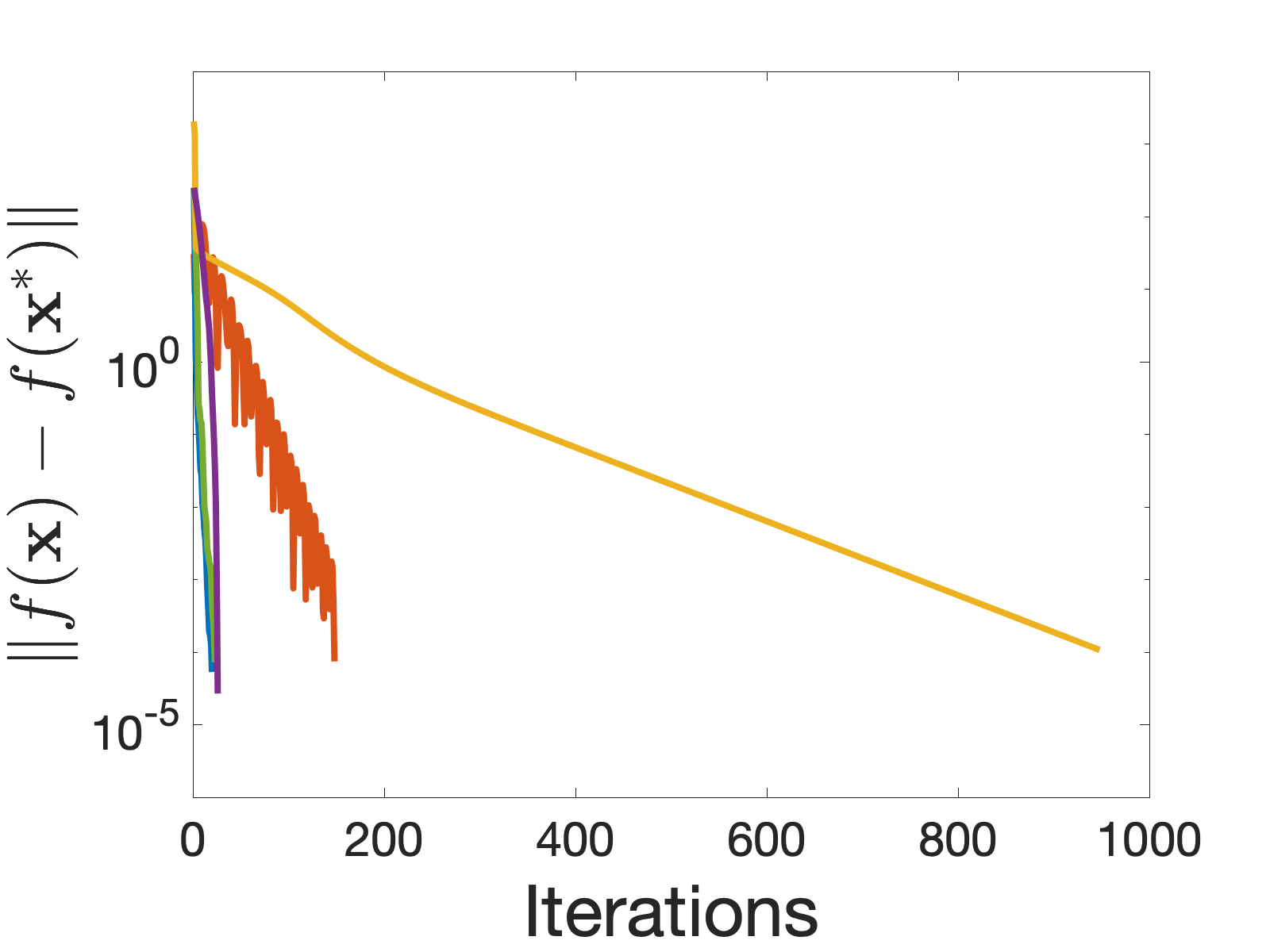}
    \caption{Himmelblau $\x(0)=(-5,-5)$}
    \label{fig:f}
\end{subfigure}
\hfill
\begin{subfigure}{0.3\textwidth}
    \includegraphics[width=\textwidth]{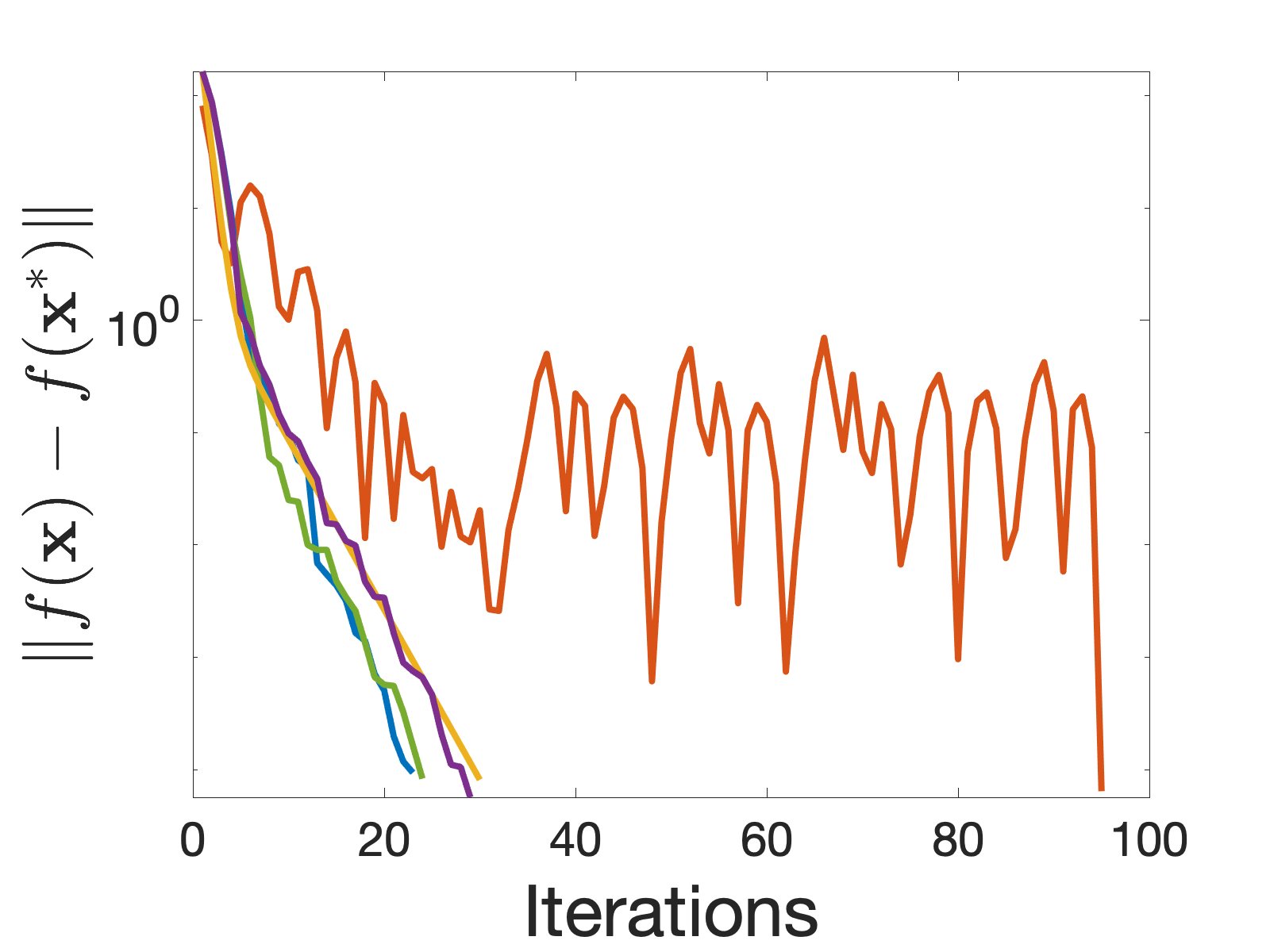}
    \caption{Booth $\x(0)=(5,5)$}
    \label{fig:g}
\end{subfigure}
\hfill
\begin{subfigure}{0.3\textwidth}
    \includegraphics[width=\textwidth]{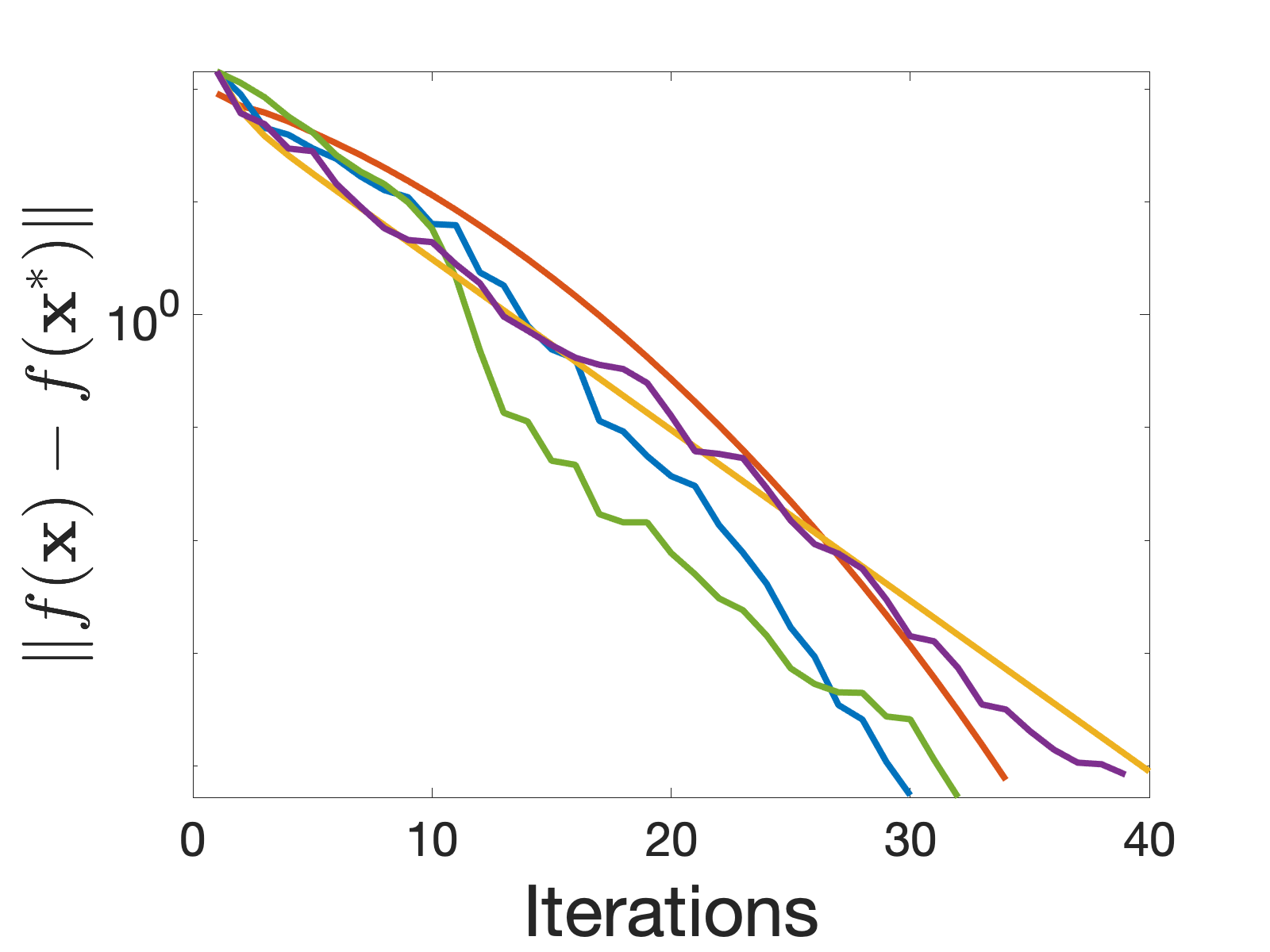}
    \caption{Booth $\x(0)=(5,-5)$}
    \label{fig:h}
\end{subfigure}
\hfill
\begin{subfigure}{0.3\textwidth}
    \includegraphics[width=\textwidth]{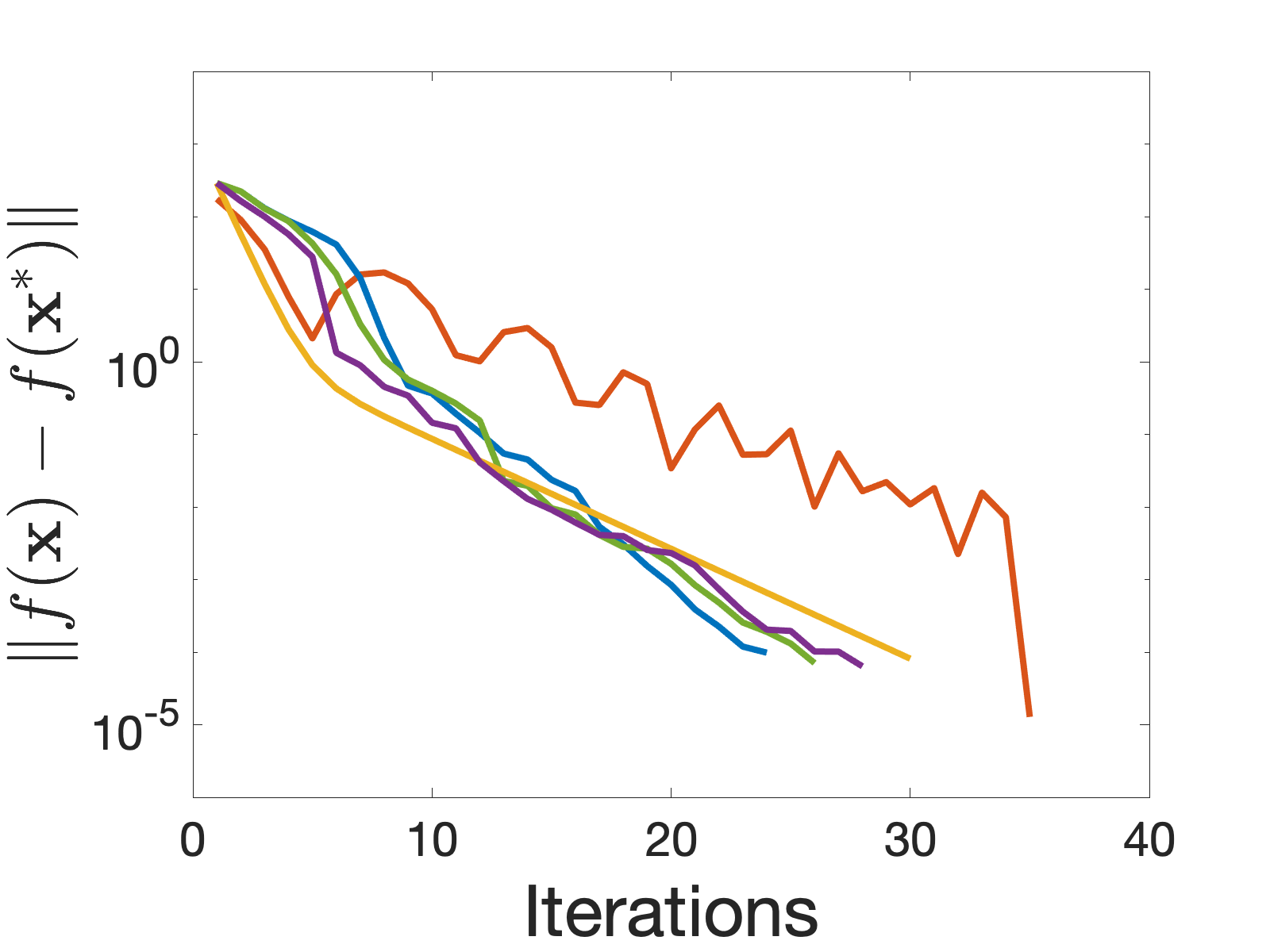}
    \caption{Booth $\x(0)=(-2,-2)$}
    \label{fig:i}
\end{subfigure}
\hfill
\begin{subfigure}{0.3\textwidth}
    \includegraphics[width=\textwidth]{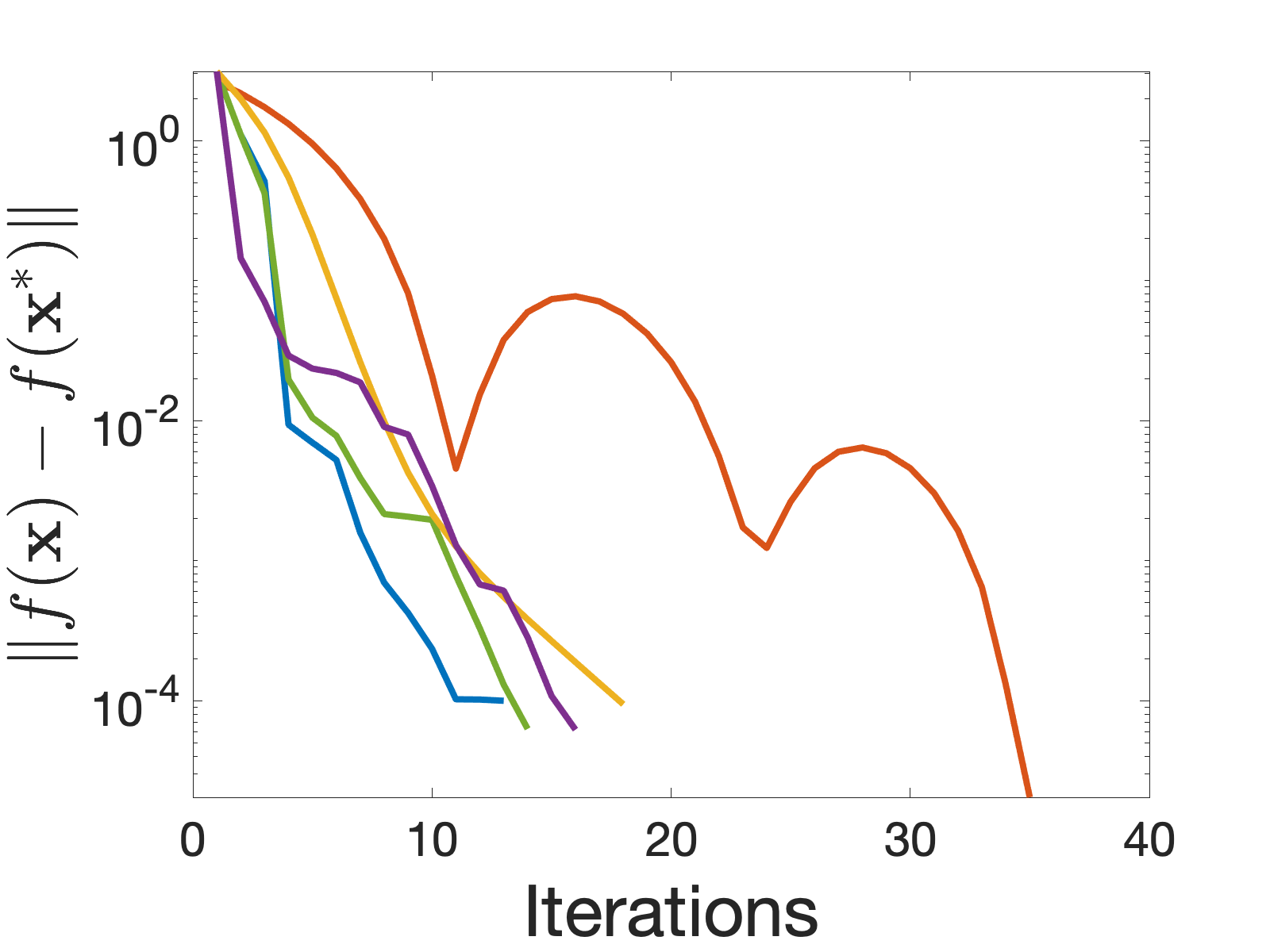}
    \caption{Three Hump $\x(0)=(1,1)$}
    \label{fig:j}
\end{subfigure}
\hfill
\begin{subfigure}{0.3\textwidth}
    \includegraphics[width=\textwidth]{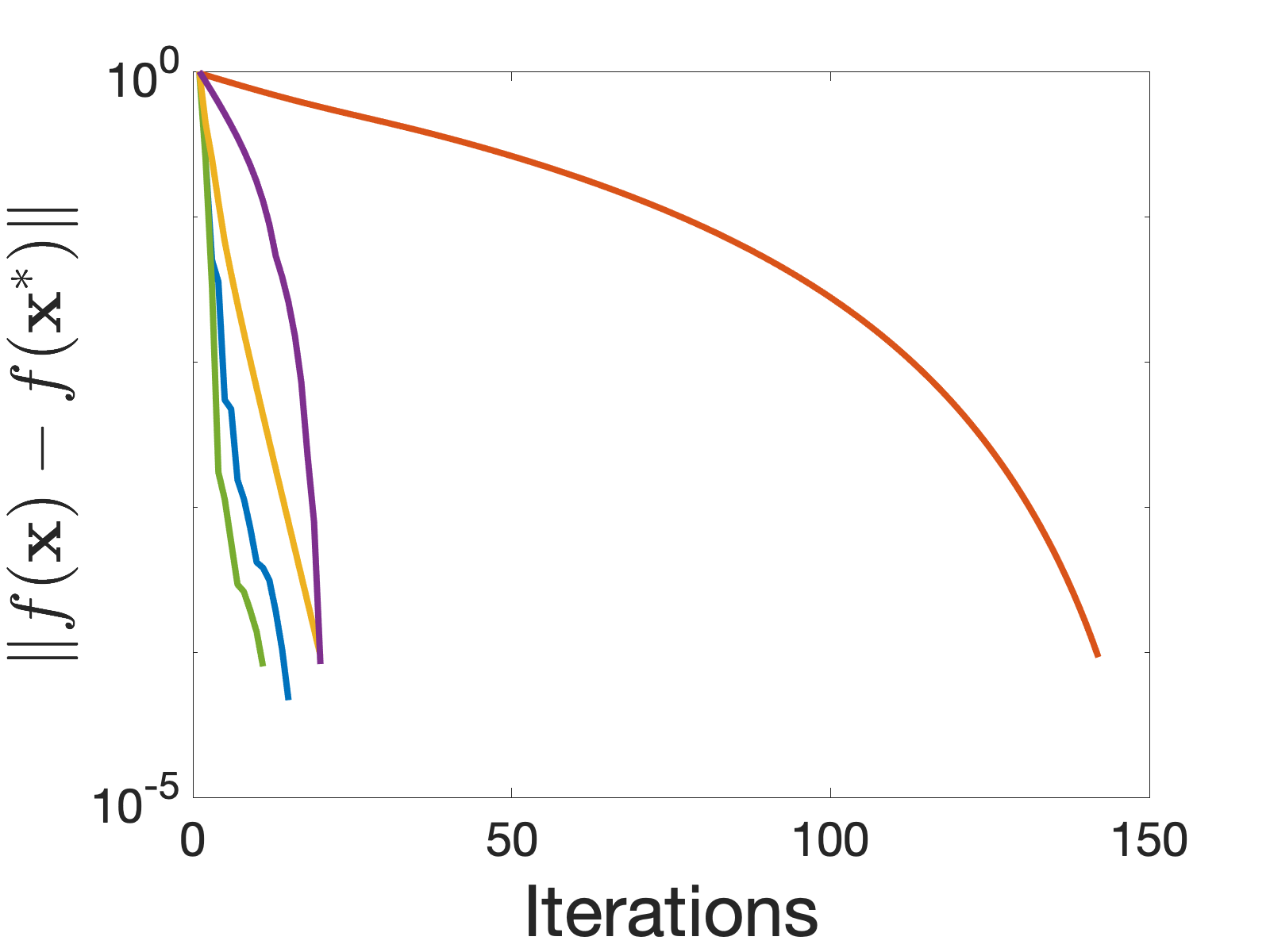}
    \caption{Three Hump $\x(0)=(0,-1)$}
    \label{fig:k}
\end{subfigure}
\hfill
\begin{subfigure}{0.3\textwidth}
    \includegraphics[width=\textwidth]{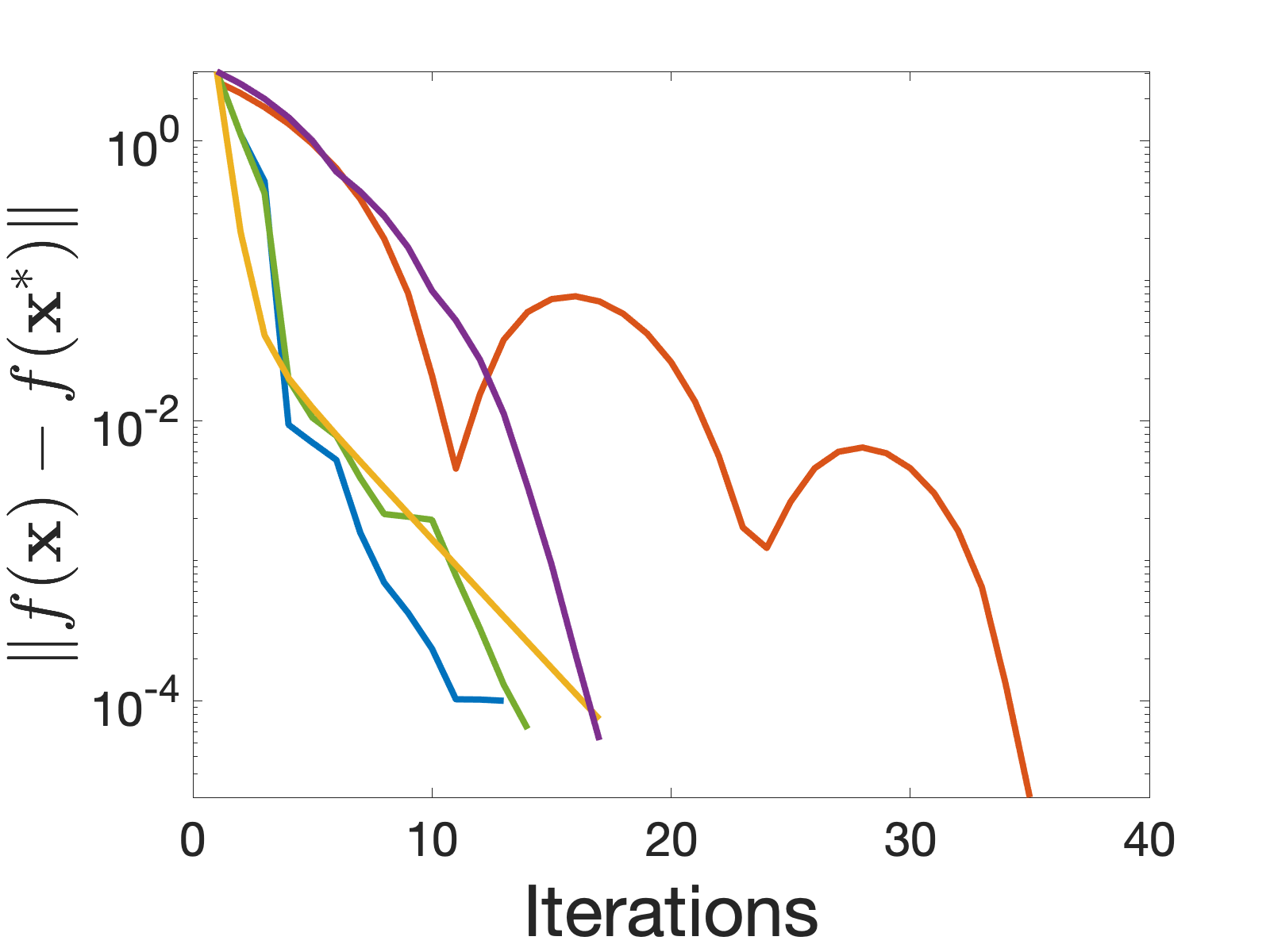}
    \caption{Three Hump $\x(0)=(-1,-1)$}
    \label{fig:l}
\end{subfigure}
\hfill
\begin{subfigure}{0.3\textwidth}
    \includegraphics[width=\textwidth]{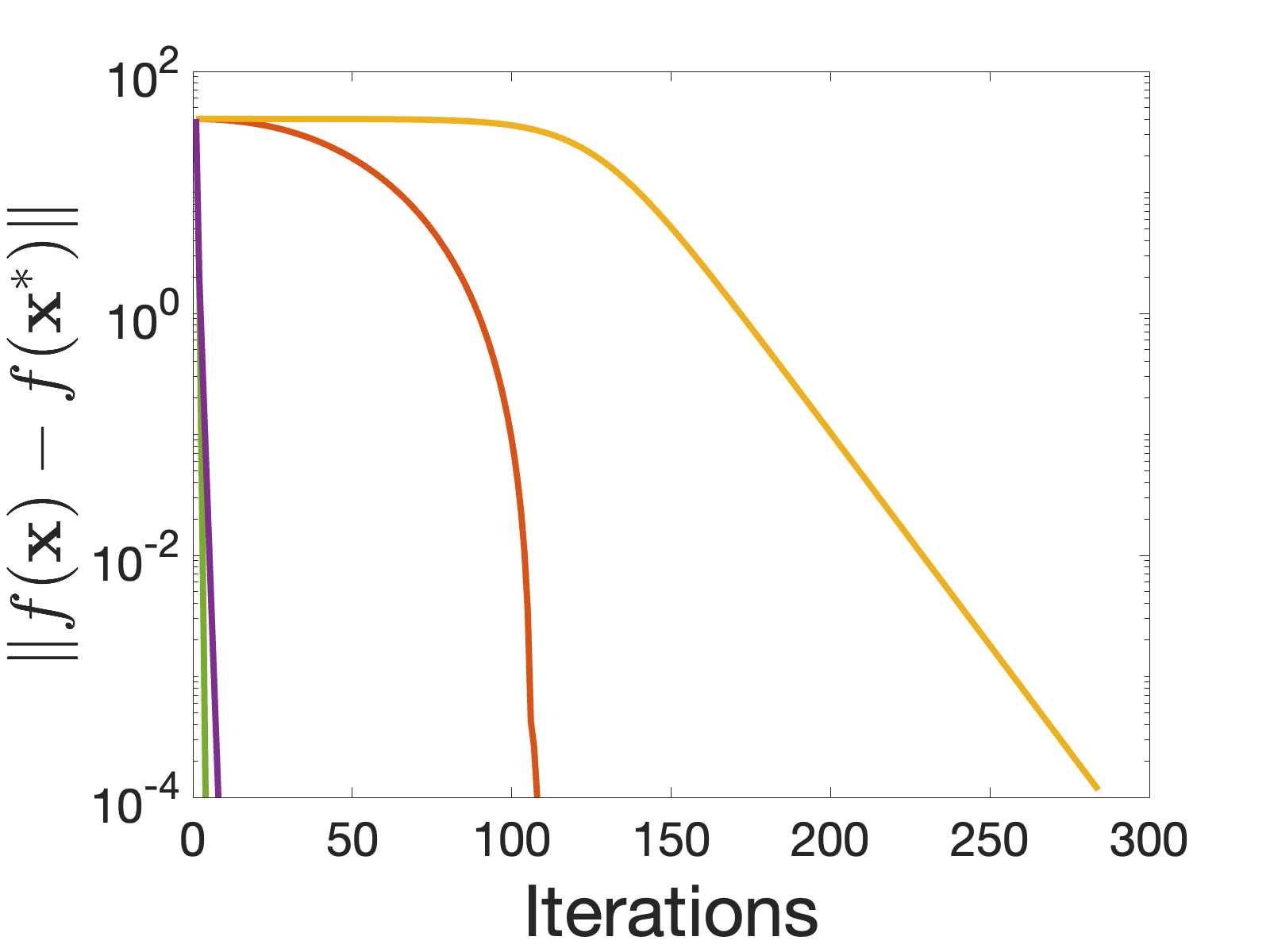}
    \caption{Rastrigin $\x(0)=(0.5, 0.5)$}
    \label{fig:m}
\end{subfigure}
\hfill
\begin{subfigure}{0.3\textwidth}
    \includegraphics[width=\textwidth]{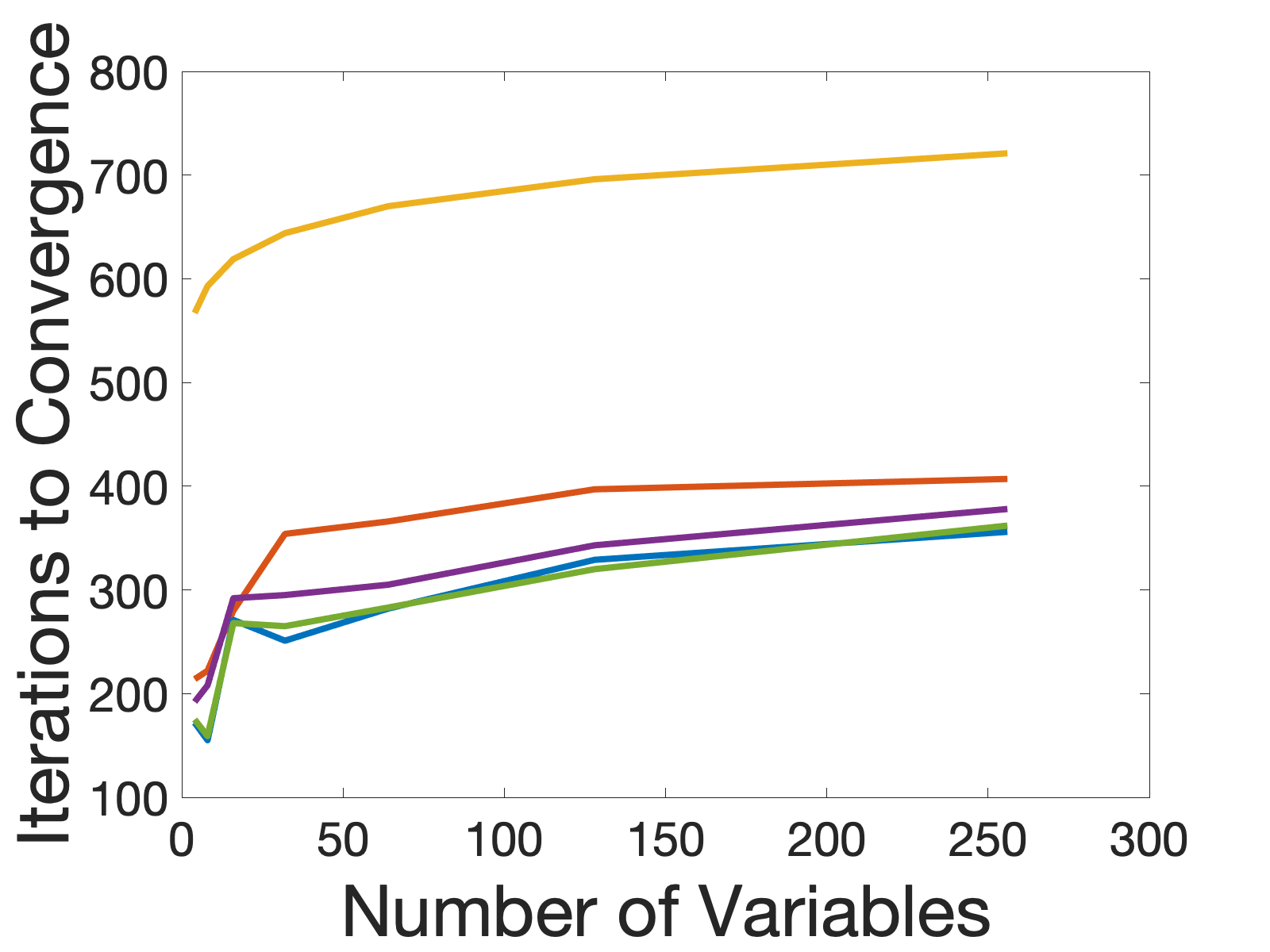}
    \caption{Extended Wood $\x(0)=2$}
    \label{fig:n}
\end{subfigure}
\hfill
\begin{subfigure}{0.3\textwidth}
    \includegraphics[width=\textwidth]{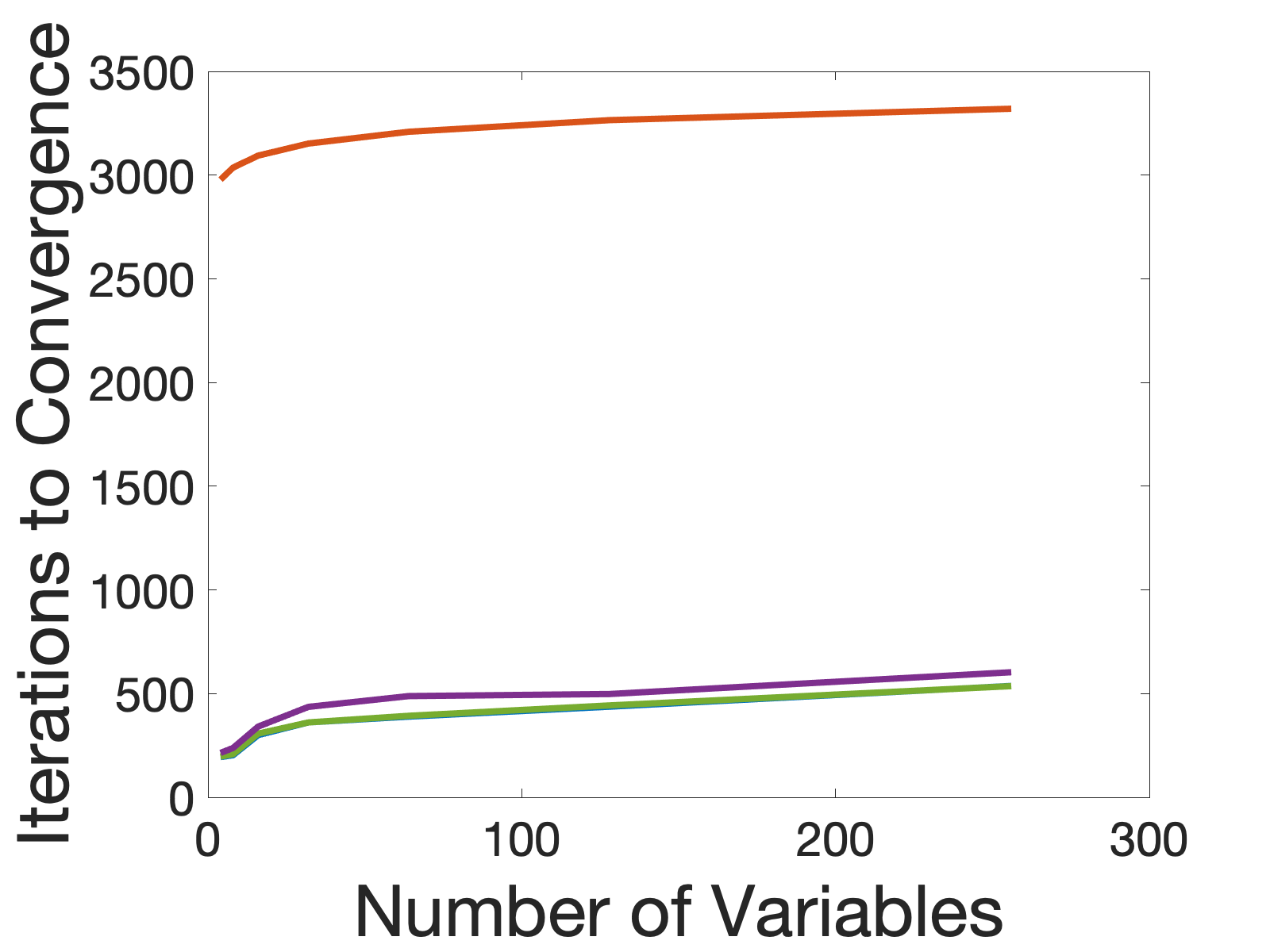}
    \caption{Extended Wood $\x(0)=10$}
    \label{fig:o}
\end{subfigure}

\caption{Test Functions}
\label{fig:test funcs}
\end{figure}

All the test function experiments can be found in Figure \ref{fig:test funcs}. The test functions considered were:
\begin{itemize}
    \item \textbf{Rosenbrock} \cite{testfunc} A convex, uni-modal, valley-shaped function.
    
    \item \textbf{Himmelblau} \cite{0070289212} A nonconvex, multi-modal, bowl-shaped function.

    \item \textbf{Booth} \cite{testfunc} A convex, uni-modal, plate-shaped function.

    \item \textbf{Three Hump} \cite{testfunc} A nonconvex, multi-modal, valley-shaped function.

    \item \textbf{Extended Wood} \cite{andrei2008unconstrained} A nonconvex function that scales for high dimensions.

    \item \textbf{Rastrigin} \cite{testfunc} A nonconvex, highly multi-modal function.
\end{itemize}

The compared methods were:
\begin{itemize}
    \item \textbf{Full Hessian ECCO} \eqref{z true}

    \item \textbf{Approximate ECCO} \eqref{z approx}

    \item \textbf{Gradient Descent} with step sizes chosen by line search and the Armijo condition

    \item \textbf{Gradient Descent} with step sizes chosen by EATSS, i.e. uncontrolled gradient flow discretized with FE+EATSS.

    \item \textbf{Adam} \cite{kingma2014adam} equipped with the full gradient 

\end{itemize}
Each experiment was terminated when convergence in the objective function was detected, i.e. when the algorithms returned $\x_k$ such that $\Vert f(\x_k)-f(\x_{k+1})\Vert <1e-4$. 

In all experiments, ECCO used $\delta=1$, $\alpha = 0.9$, $\beta=1.1$, $\eta=0.1$. The comparison methods used a grid search to optimize the hyperparameters. For Adam, we searched for $\beta_1$ and $\beta_2$ within $[0.7, 1]$ in increments of $0.01$. For both Adam and GD, we searched for an optimal learning rate within $[0.001,\dots, 1]$ in increments of $0.005$. The Armijo parameter was within $\{1e-5, 1e-4, 1e-3, 1e-2\}$.

These results empirically show that the approximation of the optimization trajectory does not induce much error into ECCO's performance, meaning that \eqref{z approx} is a reasonable substitution for \eqref{z true}. Second, FE+EATSS often outperforms FE with line search and the Armijo condition. Finally, Adam often performs badly on these test functions, which is expected as Adam excels as batch method in machine learning applications; we include it here for context.

\subsection{Wall Clock Times}
See Table \ref{toy wall clock}.
\begin{table}[t]
\caption{\textbf{Test Functions Experiment}. Wall clock time comparison of tested optimization methods on test functions experiment, where each row is normalized to the ECCO full Hessian implementation. Note that this measures the wall clock time of each measure to convergence.}
\label{toy wall clock}
\vskip 0.15in
\begin{center}
\begin{small}
\begin{sc}
\begin{tabular}{lccccccr}
\toprule
Method  & ECCO \eqref{z true} & ECCO \eqref{z approx} &  GD+FE+Armijo & GD+FE+EATSS & Adam\\
\midrule
 Fig \ref{fig:a} & 1 & \textbf{0.982} & 1.68 & 3.15 &1.72\\
 Fig \ref{fig:b} & 1 & 1.24 & 1.31 & 1.31 & \textbf{0.91}  \\
 Fig \ref{fig:c} & 1 & \textbf{0.99} & 7.26 & 1.09 & 1.96  \\
 Fig \ref{fig:d} & \textbf{1} & 1.18 & 1.22 & 1.19 & 1.24  \\
 Fig \ref{fig:e} & 1 & \textbf{0.98} & 1.45  & 1.63 & 3.27 \\
 Fig \ref{fig:f} & 1 & 1.01 & 17.4 & \textbf{0.97} & 3.96  \\
 Fig \ref{fig:g} & \textbf{1} & 1.015 & 1.023 & 1.055 & 3.013 \\
 Fig \ref{fig:h} & \textbf{1} & 1.16 & 1.21 & 1.32  & 1.13 \\
 Fig \ref{fig:i} & \textbf{1} & 1.08 & 1.12 & 1.19 & 1.28  \\
 Fig \ref{fig:j} & \textbf{1} & 1.03 & 1.15 & 1.14 & 2.52 \\
 Fig \ref{fig:k} & 1 & \textbf{0.93} & 1.13 & 1.13 & 6.06  \\
 Fig \ref{fig:l} & \textbf{1} & 1.01 & 1.17 & 1.17 & 2.03  \\
 Fig \ref{fig:m} & 1 & \textbf{0.99} & 15.2 & 1.04 & 6.04  \\
 Fig \ref{fig:n} ($n=256$) & 1 & \textbf{0.92} & 1.16 & 2.44 & 1.15  \\
 Fig \ref{fig:o} ($n=256$) & 1 & \textbf{0.91} & DNC & 1.16 & 5.07  \\
\bottomrule
\end{tabular}
\end{sc}
\end{small}
\end{center}
\vskip -0.1in
\end{table}

\subsection{Backtracking Line Search Details} \label{backtrack}
Across the test functions, the initial time step \eqref{LTE cond} led to a sequence of time steps that satisfied the LTE and monotonicity checks in Algorithm \ref{time step algo} on the first or second try in 78\% of the iterations. This is to be expected as \eqref{LTE cond} is known from circuit theory to often satisfy the LTE criterion \cite{rohrer1981passivity}. We can thus observe that the time step search subroutine did not, in general, lead to high per-iteration complexity. 


\section{More Details on Neural Network Experiment} \label{more nn}
The neural network to classify the MNIST dataset (\cite{deng2012mnist}) used to train a neural network was modified from \cite{nnet_code}. The first layer has 50 nodes with a ReLU activation, the second had 20 nodes and a sigmoid activation, and the output layer had 10 nodes and a softmax activation.

\subsection{Hyperparameter Selection} \label{hyper search}
The search for finding optimal hyperparameters was accomplished by discretizing the parameter space of Adam, GD and RMSProp and performing grid search. For Adam, we searched for $\beta_1$ and $\beta_2$ within $[0.7, 1]$ in increments of $0.01$. For Adam, GD and RMSProp, we searched for an optimal learning rate within $[0.001,\dots, 1]$ in increments of $0.005$. For RMSProp, we searched for an optimal decay rate within $[0.1,\dots,1]$ in increments of $0.005$.

\subsection{Wall Clock Times}\label{nn wall clock}
See Table \ref{nn wall clock}. Note that computing \eqref{z approx} has the same per-iteration complexity as a gradient descent update. The experimental results show that the EATSS routine in Algorithm \ref{time step algo} did not substantially affect the wall clock time, i.e. by orders of magnitude.

\begin{table}[t]
\caption{\textbf{Neural Network Experiment}. Wall clock time comparison of tested optimization methods on neural network experiment, normalized to the ECCO full Hessian implementation.  Note that this is based on 200 iterations of training.}
\label{nn wall clock}
\vskip 0.15in
\begin{center}
\begin{small}
\begin{sc}
\begin{tabular}{lcccr}
\toprule
Method  & Mean Normalized Time Per Iteration \\
\midrule
ECCO \eqref{z approx}  &  1 \\
SGD                    &  0.93 \\
Adam                   &  0.96 \\
\bottomrule
\end{tabular}
\end{sc}
\end{small}
\end{center}
\vskip -0.1in
\end{table}

\subsection{More Robustness Experiments} \label{more robust}
We tested the robustness of the neural network to the hyperparameters of the selected optimization methods by perturbing the optimal hyperparameter values within a normalized ball of some radius and recording the accuracy of the trained neural network. The hyperparameter values were sampled from a uniform distribution as $\tilde\theta \sim U(\max\{{\theta}^*-\frac{\varepsilon}{{\theta}^*}, \underline{\theta}\},\min\{{\theta}^*+\frac{\varepsilon}{{\theta}^*},\bar{\theta} \}) $, where ${\theta}^*$ was the optimal value as found by grid search. Note that this perturbation was bounded to be within the domain $[\underline{\theta} ,\bar{\theta}]$ of the hyperparameter. 

\begin{figure}
    \centering
    \includegraphics[width=0.6\linewidth]{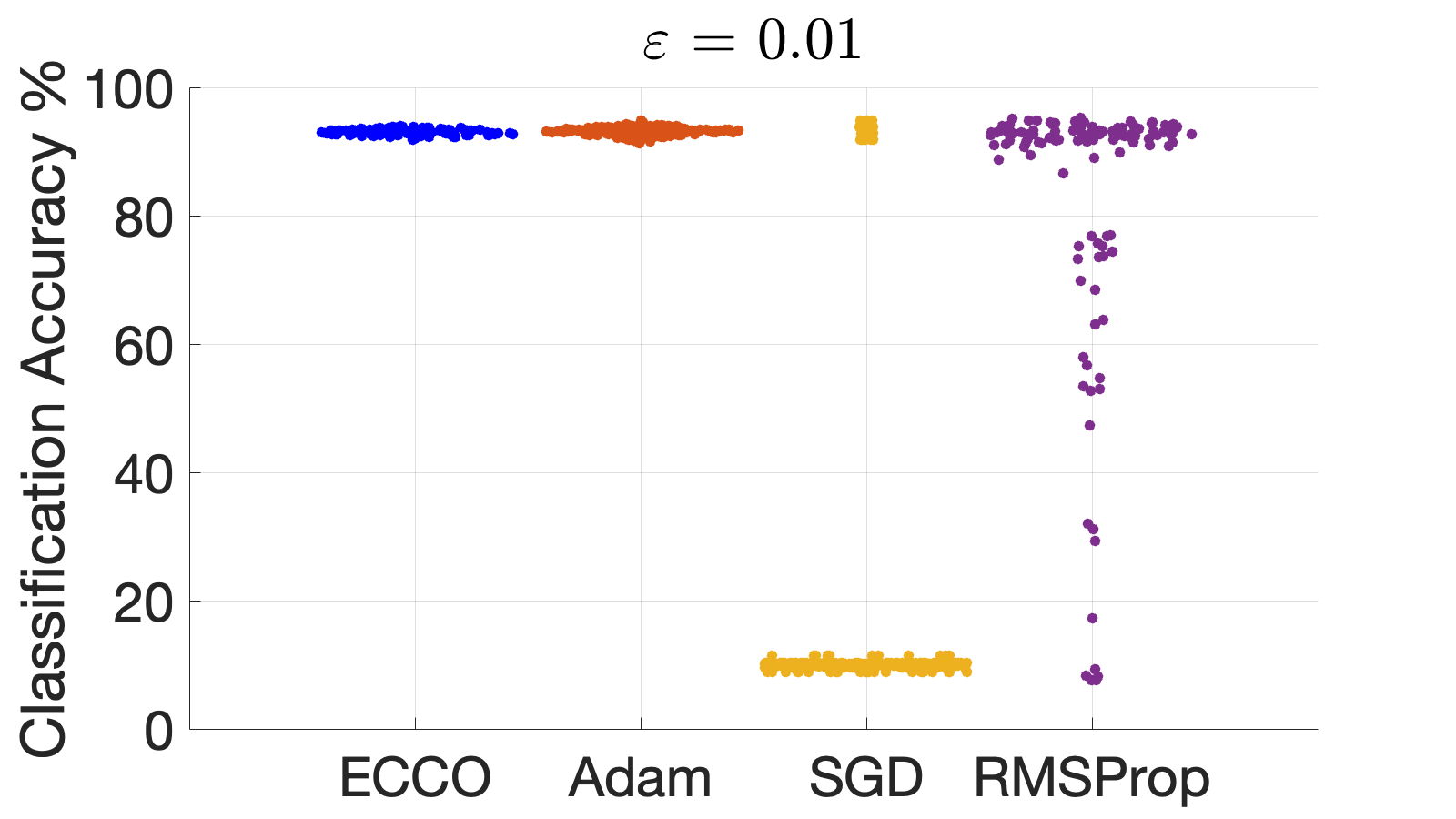}
    \caption{\small{Classification accuracy of 3-layer neural network using ECCO, Adam, gradient descent, and RMSProp with a random sampling of hyperparameters within a ball of $\varepsilon=0.01$ around the optimal hyper parameter values. Note that ECCO and Adam successfully train the neural network, while gradient descent and RMSProp are not able to consistently train in this regime.}}
    \label{fig:supp_class_acc_vareps_0.01}
\end{figure}

\begin{figure}
    \centering
    \includegraphics[width=0.6\linewidth]{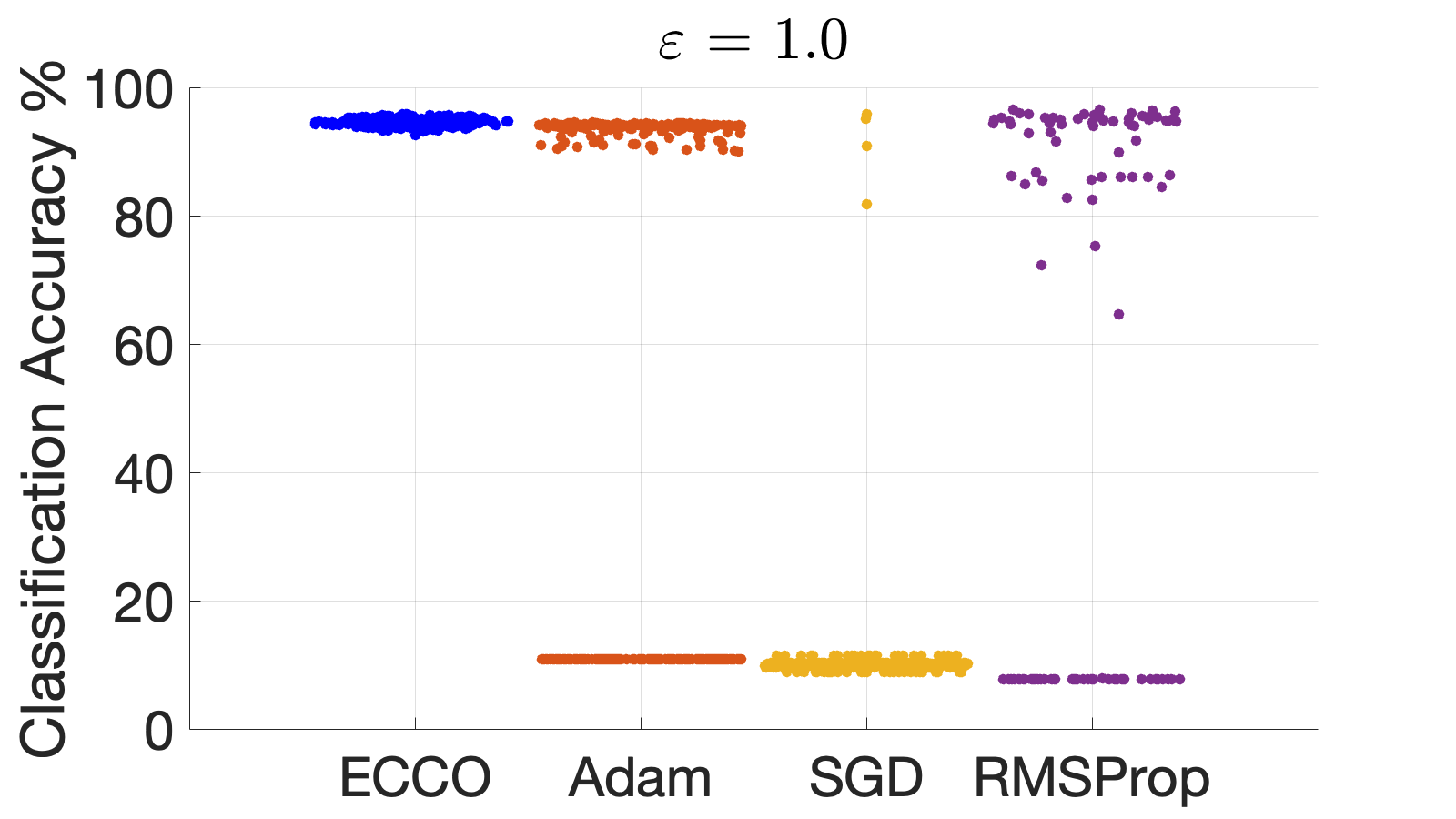}
    \caption{\small{Classification Accuracy of 3-layer neural network using ECCO, Adam, gradient descent, and RMSProp with a random sampling of hyperparameters within a ball of $\varepsilon=1$ around the optimal hyper parameter values. Note that ECCO is able to successfully train the neural network; however, none of the comparison methods are able to reliably train in this regime.}}
    \label{fig:supp_class_acc_vareps_1}
\end{figure}

In this section, we expand upon the results presented in the main paper by varying the radius for $\varepsilon\in\{0.01,1\}$ on ECCO and the comparison optimization solvers. For each fixed method and $\varepsilon$, 200 experiments were run, and the empirical classification accuracies are reported in Figures \ref{fig:supp_class_acc_vareps_0.01} and \ref{fig:supp_class_acc_vareps_1}. Note that the neural networks trained with ECCO are highly robust to the perturbations in the parameters; we found that when the neural network was trained with ECCO, we could perturb the optimal hyperparameters within a ball of $\varepsilon=1$ and lose at most 8\% in accuracy. 

ECCO's results are dramatically better than any of the comparison methods, where the accuracy suffered with small perturbations. Adam trained well for $\varepsilon=0.01$, but lost consistency and thus reliability for $\varepsilon=0.1$ and $\varepsilon=1$. Gradient descent and RMSProp similarly did not display robustness to the perturbations, implying that their methods need careful parameter tuning to be successful.

\begin{figure}
    \centering
    \includegraphics[width=0.6\linewidth]{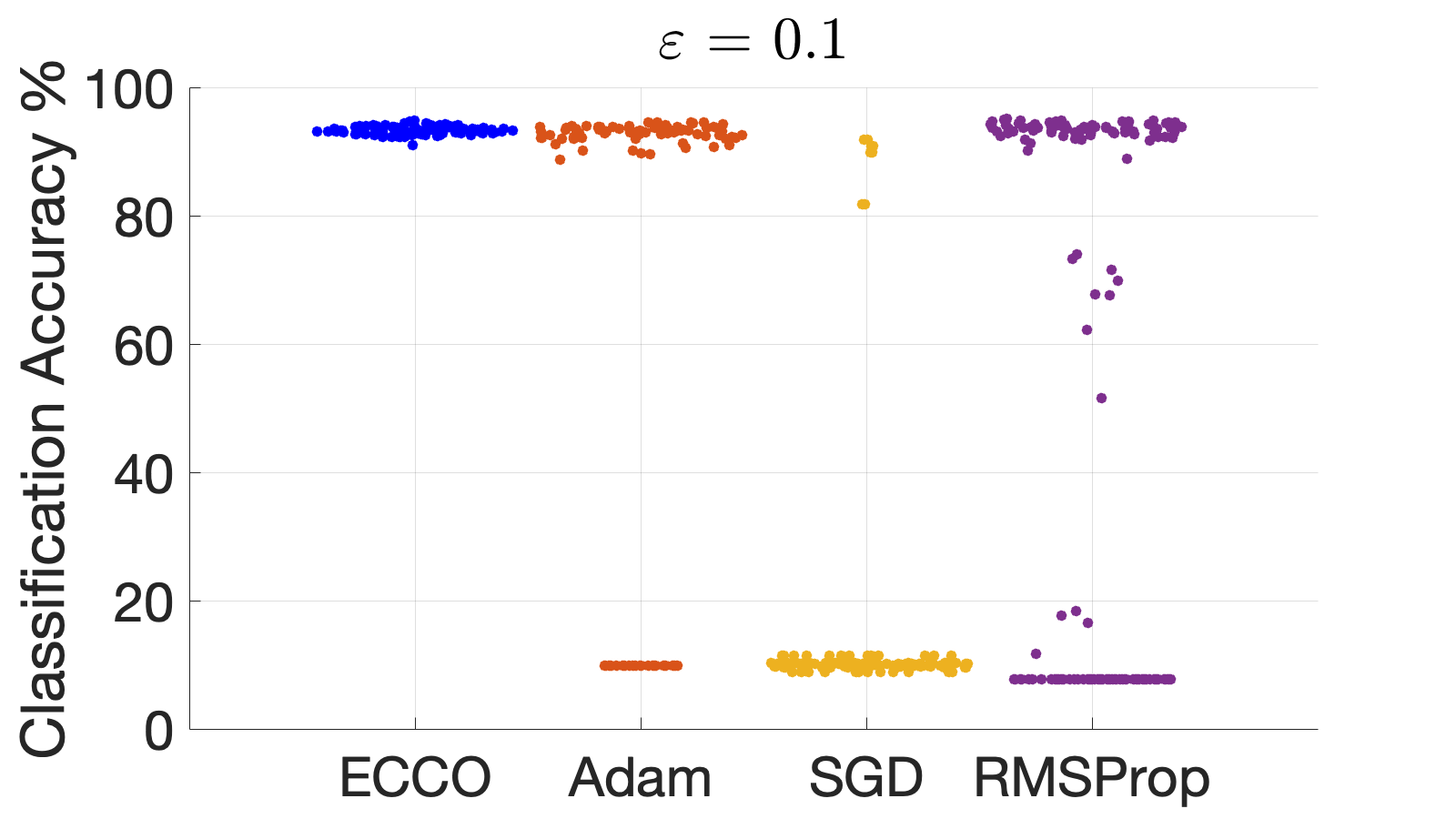}
    \caption{\small{Classification accuracy of 3-layer neural network using ECCO, Adam, gradient descent, and RMSProp with a random sampling of hyperparameters within a ball of $\varepsilon=0.1$ around the optimal hyper parameter values. Note that ECCO and Adam successfully train the neural network, while gradient descent and RMSProp are not able to consistently train in this regime.}}
    \label{fig:supp_class_acc_vareps_0.1}
\end{figure}

In addition, perturbations to the comparison methods easily caused divergence. The percentage of experiments for which the optimization methods successfully converged to a fixed point is reported in Figure \ref{diverged table}. Note that ECCO was always able to converge due to the choice of discretization, but perturbations to the optimally tuned hyperparameters for Adam, gradient descent, and RMSProp easily caused those methods to diverge. If divergence occurred, the neural network was not trained, and it yielded an accuracy of about 10\% which was equivalent to random guessing. This behavior is easily visible by observing Figures \ref{fig:supp_class_acc_vareps_0.01}, \ref{fig:supp_class_acc_vareps_0.1},  and \ref{fig:supp_class_acc_vareps_1}.

\begin{table}[t]
\caption{\textbf{Convergence of Experiments.} Percent of experiments with perturbed hyperparameters for which the optimization methods converged to an approximate fixed point. Note that ECCO was always able to converge, but perturbing the comparison methods easily caused divergence.}
\label{diverged table}
\vskip 0.15in
\begin{center}
\begin{small}
\begin{sc}
\begin{tabular}{lcccccr}
\toprule
 & ECCO Approx \eqref{z approx} & Adam & SGD & RMSProp \\
\midrule
$\varepsilon = 0.01$   & 100\% & 100\% & 61.01\% & 80.78\%  \\
$\varepsilon = 0.1$    & 100\% & 77.2\% & 26.10\% & 61.01\%  \\
$\varepsilon = 1$     & 100\% &  65.4\%& 2.51\% & 53.09\%  \\
\bottomrule
\end{tabular}
\end{sc}
\end{small}
\end{center}
\vskip -0.1in
\end{table}

These experiments suggest that ECCO may be used for neural network training for similar performance to state-of-the-art methods like Adam, but without the extensive need for hyperparameter tuning. This suggests ECCO would need less data than Adam for cross-validation procedures and is apt for generalization and distribution shift.

